\documentclass[11pt]{amsart}
\usepackage[left=2.3cm, top=2.3cm,bottom=2.3cm,right=2.3cm]{geometry}
\usepackage[T1]{fontenc}
\usepackage[utf8]{inputenc}
\usepackage{amsmath,amssymb,amsthm}
\usepackage{mathtools}
\usepackage{mathrsfs}
\usepackage{bm}
\usepackage{enumitem}
\usepackage[hidelinks]{hyperref}
\usepackage{graphicx}
\usepackage{subcaption}
\usepackage{newtxtext,newtxmath}
\usepackage{microtype}
\usepackage{booktabs,tabularx,array}
\usepackage{xcolor,fancyhdr}

\numberwithin{equation}{section}
\allowdisplaybreaks[2]

\newtheorem{theorem}{Theorem}[section]
\newtheorem{proposition}[theorem]{Proposition}
\newtheorem{lemma}[theorem]{Lemma}
\newtheorem{corollary}[theorem]{Corollary}
\theoremstyle{definition}
\newtheorem{definition}[theorem]{Definition}
\newtheorem{example}[theorem]{Example}
\theoremstyle{remark}
\newtheorem{remark}[theorem]{Remark}
\newcommand{\R}{\mathbb R}

\newcommand{\E}{\mathbb E}
\newcommand{\PP}{\mathbb P}
\newcommand{\eee}{\mathrm e}
\newcommand{\dd}{\mathop{}\!\mathrm d}

\newcommand{\normal}{\mathcal N}
\newcommand{\ind}{\mathbf 1}
\newcommand{\toD}{\xrightarrow{\mathrm d}}
\newcommand{\eqD}{\stackrel{\mathrm d}{=}}
\DeclareMathOperator{\Var}{Var}
\DeclareMathOperator{\Cov}{Cov}
\DeclareMathOperator{\tr}{tr}
\DeclareMathOperator{\relint}{relint}
\DeclareMathOperator{\spanop}{span}

\DeclareMathOperator{\Bern}{Bern}
\newcommand{\op}{\mathrm{op}}
\newcommand{\HS}{\mathrm{HS}}
\newcommand{\ip}[2]{\langle #1,#2\rangle}
\newcommand{\norm}[1]{\lVert #1\rVert}
\newcommand{\abs}[1]{\lvert #1\rvert}

\newcommand{\CCC}{\mathsf{CCC}}
\newcommand{\LS}{\mathsf{LS}}

\newcommand{\parhead}[1]{%
  \par\addvspace{0.5\baselineskip}%
  \noindent\textbf{#1}\enspace\ignorespaces}

\begin{document}

\title[CLT's for conic intrinsic volumes] {Central limit theorems for conic intrinsic volumes and Wiener chaos analysis of convex cones}
\author{Zakhar Kabluchko}
\address{
Zakhar Kabluchko:
Institut f\"ur Mathematische Stochastik,
Universit\"at M\"unster,
Orl\'eans-Ring 10,
48149 M\"unster, Germany}
\email{zakhar.kabluchko@uni-muenster.de}
\date{\today}

\begin{abstract}
Every nonempty closed convex cone $C\subseteq\R^d$ has
associated conic intrinsic volumes
$\upsilon_0(C),\ldots,\upsilon_d(C)$, which are nonnegative
and sum to one.
They therefore define an integer-valued random variable
$V_C$ with distribution $\PP\{V_C=k\}=\upsilon_k(C)$.
For each $n$, let $C_n\subseteq\R^{d_n}$ be a nonempty
closed convex cone and write $V_n=V_{C_n}$.
Assuming that $\Var(V_n)>0$ for all sufficiently large $n$,
our principal result is the equivalence
\[
 \frac{V_n-\E V_n}{\sqrt{\Var(V_n)}}
 \toD\normal(0,1)
 \quad\Longleftrightarrow\quad
 \Var(V_n)\longrightarrow\infty.
\]
A central limit theorem for $V_n$ was previously established
by Goldstein, Nourdin and Peccati
[in \emph{Gaussian Phase Transitions and Conic Intrinsic
Volumes: Steining the Steiner Formula},
Annals of Applied Probability, 2017],
under the stronger assumptions that $\E V_n\to\infty$ and
$\liminf_{n\to\infty}\Var(V_n)/\E V_n>0$.
In the present paper, we show that the divergence of
$\Var(V_n)$ alone is sufficient for the central limit
theorem and, moreover, that this condition is necessary.
Under the same condition, we also establish a uniform
local central limit theorem for the conic intrinsic volumes.

Let $\Pi_C$ denote the metric projection onto $C$ and
let $g$ be a standard Gaussian vector in $\R^d$.
The main idea of our proof is to identify the low-order
Wiener-chaos components of the function
$x\mapsto\norm{\Pi_C(x)}^2$ and to incorporate the
subtraction of its second Wiener-chaos component into
the master Steiner formula by an exponential change
of variables.

Our Wiener-chaos analysis also yields several geometric
results of independent interest.
Introducing the vector $m_C=\E\Pi_C(g)$ and the matrix
$A_C=\E[\Pi_C(g)g^{\mathsf T}]$, we show that
$x\mapsto A_Cx$ and $x\mapsto A_Cx+m_C$ are the optimal
linear and affine approximations, respectively, to
the metric projection $\Pi_C$ in Gaussian mean square.
We introduce the $L^2$-$\Pi$-metric on the space of
nonempty closed convex cones in $\R^d$ by setting
\[
 d_{\Pi,2}(C,D)
 =
 \left(\E\norm{\Pi_C(g)-\Pi_D(g)}^2\right)^{1/2}.
\]
We show that the space of nonempty closed
convex cones in $\R^d$, equipped with the metric $d_{\Pi,2}$,
admits an isometric embedding into a Hilbert-space sphere
of radius $\sqrt d/2$, with conic polarity corresponding to the
antipodal map. In particular, the diameter of this metric space is $\sqrt d$,
and the polar cone is the unique farthest cone from a given
cone. Furthermore, writing $\LS_d$ for the family of all linear
subspaces of $\R^d$, we establish the comparisons
\[
 \Var(V_C)
 \asymp \norm{m_C}^2
 \asymp \tr(A_C-A_C^2)
 \asymp d_{\Pi,2}(C,\LS_d)^2
 \asymp d_{\Pi,2}(C,-C)^2,
\]
where $\asymp$ denotes two-sided comparison up to positive
constants independent of the cone and the ambient dimension.
In particular, the intrinsic-volume random variables
are asymptotically normal after standardization if and
only if $d_{\Pi,2}(C_n,\LS_{d_n})\to\infty$ or,
equivalently, $d_{\Pi,2}(C_n,-C_n)\to\infty$.
\end{abstract}



\maketitle
\tableofcontents

\section{Introduction and main results}
\label{sec:introduction}

\subsection{Conic intrinsic volumes}

Throughout the paper, a closed convex cone in $\R^d$
is a nonempty closed set $C\subseteq\R^d$ such that
$\lambda x+\mu y\in C$ for all $x,y\in C$ and
$\lambda,\mu\geq0$.
In particular, every such cone contains the origin.

With every closed convex cone $C\subseteq\R^d$ one
associates a sequence of numbers
$\upsilon_0(C),\ldots,\upsilon_d(C)$, called its
\emph{conic intrinsic volumes}.
These quantities admit a particularly simple geometric
definition when $C$ is \emph{polyhedral}, that is,
when $C$ is the intersection of finitely many closed
half-spaces whose bounding hyperplanes pass through
the origin.
For $x\in\R^d$, let
\[
 \Pi_C(x)
 =
 \operatorname*{arg\,min}_{y\in C}\norm{x-y}
\]
denote the metric projection of $x$ onto $C$.
The minimizer exists and is unique because $C$ is
nonempty, closed, and convex.
Let $g$ be a standard Gaussian random vector in
$\R^d$, that is, a vector whose coordinates are
independent random variables with standard normal distribution
$\normal(0,1)$.
For a polyhedral cone $C\subseteq\R^d$, its conic
intrinsic volumes are defined by
\begin{equation}\label{eq:intrinsic-volumes-definition}
 \upsilon_k(C)
 =
 \PP\left\{
   \begin{array}{l}
    \Pi_C(g)\text{ belongs to the relative interior}\\
    \text{of a }k\text{-dimensional face of }C
   \end{array}
 \right\},
 \qquad k=0,\ldots,d;
\end{equation}
see \cite[Definition~2.1]{MT14}.
Here the relative interior of a convex set means
its interior within its affine span.
The relative interiors of the nonempty faces of
a polyhedral cone, including the cone itself,
form a disjoint partition of the cone.
Consequently, the numbers
$\upsilon_0(C),\ldots,\upsilon_d(C)$ are nonnegative
and satisfy
\begin{equation}\label{eq:intrinsic-volumes-normalization}
 \sum_{k=0}^d\upsilon_k(C)=1.
\end{equation}

For an arbitrary closed convex cone $C\subseteq\R^d$,
the intrinsic volumes are defined by approximation.
Let $\mathbb B^d$ denote the closed Euclidean unit
ball, and let $d_H$ denote the Hausdorff distance
between nonempty compact subsets of $\R^d$.
We use the \emph{conic Hausdorff metric} obtained
by taking the Hausdorff distance between the
intersections of cones with $\mathbb B^d$.
Choose polyhedral cones $C_m\subseteq\R^d$ such that
\[
 \lim_{m\to\infty}
 d_H(C_m\cap\mathbb B^d,C\cap\mathbb B^d)=0.
\]
Then define
\[
 \upsilon_k(C)
 =
 \lim_{m\to\infty}\upsilon_k(C_m),
 \qquad k=0,\ldots,d.
\]
The existence of such approximating sequences,
the existence of the limits, and their independence
of the chosen approximation are established in
\cite[Proposition~8.2]{MT14}; see also
\cite[Theorem~6.5.2(b)]{SchneiderWeil08}.
By passage to the limit, $\upsilon_k(C)\geq0$
for every $k=0,\ldots,d$, and
\eqref{eq:intrinsic-volumes-normalization} holds
for every closed convex cone $C$.
Throughout the paper, we use the convention that
$\upsilon_k(C)=0$ for all
$k\in\mathbb Z\setminus\{0,\ldots,d\}$.

\medskip
A systematic account of the geometry and probability
of convex cones, including the theory of conic
intrinsic volumes, is given in Schneider's monograph
\emph{Convex Cones: Geometry and Probability}
\cite{Schneider22}.
For further background, we refer to the monograph
of Schneider and Weil
\cite[Section~6.5]{SchneiderWeil08},
the dissertation of Amelunxen \cite{Amelunxen11},
and the papers of Amelunxen and Lotz
\cite{AmelunxenLotz17} and Amelunxen, Lotz, McCoy,
and Tropp \cite{ALMT14}.
Of particular relevance to the present paper is the
work of McCoy and Tropp \cite{MT14}, which develops
a probabilistic approach to conic intrinsic volumes
based on the master Steiner formula.
We shall recall the consequences of this formula
needed here in
Subsection~\ref{subsec:steiner-representation}.

\medskip
Conic intrinsic volumes arise in several areas of
mathematics, including constrained statistical
inference
\cite{BarlowEtAl72,RobertsonWrightDykstra88,SilvapulleSen05,Shapiro88},
random convex optimization and high-dimensional
inverse problems
\cite{ALMT14,MT14,AmelunxenBurgisser15},
and the combinatorics of hyperplane arrangements
\cite{KlivansSwartz11,AmelunxenLotz17}.
They also play an important role in geometric
probability and stochastic geometry, particularly
in the study of random conical tessellations
\cite{HugSchneider16},
convex hulls of random walks \cite{KVZ17a,KVZ17b},
and random polytopes, including Gaussian polytopes,
beta polytopes, and cells of tessellations arising from Poisson
point and hyperplane processes
\cite{KabluchkoZaporozhets20,Kabluchko20,Kabluchko21,
KabluchkoSteigenberger25}.

\subsection{The conic intrinsic-volume random variable}

For every closed convex cone $C\subseteq\R^d$,
the conic intrinsic volumes
$\upsilon_0(C),\ldots,\upsilon_d(C)$ are nonnegative
and sum to one; see
\eqref{eq:intrinsic-volumes-normalization}.
They therefore define an integer-valued random
variable $V_C$, called the
\emph{conic intrinsic-volume random variable}, by
\begin{equation}\label{eq:intrinsic-variable}
 \PP\{V_C=k\}=\upsilon_k(C),
 \qquad k=0,\ldots,d.
\end{equation}

Several fundamental properties of conic intrinsic
volumes admit simple probabilistic formulations.
For a closed convex cone $C\subseteq\R^d$, its
\emph{polar cone} is defined by
\[
 C^\circ
 =
 \{x\in\R^d:\ip{x}{y}\leq0
   \text{ for all }y\in C\}.
\]
If $C\subseteq\R^d$ and $D\subseteq\R^m$ are
closed convex cones, then
\begin{equation}\label{eq:polarity-product}
 V_{C^\circ}\eqD d-V_C,
 \qquad
 V_{C\times D}\eqD V_C+V_D,
\end{equation}
where $\eqD$ denotes equality in distribution.
In the second identity, the random variables
$V_C$ and $V_D$ on the right-hand side are taken
to be independent.

If $C$ is a linear subspace of dimension $k$,
then $V_C=k$ almost surely.
For every closed convex cone $C$ that is not
a linear subspace, the conic Gauss--Bonnet
relation states that
\begin{equation}\label{eq:gauss-bonnet}
 \PP\{V_C\text{ is even}\}
 =
 \PP\{V_C\text{ is odd}\}
 =
 \frac12.
\end{equation}
See, for example,
\cite{Schneider22,MT14,ALMT14}
for these identities.

A remarkable recent development is the proof
of the long-standing log-concavity conjecture
by Fu and Wang \cite[Theorem~1.1]{FuWang26}.
They established that, for every closed convex
cone $C\subseteq\R^d$,
\begin{equation}\label{eq:log-concavity}
 \upsilon_k(C)^2
 \geq
 \upsilon_{k-1}(C)\upsilon_{k+1}(C),
 \qquad 1\leq k\leq d-1.
\end{equation}

The expectation and variance of $V_C$ will be
denoted by
\begin{equation}\label{eq:mean-variance}
 \delta(C)=\E V_C,
 \qquad
 \tau^2(C)=\Var(V_C).
\end{equation}
The mean $\delta(C)$ is called the
\emph{statistical dimension} of $C$;
see \cite[Definition~2.2]{ALMT14}.
It extends the usual notion of dimension from
linear subspaces to closed convex cones:
if $L\subseteq\R^d$ is a linear subspace,
then $\delta(L)=\dim L$.
For a systematic study of the properties and
applications of the statistical dimension,
we refer to \cite[Sections~3--5]{ALMT14}.

By \eqref{eq:polarity-product}, we have
$\delta(C^\circ)=d-\delta(C)$ and
$\tau^2(C^\circ)=\tau^2(C)$.
Both the statistical dimension and the
intrinsic-volume variance are additive under
Cartesian products.

The intrinsic-volume random variable satisfies
useful concentration inequalities.
Amelunxen, Lotz, McCoy and Tropp
\cite[Theorem~6.1]{ALMT14}
established exponential concentration inequalities
for $V_C$ around its statistical dimension.
Subsequently, McCoy and Tropp
\cite[Theorem~4.5]{MT14} proved the variance bound
\begin{equation}\label{eq:variance-bound}
 \tau^2(C)
 \leq
 2\min\{\delta(C),\delta(C^\circ)\}
 =
 2\min\{\delta(C),d-\delta(C)\}.
\end{equation}
They also obtained refined exponential
concentration inequalities for $V_C$ around
$\delta(C)$; see
\cite[Corollary~4.10]{MT14}.
Further concentration results, including
sub-Gaussian tail bounds, were recently obtained
by Fern\'andez-Unzueta, Melbourne and
Palafox-Castillo
\cite[Corollary~6.2]{FMP25}.

These results show that the conic intrinsic-volume
distribution is sharply concentrated around its
statistical dimension.
These concentration estimates do not, however,
determine the limiting distribution after
centering and normalization.
The question of asymptotic normality was studied
by Goldstein, Nourdin and Peccati \cite{GNP17},
who established a central limit theorem under
additional growth assumptions on the statistical
dimension and the intrinsic-volume variance.
The assumptions of that theorem are stronger
than necessary.
In the next subsection, we state our principal
result, which gives a necessary and sufficient
condition for asymptotic normality.

\subsection{Central and local limit theorems}
We now state our main results.
Throughout the following three statements, let
$(d_n)_{n\geq1}$ be a sequence of positive integers
and let $C_n\subseteq\R^{d_n}$, $n\geq1$, be a
sequence of closed convex cones.
Write
\[
  V_n=V_{C_n},\qquad
  \delta_n=\E V_n,\qquad
  \tau_n^2=\Var(V_n).
\]

\begin{theorem}[Central limit theorem]\label{thm:clt}
If $\tau_n^2\to\infty$, then
\begin{equation}\label{eq:main-clt}
  \frac{V_n-\delta_n}{\tau_n}
  \toD\normal(0,1).
\end{equation}
\end{theorem}

Goldstein, Nourdin, and Peccati \cite{GNP17}
previously established \eqref{eq:main-clt}
under the stronger assumptions
\begin{equation}\label{eq:goldstein-nourdin-peccati-condition}
  \delta_n\to\infty,
  \qquad
  \liminf_{n\to\infty}
  \frac{\tau_n^2}{\delta_n}>0.
\end{equation}
Their more general Berry--Esseen estimate
\cite[Theorem~5.1 and Remark~5.2]{GNP17}
also covers certain regimes in which
$\tau_n^2/\delta_n\to0$.
Theorem~\ref{thm:clt} shows that no comparison between
the variance and the statistical dimension is needed:
divergence of the variance alone is sufficient.

\begin{example}[CLT holds but~\eqref{eq:goldstein-nourdin-peccati-condition} fails]
Consider $C_n=\R^{n^2}\times\R_+^n$.
By the product formula for conic intrinsic volumes,
\[
  V_n\stackrel{\mathrm d}=
  n^2+\operatorname{Bin}(n,1/2).
\]
Consequently, $\delta_n=n^2+\frac n2$, $\tau_n^2=\frac n4$,
and hence $\tau_n^2\to\infty$, $\tau_n^2/\delta_n\to0$.
Nevertheless, the classical binomial central limit
theorem applies. More generally, taking the Cartesian
product of a cone with a linear subspace shifts its
intrinsic-volume random variable by a deterministic
integer, without changing its variance or its
standardized distribution. This illustrates why
a condition involving the ratio $\tau_n^2/\delta_n$
is not intrinsic to the question of asymptotic
normality.
\end{example}

\begin{example}[$\delta_n \to\infty$ alone does not imply CLT]
If $\tau_n^2\to\infty$, then
\eqref{eq:variance-bound} implies that
$\delta_n\to\infty$. On the other hand, divergence
of $\delta_n$ alone does not suffice for the CLT.
Indeed, for the halfspace
$C_n=\R^{n-1}\times\R_+$, we have
\[
  V_n\eqD n-1+\Bern(1/2),\qquad
  \delta_n=n-\frac12,\qquad
  \tau_n^2=\frac14.
\]
Consequently, $\delta_n\to\infty$, but
$(V_n-\delta_n)/\tau_n$ takes the values
$-1$ and $1$ with equal probabilities
and is therefore not asymptotically Gaussian.
\end{example}


\medskip
The variance condition in Theorem~\ref{thm:clt}
is not only sufficient but also necessary.
Indeed, the random variables $V_n$ are integer-valued,
so the lattice spacing after normalization is
$1/\tau_n$. Convergence to a continuous Gaussian
distribution is therefore impossible along a subsequence
on which $\tau_n$ remains bounded. The following
proposition gives a more precise description of
subsequential limits in the bounded-variance regime.

\begin{proposition}[Subsequential limits at bounded variance]
\label{prop:bounded-variance}
Suppose that $0<\tau_n^2\leq M$
for some finite constant $M$ and all $n$.
Then every subsequence contains a further subsequence
along which
\[
  \frac{V_n-\delta_n}{\tau_n}
  \xrightarrow{\mathrm d}
  \frac{W-\alpha}{\tau},
\]
where $W$ is an integer-valued random variable,
$\alpha\in[0,1]$, and
$\tau\in[1/2,\sqrt M]$.
Moreover, the distribution of $W$ satisfies $\PP\{W\text{ is even}\} =\PP\{W\text{ is odd}\}=\frac12$.
In particular, every subsequential weak limit of the
standardized variables is discrete and nondegenerate.
\end{proposition}

The proof will be given in Subsection~\ref{subsec:proof-bounded-variance}. As a byproduct, we show that $\tau^2(C) \geq 1/4$ for every cone $C$ that is not a linear subspace.
Theorem~\ref{thm:clt} and
Proposition~\ref{prop:bounded-variance} imply the following result.
\begin{corollary}
Suppose that $\tau_n>0$ eventually. Then
\begin{equation}\label{eq:clt-equivalence}
  \frac{V_n-\delta_n}{\tau_n}
  \toD\normal(0,1)
  \quad\Longleftrightarrow\quad
  \tau_n^2\longrightarrow\infty.
\end{equation}
\end{corollary}
Thus, divergence of the variance gives a necessary
and sufficient condition for asymptotic normality
of conic intrinsic-volume random variables.

\medskip
Our next result strengthens the central limit theorem
to a local limit theorem.
It follows by combining Theorem~\ref{thm:clt} with the recently established log-concavity of conic intrinsic volumes \cite{FuWang26} and a general local limit theorem of Bender \cite{Bender73}.  The proof is given in Subsection~\ref{subsec:proof-central-local-clt}.
\begin{theorem}[Local central limit theorem]
\label{thm:local-clt}
If $\tau_n^2\to\infty$, then
\begin{equation}\label{eq:local-clt}
  \sup_{k\in\mathbb Z}
  \left|
    \tau_n\upsilon_k(C_n)
    -
    \frac{1}{\sqrt{2\pi}}
    \exp\left(
      -\frac{(k-\delta_n)^2}{2\tau_n^2}
    \right)
  \right|
  \longrightarrow0.
\end{equation}
\end{theorem}

\subsection{The master Steiner formula and its consequences}
\label{subsec:steiner-representation}

Let $C\subseteq\R^d$ be a nonempty closed convex cone, and
recall that $\Pi_C(x)$ denotes the metric projection of
$x\in\R^d$ onto $C$.
Let $g$ be a standard Gaussian vector in $\R^d$, and write
$\delta=\E V_C$ and $\tau^2=\Var(V_C)$.
The function
$$
H(x)=\norm{\Pi_C(x)}^2, \qquad x\in \R^d,
$$
will be central to
our proof of the central limit theorem for $V_C$.
Its relevance comes from a representation that expresses
the distribution of $H(g)$ as a mixture of chi-square
distributions, with mixing probabilities given by the
conic intrinsic volumes of $C$.

Recall that the chi-square distribution with $j$ degrees
of freedom is the distribution of the sum of the squares
of $j$ independent standard Gaussian random variables.
For $j=0$, this is the point mass at zero.
The following representation is a consequence of the
master Steiner formula of McCoy and Tropp;
see \cite[Theorem~3.1 and Corollary~3.2]{MT14}.

\begin{proposition}[Chi-square representation of the squared projection]
\label{prop:chi-square-representation}
Let $Z_1,\ldots,Z_d$ be independent standard Gaussian
random variables, independent of $V_C$.
Then
\begin{equation}\label{eq:chi-square-representation}
 H(g)=\norm{\Pi_C(g)}^2
 \eqD \sum_{i=1}^{V_C} Z_i^2,
\end{equation}
where the sum is interpreted as zero when $V_C=0$.
Thus the distribution of $H(g)$ is the mixture of the
chi-square distributions with $j$ degrees of freedom,
$0\leq j\leq d$, with respective weights $\upsilon_j(C)$.
\end{proposition}

The right-hand side of \eqref{eq:chi-square-representation}
can be constructed by first sampling $V_C$ and then,
conditionally on $V_C=j$, sampling a chi-square random
variable with $j$ degrees of freedom.

For a polyhedral cone, this representation has a direct
geometric explanation.
On the event that $\Pi_C(g)$ belongs to the relative
interior of a fixed $j$-dimensional face, it equals the
orthogonal projection of $g$ onto the span of that face.
For $j\geq1$, this event restricts the direction of the
Gaussian component in that span and the component in its
orthogonal complement, but not the length of the former.
The independence of Gaussian length and direction
therefore gives the conditional chi-square distribution
with $j$ degrees of freedom.
The case $j=0$ gives the point mass at zero.
Summing over the faces and using
\eqref{eq:intrinsic-volumes-definition} gives
\eqref{eq:chi-square-representation}.
The representation for arbitrary closed convex cones
follows by approximation by polyhedral cones.

We record two consequences that will be used repeatedly.
Both are established in \cite{MT14}; their short proofs
also clarify the distinction between the fluctuations of
$H(g)$ and those of $V_C$.

\begin{corollary}[Mean and variance of the squared projection]
\label{cor:H-moments}
For every nonempty closed convex cone $C\subseteq\R^d$,
\begin{equation}\label{eq:H-moments}
 \E H(g)=\delta,
 \qquad
 \Var(H(g))=2\delta+\tau^2.
\end{equation}
\end{corollary}

\begin{proof}
Conditionally on $V_C$, the random sum on the right-hand
side of \eqref{eq:chi-square-representation} has expectation
$V_C$ and variance $2V_C$.
Taking expectations gives $\E H(g)=\E V_C=\delta$.
The formula for conditional variance gives
$\Var(H(g))=\E[2V_C]+\Var(V_C)=2\delta+\tau^2$.
\end{proof}

These identities are
\cite[Propositions~4.3 and~4.4]{MT14}.
In particular, the intrinsic-volume variance can be
recovered from the squared projection through
$\tau^2=\Var(H(g))-2\E H(g)$.
The two terms in the variance decomposition
\eqref{eq:H-moments} have different meanings:
$2\delta$ is the expected conditional variance of the
chi-square variable, whereas $\tau^2$ is the variance of
its random number of degrees of freedom.

The same representation yields an identity between
Laplace transforms; see \cite[Proposition~4.6]{MT14}.

\begin{corollary}[Laplace-transform identity]
\label{cor:steiner-laplace}
For every nonempty closed convex cone $C\subseteq\R^d$
and every $s\in\R$,
\begin{equation}\label{eq:steiner-laplace}
 \E\eee^{sV_C}
 =
 \E\exp\left(
   \frac{1-\eee^{-2s}}{2}\,H(g)
 \right).
\end{equation}
Both expectations are finite.
\end{corollary}

\begin{proof}
For every $z<1/2$, the exponential-moment formula for a
chi-square random variable and
\eqref{eq:chi-square-representation} give
\[
 \E\eee^{zH(g)}
 =
 \sum_{j=0}^d
   \upsilon_j(C)(1-2z)^{-j/2}
 =
 \E(1-2z)^{-V_C/2}.
\]
Taking $z=(1-\eee^{-2s})/2$, which is less than $1/2$
for every real $s$, proves
\eqref{eq:steiner-laplace}.
The finite sum also proves finiteness.
\end{proof}

The chi-square representation and its consequences
play a central role in the approach of Goldstein,
Nourdin, and Peccati \cite{GNP17} to the normal
approximation of conic intrinsic volumes.
They will also be among the main ingredients of the
present paper.
A starting point is the following central limit theorem
for the squared projection, proved in
\cite[Theorem~1.1(1)]{GNP17}.

\begin{theorem}[Goldstein--Nourdin--Peccati]
\label{thm:squared-projection-GNP}
Let $C_n\subseteq\R^{d_n}$ be nonempty closed convex cones,
and let $g_n$ be a standard Gaussian vector in $\R^{d_n}$.
Write $\delta_n=\E V_{C_n}$ and
$\tau_n^2=\Var(V_{C_n})$.
If $\delta_n\to\infty$, then
\begin{equation}\label{eq:squared-projection-GNP}
 \frac{\norm{\Pi_{C_n}(g_n)}^2-\delta_n}
      {\sqrt{2\delta_n+\tau_n^2}}
 \toD \normal(0,1).
\end{equation}
\end{theorem}

Goldstein, Nourdin, and Peccati prove a stronger,
quantitative approximation in total variation;
see \cite[Proposition~3.1]{GNP17}.
We give an alternative proof of
\eqref{eq:squared-projection-GNP} in
Section~\ref{sec:squared-projection-clt}, using the
Laplace-transform estimates developed in
Section~\ref{sec:gaussian}.

Using the relation between squared projections and
conic intrinsic volumes, Goldstein, Nourdin, and Peccati
also obtain the central limit theorem for $V_{C_n}$
under the additional sufficient condition
$\liminf_{n\to\infty}\tau_n^2/\delta_n>0$;
see \cite[Theorem~1.1(2)]{GNP17} and
\eqref{eq:goldstein-nourdin-peccati-condition}.
This condition ensures that the standard deviations
$\sqrt{2\delta_n+\tau_n^2}$ and $\tau_n$ are comparable.

Under the sole assumption $\tau_n^2\to\infty$ of
Theorem~\ref{thm:clt}, however, the ratio
$\tau_n^2/\delta_n$ may tend to zero.
In this case, the variance of the squared projection
is dominated by the conditional chi-square variance
$2\delta_n$.
The convergence in \eqref{eq:squared-projection-GNP}
then describes fluctuations on a scale much larger
than the standard deviation $\tau_n$ of $V_{C_n}$,
and does not by itself provide the precision needed
for the intrinsic-volume central limit theorem.
The difficulty in passing from the squared-projection
CLT to the intrinsic-volume CLT is discussed in more
detail in
Subsection~\ref{subsec:intrinsic-volume-clt-difficulty}.

An extreme illustration is provided by a linear
subspace $L$.
Here $V_L=\dim L$ is deterministic, whereas
$\norm{\Pi_L(g)}^2$ has the chi-square distribution
with $\dim L$ degrees of freedom.
Thus the squared projection can have arbitrarily
large fluctuations even when the intrinsic-volume
random variable has none.

To overcome this difficulty, we need a more detailed
description of the fluctuations of the squared Gaussian
projection.
The main idea of our approach is to identify the
individual components in the Wiener-chaos expansion
of the function $H$, remove the component of order two,
and incorporate this subtraction into the
Laplace-transform identity \eqref{eq:steiner-laplace}.

As we shall see in the next subsection, the analysis
of the Wiener-chaos components also leads to several
geometric results of independent interest, including
a metric interpretation of the intrinsic-volume variance.

\subsection{Wiener-chaos expansions and the intrinsic-volume variance}
\label{subsec:intro-cone-chaoses}

We now examine the fluctuations of the squared Gaussian
projection through its Wiener-chaos expansion.
The first components of this expansion have geometric
interpretations and lead to alternative formulations
of the variance condition in Theorem~\ref{thm:clt}.

Let $\gamma_d$ denote the standard Gaussian probability
measure on $\R^d$, and recall that $g$ is a random vector
with distribution $\gamma_d$.
The Hilbert space $L^2(\R^d,\gamma_d)$ consists of
measurable real-valued functions $f$ satisfying
$\E[f(g)^2]<\infty$, with functions identified when they
agree $\gamma_d$-almost everywhere.
Its inner product is
$\ip{f}{h}_{L^2}=\E[f(g)h(g)]$.

For every integer $q\geq0$, the \emph{Gaussian chaos} of
order $q$, also called the \emph{Wiener chaos} of order $q$,
is the linear subspace of $L^2(\R^d,\gamma_d)$
consisting of all polynomials of degree at most $q$
that are orthogonal to all polynomials of degree
less than $q$.
These subspaces are mutually orthogonal, and every
$f\in L^2(\R^d,\gamma_d)$ admits an expansion
$f=\sum_{q=0}^{\infty}f_q$, converging in
$L^2(\R^d,\gamma_d)$, where $f_q$ is the orthogonal
projection of $f$ onto the chaos of order $q$.
In particular, the first three components have the form
\[
 f_0(x)=c,\qquad
 f_1(x)=\ip{b}{x},\qquad
 f_2(x)=x^{\mathsf T}Bx-\tr B,
\]
where $c\in\R$, $b\in\R^d$, and $B$ is a symmetric
$d\times d$ matrix.
The centering in the last expression follows from
$\E[g^{\mathsf T}Bg]=\tr B$.
The properties of these expansions needed in this
paper will be recalled in Section~\ref{sec:hermite}.

Let $C\subseteq\R^d$ be a nonempty closed convex cone,
and continue to write
$H(x)=\norm{\Pi_C(x)}^2$,
$\delta=\E V_C$, and $\tau^2=\Var(V_C)$.
Since $0\leq H(x)\leq\norm{x}^2$, the function $H$
belongs to $L^2(\R^d,\gamma_d)$.
We denote its Gaussian-chaos components by $H_q$.

The first three components are determined by the mean
Gaussian projection and the average derivative of the
projection map.
We introduce these quantities through
\begin{equation}\label{eq:intro-mean-matrix}
 m_C=\E\Pi_C(g),
 \qquad
 A_C=\E[\Pi_C(g)g^{\mathsf T}]
     =\E[D\Pi_C(g)].
\end{equation}
Here $D\Pi_C(x)$ denotes the Jacobian matrix of the
metric projection.
Since $\Pi_C$ is $1$-Lipschitz, this derivative exists
almost everywhere by Rademacher's theorem.
We shall establish these and further regularity
properties in
Proposition~\ref{prop:projection-properties} below.
The equality of the two representations of $A_C$
follows from Gaussian integration by parts, as we
shall show in Proposition~\ref{prop:H-chaos}.

The same proposition identifies the first three
components of $H$ itself:
\begin{equation}\label{eq:intro-H-chaoses}
 H_0(x)=\delta,\qquad
 H_1(x)=2\ip{m_C}{x},\qquad
 H_2(x)=x^{\mathsf T}A_Cx-\delta.
\end{equation}
Moreover, $A_C$ is symmetric and satisfies
\[
 0\preceq A_C\preceq I_d,
 \qquad
 \tr A_C=\delta.
\]
Here $B\preceq D$ means that $D-B$ is positive
semidefinite.
In particular, all eigenvalues of $A_C$ belong to
$[0,1]$.
The vector and matrix in \eqref{eq:intro-mean-matrix}
also determine the beginning of the chaos expansion
of the projection map itself: when the expansion is
applied to each coordinate of $\Pi_C$, its constant
component is $m_C$ and its linear component is
$x\mapsto A_Cx$.

The vector $m_C$ detects the lack of central symmetry
about the origin.
Indeed, as we shall prove in
Corollary~\ref{cor:cone-orthogonal-decomposition} below,
\[
 m_C=0
 \quad\Longleftrightarrow\quad
 C=-C
 \quad\Longleftrightarrow\quad
 C\text{ is a linear subspace}.
\]
This equivalence 
was already observed by Goldstein, Nourdin, and Peccati
\cite[Remark~4.1(b)]{GNP17}.
Reflection through the origin gives
$m_{-C}=-m_C$ and $A_{-C}=A_C$; these identities are
special cases of the transformation rules that will
be established in
Proposition~\ref{prop:chaos-transformations}.

When $m_C\neq0$, its direction has a concrete geometric
meaning.
Proposition~\ref{prop:mean-projection-and-kernels}
below will show that $m_C$ belongs to both $C$ and $-C^\circ$.
Consequently,
\[
 \R_+m_C
 \subseteq C
 \subseteq
 \{x\in\R^d:\ip{x}{m_C}\geq0\}.
\]
Thus the cone contains the ray in the direction of
$m_C$ and is contained in the halfspace with inward
normal $m_C$ whose boundary passes through the origin.

There is also an analytic interpretation of this
direction.
As we shall prove in
Proposition~\ref{prop:projection-gaussian-perturbations},
the expected squared projection under a small
translation of the Gaussian vector satisfies
\[
 \E H(g+h)
 =
 \delta+2\ip{m_C}{h}
 +h^{\mathsf T}A_Ch
 +o(\norm{h}^2),
 \qquad h\to0.
\]
For every unit vector $u$, the derivative
$\left.\frac{\dd}{\dd t}\E H(g+tu)\right|_{t=0}$
therefore equals $2\ip{m_C}{u}$.
If $m_C\neq0$, this derivative is uniquely maximized
by $u=m_C/\norm{m_C}$.
Hence the mean projection specifies the direction
of steepest increase of the expected squared
projection norm under infinitesimal translations.

The same expansion identifies $2A_C$ as the Hessian
at the origin of the function $h\mapsto\E H(g+h)$.
The matrix $A_C$ also describes the response to
changes in the Gaussian covariance.
For a unit vector $u$, changing the covariance of
the centered Gaussian input from $I_d$ to
$I_d+tuu^{\mathsf T}$ changes the expected squared
projection norm at derivative
$\ip{A_Cu}{u}$ when $t=0$.
This assertion will also be proved in
Proposition~\ref{prop:projection-gaussian-perturbations}.

The trace identity $\tr A_C=\delta$ separates the
radial and anisotropic parts of the second chaos.
More precisely, we shall establish in
Proposition~\ref{prop:radial-quadratic-chaos} below
the orthogonal decomposition
\[
 H_2(x)
 =
 \frac{\delta}{d}\bigl(\norm{x}^2-d\bigr)
 +
 x^{\mathsf T}
 \left(A_C-\frac{\delta}{d}I_d\right)x.
\]
The first term depends only on the Gaussian radius.
The second is determined by the traceless part of
$A_C$ and vanishes exactly when $A_C=(\delta/d)I_d$.
Thus the traceless part measures anisotropy at the
level of the second Gaussian chaos.
This does not imply that a cone with scalar $A_C$
is rotationally invariant.
For example, every self-dual cone, including the
nonnegative orthant, has $A_C=I_d/2$, as follows
from Proposition~\ref{prop:chaos-transformations}
below.

The coefficients $m_C$ and $A_C$ also identify the
best linear and affine approximations to the
projection map in Gaussian mean square.
As we shall prove in
Proposition~\ref{prop:projection-least-squares},
the unique optimal linear approximation to $\Pi_C$ is
$x\mapsto A_Cx$, and the unique optimal affine
approximation to $\Pi_C$ is $x\mapsto A_Cx+m_C$.
Their squared errors are
\begin{equation}\label{eq:intro-projection-approximation}
 \begin{aligned}
 \min_{B\in\R^{d\times d}}
 \E\norm{\Pi_C(g)-Bg}^2
 &=
 \tr(A_C-A_C^2),\\
 \min_{\substack{B\in\R^{d\times d}\\v\in\R^d}}
 \E\norm{\Pi_C(g)-Bg-v}^2
 &=
 \tr(A_C-A_C^2)-\norm{m_C}^2.
 \end{aligned}
\end{equation}
We write $a=\tr(A_C-A_C^2)$ for the squared error
of the optimal linear approximation.
Thus $\norm{m_C}^2$ is exactly the improvement in
squared approximation error obtained by allowing
an affine map.

We next relate these geometric quantities to the
intrinsic-volume variance.
Recall from \eqref{eq:H-moments} that
$\tau^2=\Var(H(g))-2\E H(g)$.
As we shall prove in
Proposition~\ref{prop:intrinsic-variance-chaoses}
below, the chaos expansion gives the more precise
identity
\begin{equation}\label{eq:intro-tau-chaoses}
 \tau^2
 =
 \sum_{q=1}^{\infty}
 \left(1-\frac q2\right)\E[H_q(g)^2]
 =
 2\norm{m_C}^2
 -
 \frac12\sum_{q=3}^{\infty}
 (q-2)\E[H_q(g)^2].
\end{equation}
This is a signed representation of the
intrinsic-volume variance; the underlying orthogonal
expansion is that of $H$.
The first chaos contributes positively, the
coefficient of the second chaos is zero, and all
higher components contribute nonpositively.

The vanishing coefficient of the second chaos
explains its special role in our proof.
Although this component can account for a large
part of the variance of $H(g)$, it makes no
contribution to \eqref{eq:intro-tau-chaoses}.
For a linear subspace, all centered fluctuations
of $H(g)$ belong to the second chaos, while
$\tau^2=0$.

In Propositions~\ref{prop:first-chaos-bounds}
and~\ref{prop:removal} below, we shall establish
the universal comparisons
\begin{equation}\label{eq:intro-variance-comparisons}
 \frac43\norm{m_C}^2
 \leq\tau^2\leq2\norm{m_C}^2,
 \qquad
 \frac12\tau^2\leq a\leq\frac32\tau^2.
\end{equation}
The upper bound $\tau^2\leq2\norm{m_C}^2$ was
proved by Goldstein, Nourdin, and Peccati
\cite[Theorem~4.1]{GNP17}, who also established
the lower bound
\[
 \tau^2
 \geq
 \frac{\norm{m_C}^4}{\sqrt{d\delta/2}},
 \qquad \delta>0.
\]
The lower bound $\tau^2\geq\frac43\norm{m_C}^2$
in \eqref{eq:intro-variance-comparisons} improves
their estimate by providing a lower bound
proportional to $\norm{m_C}^2$ with a universal
positive constant.

Thus the intrinsic-volume variance is comparable,
with constants independent of the cone and the
dimension, both to the squared norm of the mean
projection and to the optimal linear-approximation
error.

Motivated by \eqref{eq:intro-tau-chaoses}, we remove
the constant and second-chaos components of $H$.
The resulting function is
\begin{equation}\label{eq:intro-F}
 F(x)
 =
 H(x)-x^{\mathsf T}A_Cx
 =
 H_1(x)+\sum_{q=3}^{\infty}H_q(x),
\end{equation}
where the series is understood in
$L^2(\R^d,\gamma_d)$.
In particular, $\E F(g)=0$, its second-chaos
component vanishes, and its first-chaos component
is retained.
As we shall show in Proposition~\ref{prop:removal},
\begin{equation}\label{eq:intro-F-variance}
 \Var(F(g))=\tau^2+2a,
 \qquad
 2\tau^2\leq\Var(F(g))\leq4\tau^2.
\end{equation}
The standard deviation of $F(g)$ is therefore
comparable to $\tau$, even when the term $2\delta$
dominates the variance of $H(g)$.
This is the variance control needed for our proof
of the intrinsic-volume central limit theorem.

The first chaos is deliberately retained.
Indeed, Proposition~\ref{prop:first-chaos-bounds}
below will show that
\[
 \Var(H_1(g))
 =
 4\norm{m_C}^2
 \geq
 \frac34\Var(F(g)).
\]
Thus at least three quarters of the variance
remaining after the subtraction in
\eqref{eq:intro-F} always belongs to the first
chaos.
For circular cones with a fixed half-angle strictly
between $0$ and $\pi/2$, this proportion tends to
one as the dimension tends to infinity; see
Proposition~\ref{prop:circular-first-chaos} below.
In these examples, removing the first chaos as
well leaves a remainder with bounded variance,
although the intrinsic-volume variance diverges,
as we shall show in
Corollary~\ref{cor:circular-higher-chaoses}.
Consequently, such an additional subtraction
would discard the fluctuations whose size is
measured by $\tau$.

The comparisons in
\eqref{eq:intro-variance-comparisons} give
geometric formulations of the condition in
Theorem~\ref{thm:clt}.
Together with the bounded-variance obstruction
in Proposition~\ref{prop:bounded-variance}, they
yield the following characterization.

\begin{corollary}[Geometric characterizations of asymptotic normality]
\label{cor:intro-chaos-clt}
Let $C_n\subseteq\R^{d_n}$ be nonempty closed
convex cones, and write $V_n=V_{C_n}$,
$\delta_n=\E V_n$, and $\tau_n^2=\Var(V_n)$.
Assume that $\tau_n^2>0$ for all sufficiently
large $n$.
Then the following conditions are equivalent:
\begin{enumerate}
 \item
 $(V_n-\delta_n)/\tau_n\toD\normal(0,1)$;
 \item
 $\tau_n^2\to\infty$;
 \item
 $\norm{m_{C_n}}\to\infty$;
 \item
 $\tr(A_{C_n}-A_{C_n}^2)\to\infty$.
\end{enumerate}
\end{corollary}

\begin{proof}
Theorem~\ref{thm:clt} gives the implication from
the second condition to the first.
Conversely, if the variances do not tend to
infinity, there is a subsequence along which
they remain bounded.
Proposition~\ref{prop:bounded-variance} rules out
a standard Gaussian limit along this subsequence.
The equivalence with the third and fourth
conditions follows from
\eqref{eq:intro-variance-comparisons}.
\end{proof}

By \eqref{eq:intro-projection-approximation},
the last condition means that the optimal
Gaussian mean-square error of approximating
$\Pi_{C_n}$ by linear maps tends to infinity.
The third condition expresses the same criterion
through the mean Gaussian projection.
In the next subsection, we introduce a distance
between cones that gives these characterizations,
and the lack of central symmetry measured by
$m_C$, a quantitative metric interpretation.

\subsection{The \texorpdfstring{$L^2$-$\Pi$-metric}{L2-Pi-metric} and the geometry of convex cones}
\label{subsec:intro-projection-distance}

The approximation results of the preceding subsection
suggest comparing two cones through their projection
maps.
We use Gaussian mean square for this comparison, just
as in the approximation by linear and affine maps.

For $d\geq1$, let $\CCC_d$ denote the family of
all nonempty closed convex cones in $\R^d$.
For $C,D\in\CCC_d$, define
\begin{equation}\label{eq:intro-gaussian-projection-distance}
 d_{\Pi,2}(C,D)
 =
 \left(
   \E\norm{\Pi_C(g)-\Pi_D(g)}^2
 \right)^{1/2},
\end{equation}
where $g$ is a standard Gaussian vector in $\R^d$.
The same Gaussian vector is used in both projections.
As we shall prove in
Proposition~\ref{prop:gaussian-projection-distance}
below, this defines a metric on $\CCC_d$,
which we call the \emph{$L^2$-$\Pi$-metric}.

The distance also has a direct interpretation in
terms of uniformly distributed directions.
If $U$ is uniform on the unit sphere in $\R^d$,
positive homogeneity of the projection maps gives
\begin{equation}\label{eq:intro-spherical-projection-distance}
 \frac1d\,d_{\Pi,2}(C,D)^2
 =
 \E\norm{\Pi_C(U)-\Pi_D(U)}^2.
\end{equation}
Thus, up to a normalization factor depending only
on the ambient dimension, the $L^2$-$\Pi$-metric
can equivalently be defined using the uniform
distribution on the unit sphere.
It measures the discrepancy between the metric
projection maps in the $L^2$-sense.

The geometry of this metric has a simple description
in a Hilbert space.
Consider the vector-valued space
$L^2(\R^d,\gamma_d;\R^d)$, whose squared norm is
$\norm{f}_{L^2}^2=\E\norm{f(g)}^2$.
As we shall show in
Proposition~\ref{prop:cone-spherical-embedding}
below, the assignment $C\mapsto J_C$, where
$J_C(x)=\Pi_C(x)-x/2$, is an isometric embedding
of $(\CCC_d,d_{\Pi,2})$ into the sphere
of radius $\sqrt d/2$ in this Hilbert space.
Moreover, $J_{C^\circ}=-J_C$, so polarity
corresponds to the antipodal map.

In particular, the distances satisfy
\begin{equation}\label{eq:intro-polar-distance}
 d_{\Pi,2}(C,D)^2+d_{\Pi,2}(C^\circ,D)^2=d,
 \qquad C,D\in\CCC_d.
\end{equation}
We shall deduce in
Corollary~\ref{cor:metric-characterization-polarity}
that
\[
 \operatorname{diam}(\CCC_d,d_{\Pi,2})=\sqrt d,
 \qquad
 d_{\Pi,2}(C,D)=\sqrt d
 \quad\Longleftrightarrow\quad
 D=C^\circ.
\]
Thus the polar cone is the unique farthest cone
from $C$.
In this sense, polarity can be recovered from the
metric structure alone.

The statistical dimension has an immediate metric
interpretation.
Recalling that $\delta=\E V_C$, we have
\begin{equation}\label{eq:intro-statistical-dimension-distance}
 d_{\Pi,2}(C,\{0\})^2=\delta,
 \qquad
 d_{\Pi,2}(C,\R^d)^2=d-\delta.
\end{equation}
These identities will be recorded in
Corollary~\ref{cor:statistical-dimension-distance}
below.
We next give a corresponding interpretation of
the intrinsic-volume variance
$\tau^2=\Var(V_C)$.

Let $\LS_d$ be the family of all linear
subspaces of $\R^d$, including $\{0\}$ and $\R^d$,
and write
$d_{\Pi,2}(C,\LS_d)
=\inf_{L\in\LS_d}d_{\Pi,2}(C,L)$.
Recall from the preceding subsection the matrix
$A_C=\E[\Pi_C(g)g^{\mathsf T}]$ and the quantity
$a=\tr(A_C-A_C^2)$, which is the squared error
of the best linear approximation to $\Pi_C$
in Gaussian mean square.
For a linear subspace $L$, the approximation
identity \eqref{eq:intro-projection-approximation}
gives
\begin{equation}\label{eq:intro-cone-subspace-distance}
 d_{\Pi,2}(C,L)^2
 =
 a+\norm{A_C-P_L}_{\HS}^2.
\end{equation}
Here $P_L$ is the orthogonal projection matrix
onto $L$, and
$\norm{B}_{\HS}^2=\sum_{i,j=1}^d B_{ij}^2$
is the squared Hilbert--Schmidt norm.

Consequently, finding a nearest linear subspace
amounts to approximating $A_C$ by an orthogonal
projection matrix.
As we shall prove in
Proposition~\ref{prop:nearest-subspaces} below,
\begin{equation}\label{eq:intro-nearest-subspace-distance}
 d_{\Pi,2}(C,\LS_d)^2
 =
 \sum_{i=1}^d\min\{\lambda_i,1-\lambda_i\},
\end{equation}
where $\lambda_1,\ldots,\lambda_d$ are the
eigenvalues of $A_C$.
The nearest subspaces are obtained by taking
the orthogonal sum of all eigenspaces with
eigenvalues greater than $1/2$ and an arbitrary
subspace of the eigenspace with eigenvalue $1/2$.
Thus $A_C$ determines both the nearest linear
subspaces and their distance from $C$.

Restricting a linear approximation to be an
orthogonal projection changes the optimal
squared error by at most a universal factor.
More precisely, we shall establish in
Theorem~\ref{thm:variance-distance-subspaces}
below that
\begin{equation}\label{eq:intro-variance-distance-subspaces}
 a
 \leq d_{\Pi,2}(C,\LS_d)^2\leq2a,
 \qquad
 \frac12\tau^2
 \leq d_{\Pi,2}(C,\LS_d)^2\leq3\tau^2.
\end{equation}

For nonnegative quantities $u$ and $v$, we write
$u\asymp v$ when $c_1v\leq u\leq c_2v$ for
positive constants $c_1,c_2$ independent of
the cone and the ambient dimension.
With this notation,
\eqref{eq:intro-statistical-dimension-distance}
and \eqref{eq:intro-variance-distance-subspaces}
give the parallel
\begin{equation}\label{eq:intro-metric-mean-variance}
 \delta(C)=d_{\Pi,2}(C,\{0\})^2,
 \qquad
 \tau^2(C)\asymp d_{\Pi,2}(C,\LS_d)^2.
\end{equation}
The statistical dimension is the squared distance
to the zero cone, while the intrinsic-volume
variance is comparable to the squared distance
to the family of linear subspaces.

These comparisons also quantify the lack of
central symmetry discussed in
Subsection~\ref{subsec:intro-cone-chaoses}.
Since the centrally symmetric convex cones are
precisely the linear subspaces,
$d_{\Pi,2}(C,\LS_d)$ is the distance from
$C$ to the family of centrally symmetric cones.
Together with
\eqref{eq:intro-variance-comparisons},
\eqref{eq:intro-variance-distance-subspaces}
shows that this distance is comparable to
the norm of the mean projection
$m_C=\E\Pi_C(g)$.

Another measure of the lack of central symmetry
is the distance from $C$ to its reflection $-C$.
As we shall prove in
Corollary~\ref{cor:distance-reflected-cone}
below,
\begin{equation}\label{eq:intro-reflection-distance}
 2\norm{m_C}
 \leq d_{\Pi,2}(C,-C)
 \leq2\sqrt a.
\end{equation}
The variance comparisons from the preceding
subsection therefore imply
\begin{equation}\label{eq:intro-reflection-variance}
 2\tau^2
 \leq d_{\Pi,2}(C,-C)^2
 \leq6\tau^2.
\end{equation}
Together with our previous results, these
inequalities give
\begin{equation}\label{eq:intro-geometric-comparisons}
 \tau^2
 \asymp \norm{m_C}^2
 \asymp d_{\Pi,2}(C,\LS_d)^2
 \asymp d_{\Pi,2}(C,-C)^2.
\end{equation}
Thus the intrinsic-volume variance measures,
up to universal constants, both the squared
distance to the family of centrally symmetric
cones and the squared distance to the reflected
cone itself.

The spectral formula
\eqref{eq:intro-nearest-subspace-distance}
also identifies the cones that are farthest
from the family of linear subspaces.
In Proposition~\ref{prop:cones-farthest-from-subspaces}
below, we shall show that
$d_{\Pi,2}(C,\LS_d)\leq\sqrt{d/2}$,
with equality exactly when $A_C=I_d/2$.
In particular, every self-dual cone attains
this maximum.

For polyhedral cones, approximation by linear
subspaces has a further probabilistic
interpretation.
Let $\mathsf S_C$ be the linear span of the
face whose relative interior contains
$\Pi_C(g')$, where $g'$ is an independent
copy of $g$.
We shall prove in
Proposition~\ref{prop:random-face-spans}
below that, for every $L\in\LS_d$,
\begin{equation}\label{eq:intro-random-face-spans}
 d_{\Pi,2}(C,L)^2
 =
 \E\bigl[d_{\Pi,2}(\mathsf S_C,L)^2\bigr].
\end{equation}
The expectation on the right is over the random
face span, and the distance inside it is the
metric between the realized subspace and $L$.
Consequently, a nearest linear subspace to $C$
is precisely a Fr\'echet mean of the distribution
of $\mathsf S_C$ in the metric space
$(\LS_d,d_{\Pi,2})$.
Here a Fr\'echet mean means a subspace minimizing
the expected squared distance to the random
subspace; its dimension is also allowed to vary.

Finally, the variance comparisons give metric
characterizations of asymptotic normality,
complementing
Corollary~\ref{cor:intro-chaos-clt}.

\begin{corollary}[Metric characterizations of asymptotic normality]
\label{cor:intro-metric-clt}
Let $C_n\subseteq\R^{d_n}$ be nonempty closed
convex cones, and write $V_n=V_{C_n}$,
$\delta_n=\E V_n$, and $\tau_n^2=\Var(V_n)$.
Assume that $\tau_n^2>0$ for all sufficiently
large $n$.
Then the following conditions are equivalent:
\begin{enumerate}
 \item
 $(V_n-\delta_n)/\tau_n\toD\normal(0,1)$;
 \item
 $d_{\Pi,2}(C_n,\LS_{d_n})\to\infty$;
 \item
 $d_{\Pi,2}(C_n,-C_n)\to\infty$.
\end{enumerate}
\end{corollary}

\begin{proof}
By Corollary~\ref{cor:intro-chaos-clt}, the
first condition is equivalent to
$\tau_n^2\to\infty$.
The equivalence with the other two conditions
follows from
\eqref{eq:intro-variance-distance-subspaces}
and \eqref{eq:intro-reflection-variance}.
\end{proof}

Thus the intrinsic-volume central limit theorem
holds precisely when the cones move arbitrarily
far from the family of linear subspaces in
the $L^2$-$\Pi$-metric.
Equivalently, their distances from their own
reflections tend to infinity.
In Section~\ref{sec:proof-clt}, we explain how the
subtraction in \eqref{eq:intro-F} is incorporated into
\eqref{eq:steiner-laplace} to prove Theorem~\ref{thm:clt}.

\subsection{Organization of the paper}

In Section~\ref{sec:gaussian}, we develop a quantitative
Gaussian approximation for functions of Gaussian vectors.
The main result of this section is
Lemma~\ref{lem:gaussian-laplace}, which provides a
Laplace-transform estimate for functions with a small
Lipschitz constant of the gradient.
As a first application, in
Section~\ref{sec:squared-projection-clt}, we give an
alternative proof of the central limit theorem for the
squared Gaussian projection norm, stated in
Theorem~\ref{thm:squared-projection-GNP}.

Section~\ref{sec:hermite} recalls the relevant properties
of Hermite expansions in $L^2(\R^d,\gamma_d)$ and
Gaussian chaoses. In particular, we identify the
first three Gaussian-chaos components of a function
and establish a refined Poincar\'e inequality for
functions whose second chaos vanishes.

In Section~\ref{sec:removal}, we study the Wiener-chaos
expansions of Gaussian cone projections and their squared
norms, establish the geometric results announced in
Subsections~\ref{subsec:intro-cone-chaoses}
and~\ref{subsec:intro-projection-distance}, and investigate
the effect of removing the second Gaussian-chaos component.
We also study circular cones as examples illustrating
the relative importance of the individual chaos components.

In Section~\ref{sec:proof-clt}, we incorporate the
quadratic-chaos subtraction into the master Steiner formula
by an exponential change of variables, establish the
quantitative Laplace-transform estimate of
Theorem~\ref{thm:laplace-main}, and prove the central
and local limit theorems for conic intrinsic volumes.
The appendix contains a continuity lemma for Laplace
transforms.

\section{Gaussian approximation for functions of Gaussian vectors}
\label{sec:gaussian}

We write $\gamma_d$ for the standard Gaussian probability measure
on $\R^d$, given by
$\gamma_d(\dd x)=(2\pi)^{-d/2}\eee^{-\|x\|^2/2}\dd x$,
and let $g$ be a random vector with distribution $\gamma_d$.
For $1\leq p<\infty$, the space $L^p(\R^d,\gamma_d)$ consists of
measurable real-valued functions $f$ such that
$\E|f(g)|^p<\infty$, with functions that agree almost everywhere
identified. Its norm is
$\|f\|_{L^p(\gamma_d)}=(\E|f(g)|^p)^{1/p}$.
In particular, $L^2(\R^d,\gamma_d)$ is a Hilbert space with inner
product
$\langle f,h\rangle_{L^2(\gamma_d)}=\E[f(g)h(g)]$.
For a real or Euclidean-valued random variable $Y$, we also use
$\|Y\|_{L^p}=(\E\|Y\|^p)^{1/p}$, with the absolute value in the
real-valued case.

A central part of our proof is the analysis of functions in
$L^2(\R^d,\gamma_d)$. In this section, we study how the derivatives
of a function $f$ control the distribution of $f(g)$.
An affine function of $g$ is Gaussian. Our main approximation
result quantifies this observation for functions whose gradients
have a small Lipschitz constant.

More precisely, Lemma~\ref{lem:gaussian-laplace} will show the
following. Suppose that $f\in C^1(\R^d)$ satisfies
$\E f(g)=0$ and $\E\|\nabla f(g)\|^2\leq M^2$, and that
$\nabla f$ has Lipschitz constant at most
$\varepsilon\in[0,1]$.
Then there are constants $r>0$ and $K_M<\infty$, depending only
on $M$, such that
\[
 \left|
   \E\eee^{t f(g)}
   -
   \eee^{t^2\Var(f(g))/2}
 \right|
 \leq K_M\varepsilon|t|^3,
 \qquad |t|\leq r.
\]
The constants are independent of the dimension.
The formulation in terms of Laplace transforms is useful for
our purposes because \eqref{eq:steiner-laplace} relates the
Laplace transform of $V_C$ to that of a squared Gaussian
projection.

The argument uses the method of second-order Poincar\'e
inequalities developed by Chatterjee~\cite{Cha09} and
Nourdin, Peccati, and Reinert~\cite{NPR09}.
We give the details needed for the estimate above.
We first specify the regularity assumptions and establish
Gaussian integration by parts.
We then use interpolation between Gaussian vectors to obtain
a covariance identity and the Gaussian Poincar\'e inequality.
Finally, moment estimates and a second application of the
covariance identity yield the approximation of Laplace
transforms.

\subsection{Weak derivatives and smooth approximation}
\label{subsec:weak-derivatives}

The functions used later have Lipschitz gradients, but their
second derivatives need not be continuous.
We therefore formulate the analytic identities using weak
derivatives. These are defined with respect to Lebesgue
measure; Gaussian measure enters through the integrability
assumptions.

We write $C_c^\infty(\R^d)$ for the space of infinitely
differentiable functions with compact support.

\begin{definition}[Weak derivatives]
\label{def:weak-derivatives}
Let $k\geq1$ be an integer, and let $f$ be locally integrable
on $\R^d$ with respect to Lebesgue measure.
For indices $i_1,\ldots,i_k\in\{1,\ldots,d\}$, a locally
integrable function $h$ is called the weak derivative
$\partial_{i_1}\cdots\partial_{i_k}f$ if
\begin{equation}\label{eq:weak-derivative-definition}
 \int_{\R^d}
   f(x)\,
   \partial_{i_1}\cdots\partial_{i_k}\varphi(x)\,\dd x
 =
 (-1)^k
 \int_{\R^d}
   h(x)\varphi(x)\,\dd x
 \qquad
 \text{for every }\varphi\in C_c^\infty(\R^d).
\end{equation}
We say that $f$ has weak derivatives through order $k$ if these
derivatives exist for every choice of indices and every order
from $1$ to $k$.
\end{definition}

Weak derivatives are unique up to equality almost everywhere.
For a $C^k$ function, they agree with the corresponding
classical derivatives.
Since the density of $\gamma_d$ is strictly positive,
Lebesgue measure and $\gamma_d$ have the same null sets.
Moreover, on each compact set the Gaussian density is bounded
above and bounded away from zero.
Consequently, every function in $L^p(\R^d,\gamma_d)$,
$1\leq p<\infty$, belongs locally to the corresponding
Lebesgue $L^p$ space and is locally integrable.

We use the notation
$\nabla f=(\partial_1f,\ldots,\partial_df)$ for the weak gradient
and
$\nabla^2f=(\partial_i\partial_jf)_{i,j=1}^d$ for the weak
Hessian, whenever these derivatives exist.
Weak mixed derivatives commute, so the weak Hessian is
symmetric almost everywhere.

A locally Lipschitz function is differentiable almost
everywhere by Rademacher's theorem, and its classical
derivatives represent its weak first derivatives.
In particular, if $f\in C^1(\R^d)$ and $\nabla f$ is locally
Lipschitz, then $f$ has weak second derivatives.
This observation will allow us to apply the results of this
section to squared projection functions.

The following approximation lemma allows identities proved
for smooth functions with compact support to be extended to
functions with weak derivatives.

\begin{lemma}[Smooth approximation]
\label{lem:smooth-approximation}
Let $k\geq1$ and $1\leq p<\infty$.
Suppose that $f$ and all its weak derivatives of orders
$1,\ldots,k$ belong to $L^p(\R^d,\gamma_d)$.
Then there is a sequence $f_m\in C_c^\infty(\R^d)$ such that
\begin{equation}\label{eq:smooth-approximation}
 \begin{aligned}
  \|f_m-f\|_{L^p(\gamma_d)}
    &\longrightarrow0,\\
  \bigl\|
    \partial_{i_1}\cdots\partial_{i_\ell}f_m
    -
    \partial_{i_1}\cdots\partial_{i_\ell}f
  \bigr\|_{L^p(\gamma_d)}
    &\longrightarrow0
 \end{aligned}
\end{equation}
for every $1\leq\ell\leq k$ and every choice of
$i_1,\ldots,i_\ell\in\{1,\ldots,d\}$.
\end{lemma}

\begin{proof}
We first introduce a cutoff and then smooth by convolution.

Choose $\chi\in C_c^\infty(\R^d)$ such that
$0\leq\chi\leq1$, $\chi(x)=1$ when $\|x\|\leq1$, and
$\chi(x)=0$ when $\|x\|\geq2$.
For $R\geq1$, set $\chi_R(x)=\chi(x/R)$.
Dominated convergence gives
$\|\chi_Rf-f\|_{L^p(\gamma_d)}\to0$ as $R\to\infty$.

The same convergence holds for every derivative through
order $k$.
Indeed, the weak product rule expresses a derivative of
$\chi_Rf$ as $\chi_R$ times the corresponding derivative
of $f$, together with terms in which at least one derivative
falls on $\chi_R$.
The first term converges to the corresponding derivative
of $f$ by dominated convergence.
In each remaining term, a derivative of $\chi_R$ of some
order $j\geq1$ multiplies a lower derivative of $f$.
That derivative of $\chi_R$ is supported in
$\{x:R\leq\|x\|\leq2R\}$ and is bounded in absolute value by
a constant times $R^{-j}$.
The $L^p(\gamma_d)$ norm of the term therefore tends to zero,
because the lower derivative of $f$ belongs to
$L^p(\R^d,\gamma_d)$.
Thus $\chi_Rf$ approximates $f$, together with all the
required derivatives.

We next fix $R$ and smooth $\chi_Rf$.
Choose a nonnegative function
$\rho\in C_c^\infty(\R^d)$ supported in the unit ball and
satisfying $\int_{\R^d}\rho(x)\,\dd x=1$.
For $0<\eta<1$, put
$\rho_\eta(x)=\eta^{-d}\rho(x/\eta)$.
The function $(\chi_Rf)*\rho_\eta$ is smooth and has support
in the ball of radius $2R+1$.

The function $\chi_Rf$ and all its weak derivatives through
order $k$ belong to the ordinary Lebesgue $L^p$ space:
they have compact support, and the Gaussian density is
bounded away from zero on that support.
Convolution commutes with weak differentiation.
The continuity of translations in Lebesgue $L^p$ therefore
implies that, as $\eta\downarrow0$,
$(\chi_Rf)*\rho_\eta$ and each of its derivatives through
order $k$ converge in Lebesgue $L^p$ to $\chi_Rf$ and the
corresponding derivative.
Since all these functions are supported in a fixed compact
set, this convergence also holds in $L^p(\R^d,\gamma_d)$.

For each integer $m\geq1$, take $R=m$ and choose
$\eta_m\in(0,1)$ so small that every approximation error
from this convolution step is at most $1/m$.
There are only finitely many derivatives under consideration.
The sequence
$f_m=(\chi_mf)*\rho_{\eta_m}$
then satisfies \eqref{eq:smooth-approximation}.
\end{proof}

If $f\in C^1(\R^d)$ has a globally Lipschitz gradient, then
the assumptions of Lemma~\ref{lem:smooth-approximation} hold
with $k=2$ for every finite $p$.
Indeed, $\nabla f$ has at most linear growth, $f$ has at most
quadratic growth, and its weak second derivatives are bounded.
All these functions are therefore integrable to every finite
power under $\gamma_d$.

\subsection{Gaussian integration by parts}
\label{subsec:gaussian-integration-by-parts}

Differentiating the Gaussian density introduces a coordinate
factor. Integration by parts consequently relates moments
involving the coordinates of $g$ to expectations of
derivatives.
We record the first two identities, together with a version
involving a polynomial that will also be used in the Hermite
calculations.

\begin{lemma}[Gaussian integration by parts]
\label{lem:ibp}
Let $1<p<\infty$, and suppose that $f$ and all its weak first
derivatives belong to $L^p(\R^d,\gamma_d)$.
Then, for every $i\in\{1,\ldots,d\}$,
\begin{equation}\label{eq:ibp-first}
 \E[g_i f(g)]=\E[\partial_i f(g)].
\end{equation}
If, in addition, $f$ has weak second derivatives belonging
to $L^p(\R^d,\gamma_d)$, then, for every
$i,j\in\{1,\ldots,d\}$,
\begin{equation}\label{eq:ibp-second}
 \E\bigl[
   (g_i g_j-\mathbf{1}_{\{i=j\}})f(g)
 \bigr]
 =
 \E[\partial_i\partial_j f(g)].
\end{equation}
More generally, under the assumptions for
\eqref{eq:ibp-first}, every polynomial
$P:\R^d\to\R$ satisfies
\begin{equation}\label{eq:ibp-polynomial}
 \E[\partial_i f(g)\,P(g)]
 =
 \E\bigl[
   f(g)\bigl(g_iP(g)-\partial_iP(g)\bigr)
 \bigr].
\end{equation}
All expectations in these formulas are absolutely
convergent.
\end{lemma}

\begin{proof}
We first prove \eqref{eq:ibp-polynomial}.
For $f\in C_c^\infty(\R^d)$, ordinary integration by parts
gives
\[
 \begin{aligned}
 \E[\partial_i f(g)\,P(g)]
 &=
 (2\pi)^{-d/2}
 \int_{\R^d}
   \partial_i f(x)\,P(x)\,
   \eee^{-\|x\|^2/2}\,\dd x\\
 &=
 -(2\pi)^{-d/2}
 \int_{\R^d}
   f(x)\,
   \partial_i\bigl(P(x)\eee^{-\|x\|^2/2}\bigr)\,\dd x\\
 &=
 \E\bigl[
   f(g)\bigl(g_iP(g)-\partial_iP(g)\bigr)
 \bigr].
 \end{aligned}
\]
There is no boundary term because $f$ has compact support.

Now suppose that $f$ satisfies the first-order assumptions
of the lemma.
By Lemma~\ref{lem:smooth-approximation}, there are
$f_m\in C_c^\infty(\R^d)$ such that $f_m\to f$ and
$\partial_i f_m\to\partial_i f$ in $L^p(\R^d,\gamma_d)$.
Every polynomial belongs to
$L^{p/(p-1)}(\R^d,\gamma_d)$.
Hölder's inequality therefore permits passage to the limit
on both sides of \eqref{eq:ibp-polynomial}.
It also gives the asserted absolute integrability.
Taking $P=1$ proves \eqref{eq:ibp-first}.

Finally, suppose that the weak second derivatives of $f$
also belong to $L^p(\R^d,\gamma_d)$.
Apply \eqref{eq:ibp-first} to $\partial_i f$ in coordinate
$j$, and apply \eqref{eq:ibp-polynomial} to $f$ with
$P(x)=x_j$.
Since weak mixed derivatives commute, these identities give
\[
 \E[\partial_i\partial_j f(g)]
 =
 \E[g_j\,\partial_i f(g)]
 =
 \E\bigl[
   (g_i g_j-\mathbf{1}_{\{i=j\}})f(g)
 \bigr],
\]
which proves \eqref{eq:ibp-second}.
\end{proof}

\subsection{Gaussian covariance interpolation and Poincar\'e's inequality}
\label{subsec:gaussian-covariance}

We recall two classical results from Gaussian analysis:
the Gaussian covariance identity
(see \cite[Equation~(1.1)]{BGH01}) and the Gaussian
Poincar\'e inequality
(see \cite[Theorem~1.6.4]{Bog98}).
We give proofs under the weak differentiability assumptions
needed below.
In addition to yielding \eqref{eq:poincare}, the
representation \eqref{eq:covariance} will be used in the
proof of Lemma~\ref{lem:gaussian-laplace}.

Let $g'$ be an independent copy of the standard Gaussian
vector $g$.
For $f,h\in L^2(\R^d,\gamma_d)$, the covariance
$\Cov(f(g),h(g))$ is the difference between
$\E[f(g)h(g)]$ and $\E[f(g)h(g')]$.
This suggests interpolating between $g'$ and $g$ while
preserving the standard Gaussian distribution.
For $r\in[0,1]$, define
$g_r=rg+\sqrt{1-r^2}\,g'$.
Each $g_r$ is standard Gaussian, and the endpoint values
are $g_0=g'$ and $g_1=g$.
The following identity expresses the covariance as the
integral of the change in $\E[f(g)h(g_r)]$ along this
interpolation.

\begin{lemma}[Gaussian covariance identity]
\label{lem:covariance}
Let $f,h:\R^d\to\R$ have weak first partial derivatives,
and suppose that $f,h$ and all these derivatives belong
to $L^2(\R^d,\gamma_d)$.
Then
\begin{equation}\label{eq:covariance}
 \Cov(f(g),h(g))
 =
 \int_0^1
   \E\langle\nabla f(g),\nabla h(g_r)\rangle
 \,\dd r.
\end{equation}
\end{lemma}

\begin{proof}
We first assume that $f,h\in C_c^\infty(\R^d)$.
For $0<r<1$, differentiation under the expectation gives
\[
 \frac{\dd}{\dd r}\E[f(g)h(g_r)]
 =
 \sum_{i=1}^d
 \E\left[
   f(g)\partial_i h(g_r)
   \left(
     g_i-\frac{r}{\sqrt{1-r^2}}g_i'
   \right)
 \right].
\]
This differentiation is justified on every compact
subinterval of $(0,1)$ by the boundedness of $f$ and
$\nabla h$ and the integrability of the Gaussian
coordinates.

The pair $(g,g')$ is a standard Gaussian vector in
$\R^{2d}$, so \eqref{eq:ibp-first} applies to either set
of coordinates.
Integration by parts in $g_i$ gives
\[
 \begin{aligned}
 \E[g_i f(g)\partial_i h(g_r)]
 &=
 \E[\partial_i f(g)\partial_i h(g_r)]
 +
 r\,\E[f(g)\partial_i^2 h(g_r)],
 \end{aligned}
\]
whereas integration by parts in $g_i'$ gives
\[
 \E[g_i' f(g)\partial_i h(g_r)]
 =
 \sqrt{1-r^2}\,
 \E[f(g)\partial_i^2 h(g_r)].
\]
Substituting these two identities into the derivative
above cancels the terms involving second derivatives
of $h$.
We obtain
\begin{equation}\label{eq:covariance-derivative}
 \frac{\dd}{\dd r}\E[f(g)h(g_r)]
 =
 \E\langle\nabla f(g),\nabla h(g_r)\rangle,
 \qquad 0<r<1.
\end{equation}

We next justify integration up to the endpoints.
Since $g_r$ has the same distribution as $g$, the
Cauchy--Schwarz inequality gives
\[
 \left|
   \E\langle\nabla f(g),\nabla h(g_r)\rangle
 \right|
 \leq
 \|\nabla f(g)\|_{L^2}
 \|\nabla h(g_r)\|_{L^2}\\
 =
 \|\nabla f(g)\|_{L^2}
 \|\nabla h(g)\|_{L^2}.
\]
Thus the derivative in \eqref{eq:covariance-derivative}
is bounded uniformly in $r$.
We may integrate that identity over a compact subinterval
of $(0,1)$ and then let its endpoints tend to $0$ and $1$.
The boundedness and continuity of $f$ and $h$ justify
the corresponding limits of the expectations by
dominated convergence.
Consequently,
\[
 \E[f(g)h(g)]-\E[f(g)h(g')]
 =
 \int_0^1
   \E\langle\nabla f(g),\nabla h(g_r)\rangle
 \,\dd r.
\]
Independence gives
$\E[f(g)h(g')]=\E[f(g)]\,\E[h(g)]$,
which proves \eqref{eq:covariance} for smooth functions
with compact support.

Now let $f,h$ satisfy the assumptions of the lemma.
By Lemma~\ref{lem:smooth-approximation}, choose
$f_m,h_m\in C_c^\infty(\R^d)$ that converge to $f,h$,
respectively, together with their weak gradients in
Gaussian $L^2$.
The covariances
$\Cov(f_m(g),h_m(g))$ then converge to
$\Cov(f(g),h(g))$ by the Cauchy--Schwarz inequality.

For the right-hand side of \eqref{eq:covariance}, the
same inequality gives, uniformly in $r\in[0,1]$,
\[
 \begin{aligned}
 &\left|
   \E\left[
     \langle\nabla f_m(g),\nabla h_m(g_r)\rangle
     -
     \langle\nabla f(g),\nabla h(g_r)\rangle
   \right]
 \right|\\
 &\quad\leq
 \|\nabla f_m(g)-\nabla f(g)\|_{L^2}
 \|\nabla h_m(g)\|_{L^2}\\
 &\qquad\quad+
 \|\nabla f(g)\|_{L^2}
 \|\nabla h_m(g)-\nabla h(g)\|_{L^2}.
 \end{aligned}
\]
Here we again used that $g_r$ and $g$ have the same
distribution.
The right-hand side tends to zero, because the gradients
converge in $L^2$.
Passing to the limit in \eqref{eq:covariance} for
$f_m,h_m$ proves the result.
\end{proof}

Taking $h=f$ in \eqref{eq:covariance} yields the Gaussian
Poincar\'e inequality.

\begin{corollary}[Gaussian Poincar\'e inequality]
\label{cor:poincare}
Let $f:\R^d\to\R$ have weak first partial derivatives,
and suppose that $f$ and all these derivatives belong
to $L^2(\R^d,\gamma_d)$.
Then
\begin{equation}\label{eq:poincare}
 \Var(f(g))
 \leq
 \E\|\nabla f(g)\|^2.
\end{equation}
\end{corollary}

\begin{proof}
For every $r\in[0,1]$, the Cauchy--Schwarz inequality
and the equality in distribution of $g_r$ and $g$ give
\[
 \left|
   \E\langle\nabla f(g),\nabla f(g_r)\rangle
 \right|
 \leq
 \bigl(\E\|\nabla f(g)\|^2\bigr)^{1/2}
 \bigl(\E\|\nabla f(g_r)\|^2\bigr)^{1/2}\\
 =
 \E\|\nabla f(g)\|^2.
\]
Integrating this bound over $r\in[0,1]$ and applying
\eqref{eq:covariance} with $h=f$ proves
\eqref{eq:poincare}.
\end{proof}

The constant one in \eqref{eq:poincare} is optimal
in every dimension.
Indeed, if $f(x)=\langle a,x\rangle+b$ with
$a\in\R^d\setminus\{0\}$ and $b\in\R$, then both sides
of \eqref{eq:poincare} equal $\|a\|^2$.

Inequality~\eqref{eq:poincare} controls variances.
To compare Laplace transforms, we will also need
bounds on higher moments that imply exponential
integrability.
The next subsection establishes the required moment
inequality.

\subsection{Moment and exponential inequalities from Gaussian rotation}
\label{subsec:gaussian-rotation}

Inequality~\eqref{eq:poincare} controls the second moment
of a centered function of a Gaussian vector.
We now establish an inequality that also gives bounds
on higher moments and Laplace transforms.
Its formulation involves an arbitrary nonnegative convex
function, so the same argument applies to both powers
and exponentials.

The following inequality is due to Pisier
\cite[Theorem~2.2]{Pisier86}.
We present Maurey's proof based on Gaussian rotations,
as given in that reference.

\begin{lemma}[Maurey--Pisier inequality]
\label{lem:maurey-pisier}
Let $f:\R^d\to\R$ be locally Lipschitz and satisfy
$\E|f(g)|<\infty$.
Let $Z$ be a one-dimensional standard Gaussian random
variable independent of $g$.
Then, for every convex function
$\Psi:\R\to[0,\infty)$,
\begin{equation}\label{eq:maurey-pisier}
 \E\Psi\bigl(f(g)-\E f(g)\bigr)
 \leq
 \E\Psi\left(
   \frac{\pi}{2}\|\nabla f(g)\|Z
 \right).
\end{equation}
The gradient is understood almost everywhere, and the
expectations in \eqref{eq:maurey-pisier} are allowed
to be infinite.
\end{lemma}

Conditionally on $g$, the random variable on the
right-hand side of \eqref{eq:maurey-pisier} is Gaussian
with mean zero and standard deviation
$(\pi/2)\|\nabla f(g)\|$.
Thus the inequality controls convex functions of the
centered value $f(g)$ through Gaussian averages whose
size is determined by the gradient.
The independence of $Z$ and $g$ makes these averages
particularly easy to evaluate.

\begin{proof}
Let $g'$ be an independent copy of $g$.
Since
$f(g)-\E f(g)=\E[f(g)-f(g')\mid g]$,
conditional Jensen's inequality gives
\begin{equation}\label{eq:rotation-centering}
 \E\Psi\bigl(f(g)-\E f(g)\bigr)
 \leq
 \E\Psi\bigl(f(g)-f(g')\bigr).
\end{equation}
We will estimate the difference $f(g)-f(g')$ by
integrating a directional derivative along a Gaussian
rotation.

Recall the interpolation
$g_r=rg+\sqrt{1-r^2}\,g'$ from
Subsection~\ref{subsec:gaussian-covariance}.
In the angular parametrization,
$g_{\sin\theta}=g\sin\theta+g'\cos\theta$ for
$0\leq\theta\leq\pi/2$.
This path connects $g'$ to $g$, and its velocity is
$g\cos\theta-g'\sin\theta$.
For each fixed $\theta$, the two vectors
\[
 g_{\sin\theta}
 =
 g\sin\theta+g'\cos\theta,
 \qquad
 g\cos\theta-g'\sin\theta
\]
are independent standard Gaussian vectors.
Indeed, together they form an orthogonal transformation
of the standard Gaussian vector $(g,g')$ in $\R^{2d}$.

We first justify differentiation along this path.
For every realization of $g,g'$, the path is smooth
and has compact image.
Since $f$ is locally Lipschitz, the function
$\theta\mapsto f(g_{\sin\theta})$ is Lipschitz on
$[0,\pi/2]$ and hence absolutely continuous.
By Rademacher's theorem, the set where $f$ is not
differentiable has Lebesgue measure zero.
For every fixed $\theta$, the vector $g_{\sin\theta}$
has distribution $\gamma_d$, so it belongs to this
exceptional set with probability zero.
Fubini's theorem therefore shows that, almost surely,
$f$ is differentiable at $g_{\sin\theta}$ for almost
every $\theta$.
The chain rule and the fundamental theorem of calculus
then give
\[
 f(g)-f(g')
 =
 \int_0^{\pi/2}
 \left\langle
   \nabla f(g_{\sin\theta}),
   g\cos\theta-g'\sin\theta
 \right\rangle
 \,\dd\theta
 \qquad\text{almost surely}.
\]

Apply Jensen's inequality to this integral, using the
uniform probability measure $(2/\pi)\,\dd\theta$ on
$[0,\pi/2]$.
After taking expectations, we obtain
\[
 \E\Psi\bigl(f(g)-f(g')\bigr)
 \leq
 \frac{2}{\pi}
 \int_0^{\pi/2}
 \E\Psi\left(
   \frac{\pi}{2}
   \left\langle
     \nabla f(g_{\sin\theta}),
     g\cos\theta-g'\sin\theta
   \right\rangle
 \right)
 \,\dd\theta.
\]
Nonnegativity of $\Psi$ permits the interchange of
expectation and integration even when the resulting
expectation is infinite.

For every fixed $\theta$, the velocity
$g\cos\theta-g'\sin\theta$ is standard Gaussian and
independent of $g_{\sin\theta}$.
Conditionally on $g_{\sin\theta}$, the inner product
in the last display is therefore Gaussian with mean
zero and variance
$\|\nabla f(g_{\sin\theta})\|^2$.
Since $g_{\sin\theta}$ has the same distribution as $g$,
\[
 \E\Psi\left(
   \frac{\pi}{2}
   \left\langle
     \nabla f(g_{\sin\theta}),
     g\cos\theta-g'\sin\theta
   \right\rangle
 \right)
 =
 \E\Psi\left(
   \frac{\pi}{2}\|\nabla f(g)\|Z
 \right).
\]
The last expression is independent of $\theta$.
Substituting it into the preceding integral and using
\eqref{eq:rotation-centering} proves
\eqref{eq:maurey-pisier}.
\end{proof}

The essential feature of the rotation is that its
position and velocity are independent standard Gaussian
vectors at each fixed parameter value.
Conditionally on the position, the directional derivative
along the path is Gaussian with mean zero and standard
deviation equal to the norm of the gradient at that position. This is why the bound in \eqref{eq:maurey-pisier}
contains no additional factor depending on the dimension.
The factor $\pi/2$ comes from the length of the
interval of rotation.

Taking a power function in \eqref{eq:maurey-pisier}
gives the following moment inequality.
We state it for weakly differentiable functions, using
the approximation result from Subsection~\ref{subsec:weak-derivatives}.

\begin{corollary}[Gaussian moment inequality]
\label{lem:moment}
Let $2\leq q<\infty$.
Suppose that $f:\R^d\to\R$ and all its weak first
partial derivatives belong to $L^q(\R^d,\gamma_d)$.
Then
\begin{equation}\label{eq:moment}
 \|f(g)-\E f(g)\|_{L^q}
 \leq
 \frac{\pi}{2}\sqrt{q-1}\,
 \|\nabla f(g)\|_{L^q}.
\end{equation}
\end{corollary}

\begin{proof}
First suppose that $f\in C_c^\infty(\R^d)$.
Apply \eqref{eq:maurey-pisier} with
$\Psi(u)=|u|^q$.
The independence of $Z$ and $g$ gives
\[
 \E|f(g)-\E f(g)|^q
 \leq
 \left(\frac{\pi}{2}\right)^q
 \E\|\nabla f(g)\|^q\,
 \E|Z|^q.
\]
It remains to bound the Gaussian moment.

For $q>2$, Gaussian integration by parts
\eqref{eq:ibp-first} in dimension one, applied to
$z\mapsto z|z|^{q-2}$, gives
\[
 \E|Z|^q
 =
 (q-1)\E|Z|^{q-2}
 \leq
 (q-1)\bigl(\E|Z|^q\bigr)^{(q-2)/q},
\]
where the inequality follows from H\"older's inequality.
Consequently,
$\|Z\|_{L^q}\leq\sqrt{q-1}$.
The same bound holds for $q=2$, since $\E Z^2=1$.
Taking $q$th roots proves \eqref{eq:moment} for
smooth functions with compact support.

For a general $f$ satisfying the assumptions, apply
Lemma~\ref{lem:smooth-approximation} with $p=q$ and
$k=1$.
It gives functions $f_m\in C_c^\infty(\R^d)$ converging
to $f$, together with their gradients, in Gaussian
$L^q$.
Moreover, $\E f_m(g)\to\E f(g)$, because convergence
in $L^q$ implies convergence in $L^1$ on a probability
space.
We may therefore pass to the limit in
\eqref{eq:moment} for $f_m$.
\end{proof}

Taking an exponential function in
\eqref{eq:maurey-pisier} gives a direct bound on the
Laplace transform.

\begin{corollary}[Gaussian exponential inequality]
\label{cor:gaussian-exponential}
Let $f:\R^d\to\R$ be locally Lipschitz and satisfy
$\E|f(g)|<\infty$.
Then, for every $t\in\R$,
\begin{equation}\label{eq:gaussian-exponential}
 \E\eee^{t(f(g)-\E f(g))}
 \leq
 \E\eee^{\frac{\pi^2t^2}{8}\|\nabla f(g)\|^2}.
\end{equation}
The inequality holds with possibly infinite
expectations.
In particular, if $f$ is globally $L$-Lipschitz, then
\begin{equation}\label{eq:lipschitz-exponential}
 \E\eee^{t(f(g)-\E f(g))}
 \leq
 \eee^{\pi^2t^2L^2/8},
 \qquad t\in\R.
\end{equation}
\end{corollary}

\begin{proof}
Apply \eqref{eq:maurey-pisier} with
$\Psi(u)=\eee^{tu}$.
Conditioning on $g$ and evaluating the
one-dimensional Gaussian Laplace transform gives
\[
 \E\eee^{t(f(g)-\E f(g))}
 \leq
 \E\eee^{(\pi t/2)\|\nabla f(g)\|Z}
 =
 \E\eee^{\frac{\pi^2t^2}{8}\|\nabla f(g)\|^2},
\]
which proves \eqref{eq:gaussian-exponential}.

If $f$ is globally $L$-Lipschitz, then it has at most
linear growth and is therefore integrable under
$\gamma_d$.
Moreover, $\|\nabla f(x)\|\leq L$ almost everywhere.
Substituting this bound into
\eqref{eq:gaussian-exponential} proves
\eqref{eq:lipschitz-exponential}.
\end{proof}

Inequality~\eqref{eq:gaussian-exponential} reduces
the control of the Laplace transform of a centered
function to the control of exponential moments of
its squared gradient norm.
In the next subsection, we will use this inequality
for functions with globally Lipschitz gradients.
We will bound its right-hand side by applying
\eqref{eq:maurey-pisier} once more, to the function
$x\mapsto\|\nabla f(x)\|$.

\subsection{A small Hessian gives a Gaussian approximation}
\label{subsec:gaussian-laplace}

We now combine the preceding estimates to approximate the
Laplace transform of a function of a Gaussian vector.
If $f$ is affine, then $f(g)$ is Gaussian and its gradient
is constant. We will show that a small Lipschitz constant
of $\nabla f$, together with a bound on
$\E\|\nabla f(g)\|^2$, gives a quantitative approximation
by a Gaussian variable with the same mean and variance.

This is the method of second-order Gaussian Poincar\'e
inequalities introduced by Chatterjee
\cite[Theorem~2.2 and Lemma~5.3]{Cha09}
and extended to general Wiener functionals by
Nourdin, Peccati, and Reinert
\cite[Theorem~1.1]{NPR09}.
Goldstein, Nourdin, and Peccati~\cite{GNP17}
applied this method to squared Gaussian projections.
We develop a formulation in terms of Laplace transforms,
using the exponential-integrability estimates from
Subsection~\ref{subsec:gaussian-rotation}.
For related applications of Stein identities to
moment-generating functions and moderate deviations,
see \cite{CFS13}.

For a matrix $B\in\R^{d\times d}$, its operator norm is
$\|B\|_{\mathrm{op}}=\sup_{\|x\|=1}\|Bx\|$.
If $f\in C^1(\R^d)$ has a globally
$\varepsilon$-Lipschitz gradient, then its weak Hessian
satisfies
$\|\nabla^2f(x)\|_{\mathrm{op}}\leq\varepsilon$
almost everywhere.
Conversely, for a $C^1$ function with locally Lipschitz
gradient, this almost-everywhere Hessian bound implies
that the gradient is globally $\varepsilon$-Lipschitz.
Indeed, convolving the gradient with a smooth
nonnegative kernel of integral one preserves the bound
on its derivative. Each smoothed gradient is then
$\varepsilon$-Lipschitz by integration along line
segments. Local uniform convergence of the smoothed
gradients gives the same bound for $\nabla f$.

The approximation mechanism is explained by Gaussian
integration by parts. A centered Gaussian random
variable $Y$ with variance $v$ satisfies
\begin{equation}\label{eq:gaussian-stein-identity}
 \E[Y\psi(Y)]
 =
 v\,\E[\psi'(Y)],
\end{equation}
for example when $\psi$ is continuously differentiable
with bounded derivative, and also for exponential
test functions.
Taking $\psi(y)=\eee^{ty}$ shows that the derivative
of its Laplace transform is $vt$ times the transform
itself. Together with its value one at zero, this
determines the transform as $\eee^{vt^2/2}$.

For a centered function $f(g)$, the covariance identity
\eqref{eq:covariance} yields a similar formula in which
the constant variance is replaced by a random
coefficient. The expectation of this coefficient is
$\Var(f(g))$. A small Lipschitz constant of $\nabla f$
will make the coefficient close to its expectation,
giving an approximate version of
\eqref{eq:gaussian-stein-identity}.

\begin{proposition}[Approximate Gaussian integration by parts]
\label{prop:approximate-ibp}
Let $f\in C^1(\R^d)$ have a globally
$\varepsilon$-Lipschitz gradient, where
$\varepsilon\geq0$, and suppose that
$\E f(g)=0$ and $\E\|\nabla f(g)\|^2\leq M^2$.

For every $\psi\in C^1(\R)$ such that
$\psi(f(g))$, $\psi'(f(g))$, and
$\psi'(f(g))\nabla f(g)$ are square integrable, we have
\begin{equation}\label{eq:approximate-Stein}
 \E[f(g)\psi(f(g))]
 =
 \E[\psi'(f(g))\Gamma_f(g)],
\end{equation}
where
\begin{equation}\label{eq:Tf-definition}
 \Gamma_f(x)
 =
 \int_0^1
 \E'\left[
   \left\langle
     \nabla f(x),
     \nabla f(ux+\sqrt{1-u^2}\,g')
   \right\rangle
 \right]
 \,\dd u.
\end{equation}
Here $g'$ is an independent standard Gaussian vector,
and $\E'$ denotes expectation with respect to $g'$
only.

The random variable $\Gamma_f(g)$ is square integrable
and satisfies
\begin{align}
 \E\Gamma_f(g)
 &=
 \Var(f(g)),
 \label{eq:Tf-mean}\\
 \sqrt{\Var(\Gamma_f(g))}
 &\leq
 \frac32\varepsilon M.
 \label{eq:Tf-variance}
\end{align}
Consequently, for every test function $\psi$ as above,
\begin{equation}\label{eq:gaussian-ibp-error}
 \left|
   \E[f(g)\psi(f(g))]
   -
   \Var(f(g))\,\E[\psi'(f(g))]
 \right|
 \leq
 \frac32\varepsilon M\,
 \|\psi'(f(g))\|_{L^2}.
\end{equation}
\end{proposition}

\begin{proof}
The Lipschitz property of $\nabla f$ implies that
$\nabla f$ has at most linear growth and $f$ has at
most quadratic growth. In particular, $f$ and its
first derivatives are square integrable under
$\gamma_d$.
The assumptions on $\psi$ ensure that $\psi\circ f$
and its gradient
$\psi'(f)\nabla f$ are also square integrable.
We may therefore apply \eqref{eq:covariance} to
$\psi\circ f$ and $f$.

Recall that $g_u=ug+\sqrt{1-u^2}\,g'$.
Since $\E f(g)=0$, the covariance identity gives
\[
 \begin{aligned}
 \E[f(g)\psi(f(g))]
 &=
 \Cov(\psi(f(g)),f(g))\\
 &=
 \int_0^1
 \E\left[
   \psi'(f(g))
   \langle\nabla f(g),\nabla f(g_u)\rangle
 \right]
 \,\dd u\\
 &=
 \E[\psi'(f(g))\Gamma_f(g)].
 \end{aligned}
\]
The last equality follows by conditioning on $g$.
The interchange of expectation and integration is
justified by the Cauchy--Schwarz inequality: the
expectation of the absolute value of the integrand
is at most
$M\|\psi'(f(g))\nabla f(g)\|_{L^2}$,
uniformly in $u$.
This proves \eqref{eq:approximate-Stein}.
Taking $\psi(y)=y$ gives \eqref{eq:Tf-mean}.

To estimate the fluctuations of $\Gamma_f(g)$, write
$\Gamma_f(x)=\langle\nabla f(x),b_f(x)\rangle$, where
\[
 b_f(x)
 =
 \int_0^1
 \E'\left[
   \nabla f(ux+\sqrt{1-u^2}\,g')
 \right]
 \,\dd u.
\]
The two gradient fields in this inner product have
Lipschitz constants that can be controlled directly.
Indeed, for $x,y\in\R^d$,
\[
 \begin{aligned}
 \|b_f(x)-b_f(y)\|
 &\leq
 \int_0^1
 \E'\left[
   \left\|
     \nabla f(ux+\sqrt{1-u^2}\,g')
     -
     \nabla f(uy+\sqrt{1-u^2}\,g')
   \right\|
 \right]
 \,\dd u\\
 &\leq
 \varepsilon\int_0^1u\|x-y\|\,\dd u
 =
 \frac{\varepsilon}{2}\|x-y\|.
 \end{aligned}
\]
Thus $b_f$ is globally
$\varepsilon/2$-Lipschitz.
Jensen's inequality and the standard Gaussian
distribution of $g_u$ also give
\[
 \E\|b_f(g)\|^2
 \leq
 \int_0^1\E\|\nabla f(g_u)\|^2\,\dd u
 =
 \E\|\nabla f(g)\|^2
 \leq M^2.
\]

Both $\nabla f$ and $b_f$ have at most linear growth.
Their inner product $\Gamma_f$ is therefore locally
Lipschitz, has at most quadratic growth, and belongs
to $L^2(\R^d,\gamma_d)$.
At almost every point, both vector fields are
differentiable. The product rule and their Lipschitz
bounds imply
\[
 \|\nabla\Gamma_f(x)\|
 \leq
 \varepsilon\|b_f(x)\|
 +
 \frac{\varepsilon}{2}\|\nabla f(x)\|.
\]
These almost-everywhere derivatives represent the
weak gradient of $\Gamma_f$, as discussed in
Subsection~\ref{subsec:weak-derivatives}.
Taking $L^2$ norms yields
\[
 \|\nabla\Gamma_f(g)\|_{L^2}
 \leq
 \varepsilon\|b_f(g)\|_{L^2}
 +
 \frac{\varepsilon}{2}\|\nabla f(g)\|_{L^2}
 \leq
 \frac32\varepsilon M.
\]
Applying \eqref{eq:poincare} to $\Gamma_f$ proves
\eqref{eq:Tf-variance}.

Finally, \eqref{eq:approximate-Stein} and
\eqref{eq:Tf-mean} give
\[
 \E[f(g)\psi(f(g))]
 -
 \Var(f(g))\,\E[\psi'(f(g))]\\
 =
 \E\left[
   \psi'(f(g))
   \bigl(\Gamma_f(g)-\E\Gamma_f(g)\bigr)
 \right].
\]
The Cauchy--Schwarz inequality and
\eqref{eq:Tf-variance} now prove
\eqref{eq:gaussian-ibp-error}.
\end{proof}

For a centered affine function
$f(x)=\langle a,x\rangle$, the gradient is constant,
$b_f(x)=a$, and
$\Gamma_f(x)=\|a\|^2=\Var(f(g))$ for every $x$.
In this case, \eqref{eq:approximate-Stein} is exactly
\eqref{eq:gaussian-stein-identity}.
In general, \eqref{eq:Tf-variance} quantifies how
close the random coefficient is to the constant
variance.

We now apply Proposition~\ref{prop:approximate-ibp}
to exponential test functions.
The required exponential integrability follows
from the inequalities established in
Subsection~\ref{subsec:gaussian-rotation}.

\begin{lemma}[Gaussian Laplace approximation]
\label{lem:gaussian-laplace}
Fix $0\leq M<\infty$.
There exist constants $r>0$ and $K_M<\infty$,
depending only on $M$, with the following property.
Let $f\in C^1(\R^d)$ have a globally
$\varepsilon$-Lipschitz gradient, where
$0\leq\varepsilon\leq1$, and suppose that
\begin{equation}\label{eq:lem:gaussian-laplace-assumption}
 \E f(g)=0,
 \qquad
 \E\|\nabla f(g)\|^2\leq M^2.
\end{equation}
Then
\begin{equation}\label{eq:gaussian-exp-moment}
 \E\eee^{4r|f(g)|}\leq K_M
\end{equation}
and
\begin{equation}\label{eq:gaussian-laplace}
 \left|
   \E\eee^{tf(g)}
   -
   \eee^{t^2\Var(f(g))/2}
 \right|
 \leq
 K_M\varepsilon|t|^3,
 \qquad |t|\leq r.
\end{equation}
The constants are independent of $d$, $f$, and
$\varepsilon$. We may choose $r\leq1/2$.
\end{lemma}

The cubic factor in \eqref{eq:gaussian-laplace}
reflects the matching of the mean and variance:
the two Laplace transforms have the same value
and the same first two derivatives at zero.
In the proof, this matching will appear through
the identity \eqref{eq:Tf-mean}.

\begin{proof}
We first establish exponential integrability and
then compare the Laplace transforms.
Throughout the proof, $K_M$ denotes a finite
constant depending only on $M$; its value may
increase from one occurrence to the next.

\parhead{Step 1: Exponential integrability.}
Since $\nabla f$ is globally
$\varepsilon$-Lipschitz, the function
$x\mapsto\|\nabla f(x)\|$ is also
$\varepsilon$-Lipschitz.
Moreover,
$\E\|\nabla f(g)\|\leq M$ by the
Cauchy--Schwarz inequality.

Apply \eqref{eq:maurey-pisier} to
$x\mapsto\|\nabla f(x)\|$ with
$\Psi(u)=\eee^{cu^2}$, where $c\geq0$.
The almost-everywhere norm of the gradient of
this function is at most $\varepsilon$.
Consequently,
\[
\E\eee^{c(\|\nabla f(g)\|-\E\|\nabla f(g)\|)^2}
 \leq
 \E\eee^{\pi^2c\varepsilon^2Z^2/4}
 =
 \frac{1}{\sqrt{1-\pi^2c\varepsilon^2/2}},
\]
provided that $\pi^2c\varepsilon^2/2<1$.
Here $Z$ is a one-dimensional standard Gaussian
variable, and the last equality follows by
integrating its density.

Since $\E f(g)=0$, inequality
\eqref{eq:gaussian-exponential} gives
$\E\eee^{tf(g)}
\leq\E\eee^{\pi^2t^2\|\nabla f(g)\|^2/8}$.
Using $(a+b)^2\leq2a^2+2b^2$ and
$\E\|\nabla f(g)\|\leq M$, and then applying
the preceding estimate with $c=\pi^2t^2/4$,
we obtain
\[
 \E\eee^{tf(g)}
 \leq
 \frac{\eee^{\pi^2t^2M^2/4}}
      {\sqrt{1-\pi^4t^2\varepsilon^2/8}},
 \qquad
 \frac{\pi^4t^2\varepsilon^2}{8}<1.
\]
These applications of
\eqref{eq:maurey-pisier} and
\eqref{eq:gaussian-exponential} are valid with
possibly infinite expectations; the displayed
bound establishes their finiteness for the
indicated parameters.

Choose $r=1/(2\pi^2)$.
Since $\varepsilon\leq1$, the denominator in the
last display is at least $1/\sqrt2$ whenever
$|t|\leq4r$.
It follows that
\[
 \E\eee^{4r|f(g)|}
 \leq
 \E\eee^{4rf(g)}+\E\eee^{-4rf(g)}
 \leq
 2\sqrt2\,\eee^{M^2/\pi^2}.
\]
This proves \eqref{eq:gaussian-exp-moment}.

\parhead{Step 2: Comparison of Laplace transforms.}
Write $m(t)=\E\eee^{tf(g)}$.
We will apply Proposition~\ref{prop:approximate-ibp}
with $\psi(x)=\eee^{tx}$ for $|t|\leq r$.

We first check the integrability assumptions.
By \eqref{eq:gaussian-exp-moment},
both $\eee^{tf(g)}$ and $t\eee^{tf(g)}$
are square integrable.
Also, Corollary~\ref{lem:moment}, applied with
$q=4$ to $x\mapsto\|\nabla f(x)\|$, gives
\[
 \|\nabla f(g)\|_{L^4}
 \leq
 \E\|\nabla f(g)\|
 +
 \frac{\pi\sqrt3}{2}\varepsilon
 \leq
 M+\frac{\pi\sqrt3}{2}\varepsilon.
\]
The Cauchy--Schwarz inequality therefore yields
\[
 \E\left[
   \eee^{2tf(g)}\|\nabla f(g)\|^2
 \right]
 \leq
 \bigl(\E\eee^{4r|f(g)|}\bigr)^{1/2}
 \bigl(\E\|\nabla f(g)\|^4\bigr)^{1/2}
 <\infty.
\]
Thus the exponential test function satisfies
all the assumptions of
Proposition~\ref{prop:approximate-ibp}.

The bound \eqref{eq:gaussian-exp-moment} also
justifies differentiation of $m$ under the
expectation for $|t|\leq r$.
Indeed, $|f(g)|\eee^{2r|f(g)|}$ is integrable
and provides a dominating function for the
difference quotients on this interval.
Using \eqref{eq:approximate-Stein} and
\eqref{eq:Tf-mean}, we obtain
\begin{equation}\label{eq:laplace-differential}
 \begin{aligned}
 &m'(t)-t\,\Var(f(g))\,m(t)
 =
 t\,\E\left[
   \bigl(\Gamma_f(g)-\Var(f(g))\bigr)
   \eee^{tf(g)}
 \right]\\
 &\quad=
 t\,\E\left[
   \bigl(\Gamma_f(g)-\Var(f(g))\bigr)
   \bigl(\eee^{tf(g)}-1\bigr)
 \right].
 \end{aligned}
\end{equation}
The last equality uses
$\E\Gamma_f(g)=\Var(f(g))$.
It supplies the additional factor of $|t|$
needed for the cubic estimate.

For $x\in\R$ and $|t|\leq r$,
\[
 |\eee^{tx}-1|
 \leq
 |t|\,|x|\,\eee^{r|x|}
 \leq
 \frac{|t|}{r}\eee^{2r|x|}.
\]
Consequently,
\eqref{eq:gaussian-exp-moment} implies
$\|\eee^{tf(g)}-1\|_{L^2}\leq K_M|t|$.
Combining this bound with
\eqref{eq:Tf-variance} and the
Cauchy--Schwarz inequality in
\eqref{eq:laplace-differential}, we find that
\[
 |m'(t)-t\,\Var(f(g))\,m(t)|
 \leq K_M\varepsilon t^2,
 \qquad |t|\leq r.
\]

Multiplication by the integrating factor
$\eee^{-t^2\Var(f(g))/2}$ gives
\[
 \left|
 \frac{\dd}{\dd t}
 \left(
   \eee^{-t^2\Var(f(g))/2}m(t)
 \right)
 \right|
 \leq
 K_M\varepsilon t^2.
\]
Integrating over the interval between $0$ and $t$,
and using $m(0)=1$, yields
\[
 \left|
   \eee^{-t^2\Var(f(g))/2}m(t)-1
 \right|
 \leq
 K_M\varepsilon|t|^3.
\]
Finally, \eqref{eq:poincare} gives
$\Var(f(g))\leq M^2$.
Multiplying by
$\eee^{t^2\Var(f(g))/2}$ and enlarging $K_M$
therefore proves
\eqref{eq:gaussian-laplace}.
\end{proof}

The proof combines two estimates with different
roles.
The covariance identity and
\eqref{eq:poincare} control the fluctuations of
the integration-by-parts coefficient, while
\eqref{eq:maurey-pisier} and
\eqref{eq:gaussian-exponential} provide the
exponential integrability needed to use
exponential test functions.
Their combination gives a bound that is uniform
in the dimension.
In Section~\ref{sec:squared-projection-clt},
we apply Lemma~\ref{lem:gaussian-laplace}
to squared Gaussian projections.

\section{A central limit theorem for squared Gaussian projections}
\label{sec:squared-projection-clt}

\subsection{Metric projection and its regularity}
\label{subsec:projection-regularity}

Throughout this subsection, $C\subseteq\R^d$ is
a nonempty closed convex cone.
Recall that $\Pi_C(x)$ is the unique point of $C$
minimizing $\|x-z\|$ over $z\in C$, and that
$H(x)=\|\Pi_C(x)\|^2$.
We use the polar cone with the nonpositive convention:
$C^\circ=\{y\in\R^d:\langle y,z\rangle\leq0
\text{ for every }z\in C\}$.

We first recall Moreau's decomposition, which
characterizes the projections onto $C$ and $C^\circ$
through an orthogonal sum.
We then establish the regularity properties needed
for Gaussian analysis.
These are standard facts from convex analysis;
in particular, the decomposition and the gradient
identity \eqref{eq:H-derivatives} are recalled in
\cite[Sections~1.2 and~7.1]{MT14}.
We include proofs to make explicit that no
polyhedral assumption is required.

\begin{proposition}[Moreau decomposition]
\label{prop:moreau}
For every $x\in\R^d$,
\begin{equation}\label{eq:moreau}
  x=\Pi_C(x)+\Pi_{C^\circ}(x),
  \qquad
  \langle\Pi_C(x),\Pi_{C^\circ}(x)\rangle=0.
\end{equation}
Conversely, if $x=p+q$ with $p\in C$,
$q\in C^\circ$, and $\langle p,q\rangle=0$,
then $p=\Pi_C(x)$ and $q=\Pi_{C^\circ}(x)$.
In particular,
\begin{equation}\label{eq:projection-energy}
  \langle x,\Pi_C(x)\rangle=\|\Pi_C(x)\|^2,
  \qquad
  \|\Pi_C(x)\|\leq\|x\|.
\end{equation}
\end{proposition}

\begin{proof}
We begin with the variational characterization
of metric projection.
For a point $p\in C$, one has $p=\Pi_C(x)$
if and only if
\begin{equation}\label{eq:projection-characterization}
  \langle x-p,z-p\rangle\leq0
  \qquad\text{for every }z\in C.
\end{equation}
For necessity, note that the segment
$p+t(z-p)$ belongs to $C$ for $0\leq t\leq1$.
The function $t\mapsto\|x-p-t(z-p)\|^2$
therefore has a minimum at $t=0$ on this interval.
Its right derivative at zero is nonnegative,
which gives \eqref{eq:projection-characterization}.
For sufficiency, expand
\[
  \|x-z\|^2-\|x-p\|^2
  =
  \|z-p\|^2-2\langle x-p,z-p\rangle
  \geq0.
\]
The same characterization applies to projection
onto the closed convex cone $C^\circ$.

Taking $p=\Pi_C(x)$ and then $z=0$ and
$z=2\Pi_C(x)$ in
\eqref{eq:projection-characterization} gives
$\langle x-\Pi_C(x),\Pi_C(x)\rangle=0$.
Substituting this equality back into
\eqref{eq:projection-characterization} shows
that $\langle x-\Pi_C(x),z\rangle\leq0$
for every $z\in C$.
Thus $x-\Pi_C(x)\in C^\circ$.

It remains to identify the projection onto
$C^\circ$ and prove the converse.
Consider any decomposition $x=p+q$ with
$p\in C$, $q\in C^\circ$, and
$\langle p,q\rangle=0$.
For $z\in C$ and $w\in C^\circ$, respectively,
\begin{align*}
  \langle x-p,z-p\rangle
  =\langle q,z\rangle\leq0,
  \qquad
  \langle x-q,w-q\rangle
  =\langle p,w\rangle\leq0.
\end{align*}
By \eqref{eq:projection-characterization},
these inequalities identify $p=\Pi_C(x)$
and $q=\Pi_{C^\circ}(x)$.
Applying this observation to
$p=\Pi_C(x)$ and $q=x-\Pi_C(x)$ proves
\eqref{eq:moreau}, together with its converse.

Finally, taking the inner product of
\eqref{eq:moreau} with $\Pi_C(x)$ gives
the first identity in
\eqref{eq:projection-energy}.
Orthogonality also gives
$\|x\|^2=\|\Pi_C(x)\|^2+\|\Pi_{C^\circ}(x)\|^2$,
which proves the second assertion.
\end{proof}

The projection itself need not be differentiable
everywhere.
Its squared norm, however, is continuously
differentiable, and its gradient is Lipschitz.
The next proposition records these properties,
together with the almost-everywhere derivative
bounds that will be used later.

\begin{proposition}[Regularity of metric projection]
\label{prop:projection-properties}
The mapping $\Pi_C$ is $1$-Lipschitz.
The function $H(x)=\|\Pi_C(x)\|^2$ is convex
and continuously differentiable, with
\begin{equation}\label{eq:H-derivatives}
  \nabla H(x)=2\Pi_C(x),
  \qquad
  \|\nabla H(x)-\nabla H(y)\|\leq2\|x-y\|,
  \qquad x,y\in\R^d.
\end{equation}
Moreover, for every $x,h\in\R^d$,
\begin{equation}\label{eq:convex-Taylor}
  0
  \leq H(x+h)-H(x)-\langle\nabla H(x),h\rangle
  \leq\|h\|^2.
\end{equation}

The mapping $\Pi_C$ is differentiable almost
everywhere with respect to Lebesgue measure,
and its derivative satisfies
\begin{equation}\label{eq:projection-derivative}
  D\Pi_C(x)=D\Pi_C(x)^{\mathsf T},
  \qquad
  0\preceq D\Pi_C(x)\preceq I_d
  \quad\text{almost everywhere}.
\end{equation}
Here $B\preceq D$ means that $D-B$ is positive
semidefinite.
\end{proposition}

\begin{proof}
We first prove the Lipschitz property.
Apply \eqref{eq:projection-characterization}
at $x$ and $y$, using each projection as a
competitor for the other.
Adding the resulting inequalities yields
\[
  \|\Pi_C(x)-\Pi_C(y)\|^2
  \leq
  \langle\Pi_C(x)-\Pi_C(y),x-y\rangle.
\]
The Cauchy--Schwarz inequality therefore shows
that $\Pi_C$ is $1$-Lipschitz.

To establish convexity and differentiability
of $H$, we use a variational representation.
Completing the square and then applying
\eqref{eq:moreau}, we obtain
\begin{align*}
  \sup_{z\in C}
  \left\{
    \langle x,z\rangle-\frac12\|z\|^2
  \right\}
  =
  \frac12\|x\|^2
  -\frac12\|x-\Pi_C(x)\|^2
  =
  \frac12\|\Pi_C(x)\|^2
  =\frac12H(x).
\end{align*}
The supremum is attained uniquely at $\Pi_C(x)$.
Since the expression inside the supremum is
affine as a function of $x$, this representation
shows that $H$ is convex.

Evaluating the supremum for $x+h$ at $\Pi_C(x)$
gives a lower bound for $H(x+h)-H(x)$.
Evaluating the supremum for $x$ at $\Pi_C(x+h)$
gives the corresponding upper bound.
Together with the Lipschitz property just proved,
these bounds yield
\begin{align*}
  0
  \leq
  \frac12\bigl(H(x+h)-H(x)\bigr)
  -\langle\Pi_C(x),h\rangle
  \leq
  \langle\Pi_C(x+h)-\Pi_C(x),h\rangle
  \leq\|h\|^2.
\end{align*}
Consequently, $H$ is differentiable and
$\nabla H(x)=2\Pi_C(x)$.
The continuity and $2$-Lipschitz property of
its gradient follow from the Lipschitz property
of $\Pi_C$.
This proves \eqref{eq:H-derivatives}.

The lower bound in \eqref{eq:convex-Taylor}
follows from convexity.
For the upper bound, integrate the gradient
along the segment from $x$ to $x+h$ and use
\eqref{eq:H-derivatives}:
\begin{align*}
  H(x+h)-H(x)-\langle\nabla H(x),h\rangle
  =
  \int_0^1
    \langle\nabla H(x+th)-\nabla H(x),h\rangle
  \,\dd t
  \leq
  \int_0^1 2t\|h\|^2\,\dd t
  =\|h\|^2.
\end{align*}

Finally, Rademacher's theorem implies that
$\Pi_C$ is differentiable almost everywhere.
Its weak derivative agrees almost everywhere
with this derivative.
Since $\Pi_C=\nabla(H/2)$, the commutation of
mixed weak derivatives shows that $D\Pi_C$
is symmetric almost everywhere.

Convexity of $H$ implies the monotonicity of
its gradient, so
$\langle\Pi_C(x+th)-\Pi_C(x),th\rangle\geq0$.
At a point where $\Pi_C$ is differentiable,
divide this inequality by $t^2$ and let $t\to0$.
We obtain
$\langle D\Pi_C(x)h,h\rangle\geq0$
for every $h\in\R^d$.
Hence $D\Pi_C(x)$ is positive semidefinite
almost everywhere.
Recall that the operator norm of a matrix $B$ is
$\|B\|_{\op}=\sup_{\|h\|=1}\|Bh\|$.
Since $\Pi_C$ is $1$-Lipschitz, we have
$\|D\Pi_C(x)\|_{\op}\leq1$
wherever the derivative exists.
Together with symmetry and positive
semidefiniteness, this proves
\eqref{eq:projection-derivative}.
\end{proof}

In particular, $\nabla^2H(x)=2D\Pi_C(x)$
almost everywhere.
Since the standard Gaussian distribution has
a density, all the almost-everywhere assertions
above also hold almost surely when evaluated
at a standard Gaussian vector.
We shall use \eqref{eq:projection-derivative}
without assuming that $D\Pi_C(x)$ is an
orthogonal projection.

As a first application, we obtain the following
variance bound, which is contained in
\cite[Theorem~4.5]{MT14}.
Its proof combines \eqref{eq:H-derivatives}
with \eqref{eq:poincare} and the moment
identities \eqref{eq:H-moments}.

\begin{lemma}[Universal variance bound]
\label{lem:universal-variance}
Let $C\subseteq\R^d$ be a nonempty closed convex cone,
and let $V_C$ be its conic intrinsic-volume
random variable.
Writing $\delta=\E V_C$ and
$\tau^2=\Var(V_C)$, we have
\begin{equation}\label{eq:universal-variance}
  \tau^2\leq2\delta.
\end{equation}
\end{lemma}

\begin{proof}
Let $g$ be a standard Gaussian random vector
in $\R^d$.
By \eqref{eq:projection-energy} and
\eqref{eq:H-derivatives},
$0\leq H(x)\leq\|x\|^2$ and
$\|\nabla H(x)\|\leq2\|x\|$.
Thus $H$ and its gradient are square integrable
under the Gaussian law, so
\eqref{eq:poincare} applies.
Using \eqref{eq:H-moments} and
\eqref{eq:H-derivatives}, we obtain
\begin{align*}
  2\delta+\tau^2
  =\Var(H(g))
  \leq\E\|\nabla H(g)\|^2
  =4\E H(g)
  =4\delta.
\end{align*}
Subtracting $2\delta$ proves
\eqref{eq:universal-variance}.
\end{proof}

Together, \eqref{eq:H-moments} and
\eqref{eq:universal-variance} give
$2\delta\leq\Var(H(g))\leq4\delta$.
In the next subsection, we use
\eqref{eq:H-derivatives} and
\eqref{eq:gaussian-laplace} to prove a central
limit theorem for the centered and normalized
squared projection norm.

\subsection{A central limit theorem for squared Gaussian projections}
\label{subsec:squared-projection-clt}

By \eqref{eq:H-derivatives}, the gradient of
$x\mapsto\|\Pi_C(x)\|^2$ has Lipschitz constant
at most two, independently of the cone and the
ambient dimension.
On the other hand, \eqref{eq:H-moments} shows
that its variance under the Gaussian law tends
to infinity whenever $\delta$ tends to infinity.
After centering and dividing by the standard
deviation, we therefore obtain a function whose
gradient has a small Lipschitz constant.
We shall verify that the expected squared norm
of its gradient remains bounded, so that
Lemma~\ref{lem:gaussian-laplace} applies.

The following result was established by
Goldstein, Nourdin, and Peccati
\cite[Theorem~1.1, part~1]{GNP17}.
We give a proof using \eqref{eq:gaussian-laplace}.

\begin{proposition}[Central limit theorem for squared Gaussian projections]
\label{prop:squared-projection-clt}
Let $C_n\subseteq\R^{d_n}$ be a sequence of
nonempty closed convex cones, and let $g_n$
be a standard Gaussian random vector in
$\R^{d_n}$.
Write $\delta_n=\E V_{C_n}$ and
$\tau_n^2=\Var(V_{C_n})$.
If $\delta_n\to\infty$, then
\begin{equation}\label{eq:squared-projection-clt}
  \frac{
    \|\Pi_{C_n}(g_n)\|^2-\delta_n
  }{
    \sqrt{2\delta_n+\tau_n^2}
  }
  \toD\normal(0,1).
\end{equation}
\end{proposition}

\begin{proof}
Since $\delta_n\to\infty$, we may restrict
attention to indices for which $\delta_n\geq2$.
To apply Lemma~\ref{lem:gaussian-laplace},
consider the centered and normalized function
\[
  f_n(x)
  =
  \frac{
    \|\Pi_{C_n}(x)\|^2-\delta_n
  }{
    \sqrt{2\delta_n+\tau_n^2}
  },
  \qquad x\in\R^{d_n}.
\]
By \eqref{eq:H-moments},
$\E f_n(g_n)=0$ and $\Var(f_n(g_n))=1$.

Equation~\eqref{eq:H-derivatives} shows that
$f_n$ is continuously differentiable and gives
its gradient explicitly:
\[
  \nabla f_n(x)
  =
  \frac{2\Pi_{C_n}(x)}
       {\sqrt{2\delta_n+\tau_n^2}}.
\]
Consequently, using \eqref{eq:H-moments} again,
\begin{align*}
  \E\|\nabla f_n(g_n)\|^2
  =
  \frac{
    4\E\|\Pi_{C_n}(g_n)\|^2
  }{
    2\delta_n+\tau_n^2
  }
  =
  \frac{4\delta_n}{2\delta_n+\tau_n^2}
  \leq2.
\end{align*}
Thus the expected squared norm of the gradient
is bounded uniformly in $n$.

Moreover, $\Pi_{C_n}$ is $1$-Lipschitz by
Proposition~\ref{prop:projection-properties}.
Hence
\[
  \operatorname{Lip}(\nabla f_n)
  \leq
  \frac{2}{\sqrt{2\delta_n+\tau_n^2}}
  \leq
  \sqrt{\frac{2}{\delta_n}}
  \longrightarrow0.
\]
Because $\delta_n\geq2$, the displayed bound
is at most one.
The hypotheses of
Lemma~\ref{lem:gaussian-laplace} therefore hold
with $M=\sqrt2$ and
$\varepsilon=2/\sqrt{2\delta_n+\tau_n^2}$.

Since $\Var(f_n(g_n))=1$,
equation~\eqref{eq:gaussian-laplace} yields
universal constants $r>0$ and $K<\infty$
such that
\begin{equation}\label{eq:squared-projection-laplace}
  \left|
    \E\eee^{t f_n(g_n)}-\eee^{t^2/2}
  \right|
  \leq
  \frac{2K|t|^3}{\sqrt{2\delta_n+\tau_n^2}}
  \leq
  \frac{\sqrt2K|t|^3}{\sqrt{\delta_n}},
  \qquad |t|\leq r.
\end{equation}
In particular, the moment-generating functions
converge to that of a standard Gaussian
random variable on a fixed neighborhood of zero.
Lemma~\ref{lem:laplace-continuity}, stated
and proved in the appendix, now gives
\eqref{eq:squared-projection-clt}.
\end{proof}

The conclusion concerns the squared projection
norm, rather than the intrinsic-volume random
variable $V_{C_n}$.
To prove Theorem~\ref{thm:clt}, we will need
additional information about the Wiener-chaos
expansion of the squared projection norm.
We develop the necessary tools in the next
section and apply them in Section~\ref{sec:removal}.

\section{Hermite expansions and Gaussian chaoses}
\label{sec:hermite}

The Gaussian approximation developed in Section~\ref{sec:gaussian}
and applied to squared Gaussian projections in
Section~\ref{sec:squared-projection-clt} relies on
interpolation between Gaussian vectors.
We now introduce a different tool: the orthogonal
decomposition of functions on Gaussian space into
Hermite polynomials, also called the Wiener-chaos
expansion.

Our main purpose is to identify the quadratic
component of a function of a Gaussian vector and
to establish a refinement of the Gaussian Poincar\'e
inequality when this component vanishes.
Both results will be used in the proof of the
central limit theorem for conic intrinsic volumes.

\subsection{Hermite polynomials and Gaussian chaos expansions}
\label{subsec:hermite}

Let $\gamma_d$ denote the standard Gaussian
probability measure on $\R^d$, and let $g$ be a
random vector with distribution $\gamma_d$.
We write $L^2(\R^d,\gamma_d)$ for the Hilbert
space of functions $f$ satisfying
$\E[f(g)^2]<\infty$, with functions identified
when they agree almost everywhere with respect
to $\gamma_d$.

The probabilists' Hermite polynomials
$\operatorname{He}_0(x)=1$,
$\operatorname{He}_1(x)=x$,
$\operatorname{He}_2(x)=x^2-1$, and so on,
are defined by the generating function
\[
  \eee^{tx-t^2/2}
  =
  \sum_{k=0}^\infty
    \operatorname{He}_k(x)\frac{t^k}{k!}.
\]
For a standard Gaussian random variable $Z$,
they satisfy
\[
  \E[\operatorname{He}_k(Z)\operatorname{He}_\ell(Z)]
  =
  k!\ind_{\{k=\ell\}},
  \qquad
  \operatorname{He}_k'(x)
  =
  k\operatorname{He}_{k-1}(x).
\]

For a multi-index
$\nu=(\nu_1,\ldots,\nu_d)\in\mathbb N_0^d$,
write $|\nu|=\nu_1+\cdots+\nu_d$ and
$\nu!=\nu_1!\cdots\nu_d!$.
The normalized products
\[
  \mathcal H_\nu(x)
  =
  \frac{
    \prod_{i=1}^d\operatorname{He}_{\nu_i}(x_i)
  }{\sqrt{\nu!}}
\]
form an orthonormal basis of
$L^2(\R^d,\gamma_d)$.

For completeness, let us recall why this system
is complete.
The Hermite polynomials span the same finite-degree
polynomial spaces as ordinary monomials.
If $f\in L^2(\R^d,\gamma_d)$ is orthogonal to
every polynomial, all derivatives at zero of
\[
  z\longmapsto
  \E\left[f(g)\eee^{\langle z,g\rangle}\right],
  \qquad z\in\mathbb C^d,
\]
vanish.
This function is entire, by Cauchy--Schwarz
and Gaussian exponential integrability,
so it vanishes identically.
Its restriction to purely imaginary arguments
is the Fourier transform of the finite signed
measure $f\,\gamma_d$.
Uniqueness of the Fourier transform implies
that $f=0$ almost everywhere, proving completeness.

Consequently, every $f\in L^2(\R^d,\gamma_d)$
admits a unique orthogonal expansion
\begin{equation}\label{eq:hermite-expansion}
  f=\sum_{q=0}^\infty f_q,
  \qquad
  f_q(x)
  =
  \sum_{|\nu|=q}
    \E[f(g)\mathcal H_\nu(g)]\mathcal H_\nu(x).
\end{equation}
The series converges in $L^2(\R^d,\gamma_d)$,
not necessarily pointwise.
The function $f_q$ is called the Gaussian-chaos
component, or Wiener-chaos component, of order $q$.
In particular, $f_0=\E f(g)$ is constant,
$f_1$ is linear, and $f_2$ is a centered
quadratic polynomial.

The orthogonality of the components immediately
gives a decomposition of the variance.
Similar identities hold for the expected squared
norms of the gradient and the Hessian.

For a matrix $B\in\R^{d\times d}$, we denote
its Hilbert--Schmidt norm by $\|B\|_{\HS}$,
where
\[
  \|B\|_{\HS}^2
  =
  \tr(B^{\mathsf T}B)
  =
  \sum_{i,j=1}^d B_{ij}^2.
\]
\begin{proposition}[Hermite identities]
\label{prop:hermite-identities}
Let $f\in L^2(\R^d,\gamma_d)$, and suppose
that all its first and second weak partial
derivatives belong to $L^2(\R^d,\gamma_d)$.
Let $f=\sum_{q\geq0}f_q$ be its expansion
in \eqref{eq:hermite-expansion}.
Then
\begin{align}
  \Var(f(g))
  &=
  \sum_{q\geq1}\E[f_q(g)^2],
  \label{eq:hermite-variance}\\
  \E\|\nabla f(g)\|^2
  &=
  \sum_{q\geq1}q\,\E[f_q(g)^2],
  \label{eq:hermite-gradient}\\
  \E\|\nabla^2f(g)\|_{\HS}^2
  &=
  \sum_{q\geq2}q(q-1)\,\E[f_q(g)^2].
  \label{eq:hermite-hessian}
\end{align}
\end{proposition}

\begin{proof}
The first identity follows directly from
\eqref{eq:hermite-expansion}, orthogonality,
and Parseval's identity.
Indeed, $f-f_0$ is the orthogonal sum
of the components of orders $q\geq1$.

We next prove \eqref{eq:hermite-gradient}.
Let $e_i$ denote the $i$th coordinate unit vector
in $\mathbb N_0^d$.
The Hermite recurrence
\[
  x\operatorname{He}_k(x)
  -\operatorname{He}_k'(x)
  =
  \operatorname{He}_{k+1}(x)
\]
gives
\[
  x_i\mathcal H_\nu(x)
  -\partial_i\mathcal H_\nu(x)
  =
  \sqrt{\nu_i+1}\,\mathcal H_{\nu+e_i}(x).
\]
Applying \eqref{eq:ibp-polynomial} with
$P=\mathcal H_\nu$ therefore yields
\begin{equation}\label{eq:hermite-derivative-coefficients}
  \E[\partial_i f(g)\mathcal H_\nu(g)]
  =
  \sqrt{\nu_i+1}\,
  \E[f(g)\mathcal H_{\nu+e_i}(g)].
\end{equation}

Applying Parseval's identity to each
$\partial_i f$ and summing over $i$ now gives
\begin{align*}
  \E\|\nabla f(g)\|^2
  =
  \sum_{i=1}^d\sum_{\nu\in\mathbb N_0^d}
    (\nu_i+1)
    \left(
      \E[f(g)\mathcal H_{\nu+e_i}(g)]
    \right)^2
  =
  \sum_{q\geq1}q\,\E[f_q(g)^2].
\end{align*}
In the last equality, each coefficient
corresponding to a multi-index of total
degree $q$ occurs with total weight $q$.
This proves \eqref{eq:hermite-gradient}.

Finally, apply
\eqref{eq:hermite-derivative-coefficients}
twice to the second weak derivatives of $f$.
Parseval's identity, summed over all ordered
pairs of indices, shows that a component
of degree $q$ contributes with weight
$q(q-1)$ to the expected squared
Hilbert--Schmidt norm of the Hessian.
This proves \eqref{eq:hermite-hessian}.

The argument uses only Gaussian integration
by parts and Parseval's identity.
In particular, it does not require
term-by-term differentiation of the
infinite series in \eqref{eq:hermite-expansion}.
\end{proof}

\parhead{Another proof of Gaussian Poincar\'e's inequality.}
The Hermite identities give another proof of the
Gaussian Poincar\'e inequality \eqref{eq:poincare}.
For this approach in dimension one, see
\cite{Chernoff81}.
Indeed, for a function satisfying the assumptions
of Proposition~\ref{prop:hermite-identities},
\[
  \Var(f(g))
  =
  \sum_{q\geq1}\E[f_q(g)^2]
  \leq
  \sum_{q\geq1}q\,\E[f_q(g)^2]
  =
  \E\|\nabla f(g)\|^2.
\]
The assumption on the second weak derivatives is
not needed here, since the proof of
\eqref{eq:hermite-gradient} requires only that
$f$ and its first weak derivatives belong to
$L^2(\R^d,\gamma_d)$.

Thus the Gaussian Poincar\'e inequality
reflects the elementary fact that every
nonconstant Gaussian-chaos component has
order at least one.

More importantly for our purposes,
\eqref{eq:hermite-gradient} and
\eqref{eq:hermite-hessian} allow us to
obtain a stronger estimate when one of
the chaos components vanishes.

\subsection{The first three Gaussian chaoses}
\label{subsec:second-chaos}
The zeroth Gaussian chaos consists of constant functions,
and the first Gaussian chaos consists of linear functions
$x\mapsto\langle b,x\rangle$, where $b\in\R^d$.
The second Gaussian chaos consists of the centered
quadratic polynomials $x^{\mathsf T}Bx-\tr B$,
where $B$ ranges over all symmetric $d\times d$ matrices.
Equivalently, it is spanned by $x_i^2-1$, $1\leq i\leq d$,
and $x_ix_j$, $1\leq i<j\leq d$.

The next lemma identifies the first three components
of a Gaussian-chaos expansion in terms of the
expectations of the function, its gradient, and its Hessian.

\begin{lemma}[Projection onto the first three Gaussian chaoses]
\label{lem:second-chaos-projection}
Let $f\in L^2(\R^d,\gamma_d)$, and suppose
that its first and second weak partial derivatives
belong to $L^2(\R^d,\gamma_d)$.
Then the zeroth and first Gaussian-chaos components of $f$
are given by
\begin{equation}\label{eq:zeroth-first-chaos-formulas}
  f_0=\E[f(g)],
  \qquad
  f_1(x)=\langle\E[\nabla f(g)],x\rangle,
\end{equation}
and its second Gaussian-chaos component is given by
\begin{equation}\label{eq:second-chaos-formula}
  f_2(x)
  =
  \frac12\left(
    x^{\mathsf T}\E[\nabla^2f(g)]x
    -\tr\bigl(\E[\nabla^2f(g)]\bigr)
  \right).
\end{equation}
In particular, $f_1=0$ if and only if
$\E[\nabla f(g)]=0$, and $f_2=0$ if and only if
$\E[\nabla^2f(g)]=0$.
\end{lemma}

\begin{proof}
The formula for $f_0$ follows immediately from
the fact that the zeroth Gaussian chaos consists
of constant functions.

The coordinate functions $x_1,\ldots,x_d$ form
an orthonormal basis of the first Gaussian chaos.
By Gaussian integration by parts \eqref{eq:ibp-first},
the coefficient of $x_i$ in the orthogonal projection
of $f$ onto the first Gaussian chaos is
$\E[g_if(g)]=\E[\partial_i f(g)]$.
This proves \eqref{eq:zeroth-first-chaos-formulas}.

The polynomials $x_i^2-1$, $1\leq i\leq d$,
and $x_ix_j$, $1\leq i<j\leq d$, form an
orthogonal basis of the second Gaussian chaos.
Their squared $L^2(\R^d,\gamma_d)$ norms
are $2$ and $1$, respectively.

By Gaussian integration by parts
\eqref{eq:ibp-second}, the coefficient of
$x_i^2-1$ in the orthogonal projection is
$\E[\partial_i^2f(g)]/2$.
Similarly, the coefficient of $x_ix_j$
for $i<j$ is $\E[\partial_i\partial_jf(g)]$.
Summing these terms and using the symmetry
of the weak Hessian gives
\eqref{eq:second-chaos-formula}.

The final assertions follow from the uniqueness
of the coefficients in the respective orthogonal bases.
\end{proof}

\subsection{A refined Gaussian Poincar\'e inequality}
The absence of certain components in a Gaussian-chaos
expansion allows us to derive additional inequalities
relating the variance of a function to the squared
$L^2$-norms of its gradient and Hessian.
The ordinary Poincar\'e inequality bounds the variance
of $f(g)$ by the expected squared norm of its gradient.
Conversely, the expected squared norm of the gradient
can be bounded by the variance of $f(g)$ and a correction
term involving the Hessian.
In the next section, we shall need a refined version
of this reverse inequality for functions whose second
Gaussian-chaos component vanishes.

\begin{lemma}[A refined Gaussian Poincar\'e inequality]
\label{lem:hermite-gap}
Let $f\in L^2(\R^d,\gamma_d)$, and suppose
that its first and second weak partial derivatives
belong to $L^2(\R^d,\gamma_d)$.
If $\E[\nabla^2f(g)]=0$, then
\begin{equation}\label{eq:hermite-gap}
  \E\|\nabla f(g)\|^2
  \leq
  \Var(f(g))
  +\frac13\E\|\nabla^2f(g)\|_{\HS}^2.
\end{equation}
\end{lemma}

\begin{proof}
Let $f=\sum_{q\geq0}f_q$ be the orthogonal
Gaussian-chaos expansion of $f$.
By Lemma~\ref{lem:second-chaos-projection},
the assumption $\E[\nabla^2f(g)]=0$
is equivalent to $f_2=0$.

Using the identities
\eqref{eq:hermite-variance}--\eqref{eq:hermite-hessian},
the desired inequality becomes
\[
  \sum_{q\geq1}q\,\E[f_q(g)^2]
  \leq
  \sum_{q\geq1}\E[f_q(g)^2]
  +\frac13\sum_{q\geq2}
    q(q-1)\,\E[f_q(g)^2].
\]
Since $f_2=0$, it suffices to verify
the scalar inequality
\[
  q\leq1+\frac{q(q-1)}3
  \qquad\text{for }q=1\text{ and all }q\geq3.
\]
Indeed, the difference between the right-hand
side and the left-hand side equals
$(q-1)(q-3)/3$, which is nonnegative
for these values of $q$.

Multiplying the scalar inequality by
$\E[f_q(g)^2]$ and summing over all
nonvanishing chaos components proves
\eqref{eq:hermite-gap}.
\end{proof}

The absence of the second chaos is essential
for \eqref{eq:hermite-gap}.
Indeed, the scalar inequality used in the proof
fails for $q=2$.
On the other hand, the first chaos is allowed:
the inequality holds with equality for $q=1$.
Thus the lemma provides a stronger estimate
without requiring the linear Gaussian component
to vanish.

In the next section, we apply
Lemma~\ref{lem:second-chaos-projection}
to identify the low-order Gaussian-chaos components
of the projection map $\Pi_C$ and of the squared
projection norm $H(x)=\|\Pi_C(x)\|^2$.
We investigate the geometric meaning of these components
and their relationship with the intrinsic-volume variance.
After removing the second Gaussian-chaos component
of $H$, we use \eqref{eq:hermite-gap} to control
the variance and derivative energies of the
remaining function.

\section{Wiener chaoses and the geometry of convex cones}
\label{sec:removal}

Section~\ref{sec:squared-projection-clt} establishes a central limit
theorem for the squared Gaussian projection $\norm{\Pi_C(g)}^2$.
To prove the corresponding theorem for the intrinsic-volume random
variable $V_C$ under the sole assumption that its variance tends to
infinity, we need a more precise description of these Gaussian
fluctuations. In this section, we study the Wiener-chaos expansions
of both $\Pi_C(g)$ and $\norm{\Pi_C(g)}^2$, using the results of
Section~\ref{sec:hermite}.

The first coefficients of these expansions have direct geometric
interpretations. They describe the mean Gaussian projection and
the best linear approximation to the projection map, and they
lead to a distance between cones defined by their Gaussian
projections. We shall use these interpretations to relate
$\Var(V_C)$ to approximation by linear subspaces. We shall also
show that, after centering the squared projection and removing
its second Gaussian-chaos component, the variance of the
remaining function is comparable to $\Var(V_C)$. The estimates
obtained here will be used in the next section to prove the
central limit theorem for $V_C$.

\subsection{From the squared-projection CLT to the
intrinsic-volume CLT}
\label{subsec:intrinsic-volume-clt-difficulty}

Let $C_n\subset\R^{d_n}$ be a sequence of closed convex cones,
and write $V_n=V_{C_n}$, $\delta_n=\E V_n$, and
$\tau_n^2=\Var(V_n)$. Suppose that $\tau_n\to\infty$.
Our eventual goal is to prove that
$(V_n-\delta_n)/\tau_n\toD\normal(0,1)$.
By Lemma~\ref{lem:laplace-continuity}, it is enough to establish
\begin{equation}\label{eq:laplace-clt-goal}
 \E\exp\left(t\frac{V_n-\delta_n}{\tau_n}\right)
 \longrightarrow \eee^{t^2/2}
 \qquad\text{for every } |t|\leq t_0,
\end{equation}
where $t_0>0$ is independent of $n$.

Let $g_n$ be a standard Gaussian vector in $\R^{d_n}$, and put
$H_n(x)=\norm{\Pi_{C_n}(x)}^2$.
By \eqref{eq:H-moments}, we have $\E H_n(g_n)=\delta_n$ and
$\Var(H_n(g_n))=2\delta_n+\tau_n^2$.
Moreover, Lemma~\ref{lem:universal-variance} gives
$\tau_n^2\leq 2\delta_n$, so our assumption implies
$\delta_n\to\infty$. Proposition~\ref{prop:squared-projection-clt}
therefore applies to this sequence.

It is tempting to combine that proposition directly with
\eqref{eq:steiner-laplace}. The difficulty is that the two
random variables are normalized differently.
The estimate \eqref{eq:squared-projection-laplace} controls
the Laplace transform of
\[
 \frac{H_n(g_n)-\delta_n}
      {\sqrt{2\delta_n+\tau_n^2}}
\]
for parameters in a fixed neighborhood of zero.
On the other hand, substituting $s=t/\tau_n$ into
\eqref{eq:steiner-laplace} gives
\begin{equation}\label{eq:steiner-at-standardized-parameter}
 \E\exp\left(t\frac{V_n-\delta_n}{\tau_n}\right)
 =
 \eee^{-t\delta_n/\tau_n}
 \E\exp\left(
   \frac{1-\eee^{-2t/\tau_n}}{2}\,H_n(g_n)
 \right).
\end{equation}
After centering $H_n(g_n)$, the expectation on the right
requires the Laplace transform of the normalized squared
projection at the parameter
\begin{equation}\label{eq:projection-laplace-parameter}
 \frac{\sqrt{2\delta_n+\tau_n^2}}{2}
       \bigl(1-\eee^{-2t/\tau_n}\bigr)
 \sim
 t\sqrt{1+\frac{2\delta_n}{\tau_n^2}},
 \qquad n\to\infty,
\end{equation}
for every fixed $t\ne0$.

The ratio $\delta_n/\tau_n^2$ need not be bounded.
Consequently, the absolute value of the parameter in
\eqref{eq:projection-laplace-parameter} may tend to infinity,
even when $t$ is fixed and arbitrarily small.
Thus \eqref{eq:squared-projection-laplace} does not directly
provide the estimate required in
\eqref{eq:steiner-at-standardized-parameter} under the sole
assumption $\tau_n\to\infty$.

Linear directions explain why this mismatch can occur.
If $L$ is a $k$-dimensional linear subspace and $D$ is another
closed convex cone, then \eqref{eq:polarity-product} gives
$V_{L\times D}\eqD k+V_D$.
Taking the Cartesian product with $L$ therefore increases the
statistical dimension by $k$ and leaves the intrinsic-volume
variance unchanged. By adding sufficiently many linear
directions, we can make the ratio of the statistical dimension
to the intrinsic-volume variance arbitrarily large.

The same observation also points to the Gaussian-chaos
component that should be removed.

\begin{example}[A linear subspace]
\label{ex:subspace-chaos-motivation}
Let $L\subset\R^d$ be a linear subspace of dimension $k$.
Choose orthonormal coordinates in which $L$ is spanned by
the first $k$ coordinate vectors.
For a standard Gaussian vector $g=(g_1,\ldots,g_d)$,
\[
 \norm{\Pi_L(g)}^2-k
 =
 \sum_{j=1}^k(g_j^2-1).
\]
Each summand belongs to the second Gaussian chaos.
Thus all the fluctuations of the squared Gaussian projection
belong to the second Gaussian chaos, whereas $V_L=k$ almost
surely and has no fluctuations.
\end{example}

This example identifies a Gaussian component that can be
arbitrarily large even when the intrinsic-volume random
variable is deterministic. For a general cone $C$, we shall
therefore center $H(g)=\norm{\Pi_C(g)}^2$ and subtract its
second Gaussian-chaos component.

The first Gaussian-chaos component will be retained.
As we shall see, it can account for almost all of the
variance that remains after the quadratic component has
been removed. We shall identify these components explicitly
and prove that the resulting remainder has variance
comparable, up to universal constants, to $\Var(V_C)$.
This comparison will allow the estimates used in the next
section to depend on the intrinsic-volume variance, without
requiring a bound on the ratio
$\E V_C/\Var(V_C)$.

\subsection{Wiener-chaos expansions of Gaussian projections}
\label{subsec:quadratic-component}

Throughout this subsection, let $C\subset\R^d$ be a nonempty
closed convex cone, let $g$ be a standard Gaussian vector in
$\R^d$, and write $\delta=\E V_C$ and
$H(x)=\norm{\Pi_C(x)}^2$.

We regard $H$ as an element of $L^2(\R^d,\gamma_d)$ and
$\Pi_C$ as an element of $L^2(\R^d,\gamma_d;\R^d)$.
In the latter space, the squared norm of a function $u$ is
$\E\norm{u(g)}^2$.
The bound $\norm{\Pi_C(x)}\leq\norm{x}$ ensures that both
functions belong to the indicated spaces.
We denote their components in the $q$-th Gaussian chaos by
$H_q$ and $(\Pi_C)_q$, respectively.
For the vector-valued function $\Pi_C$, the chaos projection
is taken coordinate by coordinate.

The first components of these expansions are determined by
the mean Gaussian projection and the average derivative of
the projection map. Accordingly, we write
\begin{equation}\label{eq:A-definition}
 m=\E\Pi_C(g),
 \qquad
 A=\E[D\Pi_C(g)].
\end{equation}
The derivative in \eqref{eq:A-definition} exists almost
everywhere and is bounded, by
Proposition~\ref{prop:projection-properties}, so its expectation
is well defined. When several cones are considered, we write
$m_C$ and $A_C$ to indicate the dependence on $C$.

\begin{proposition}[Gaussian-chaos components of the projection]
\label{prop:H-chaos}
The projection map and its squared norm have the expansions
\begin{align}
 \Pi_C(x)
 &=
 m+Ax+\sum_{q=2}^{\infty}(\Pi_C)_q(x),
 \label{eq:projection-chaos-expansion}\\
 H(x)
 &=
 \delta+2\ip{m}{x}
 +\bigl(x^{\mathsf T}Ax-\delta\bigr)
 +\sum_{q=3}^{\infty}H_q(x).
 \label{eq:H-chaos-expansion}
\end{align}
The series converge in
$L^2(\R^d,\gamma_d;\R^d)$ and
$L^2(\R^d,\gamma_d)$, respectively.
Their components are related by
\begin{equation}\label{eq:projection-chaos-derivative}
 (\Pi_C)_q(x)=\frac12\nabla H_{q+1}(x),
 \qquad q\geq0.
\end{equation}

The matrix $A$ is symmetric and satisfies
\begin{equation}\label{eq:A-trace}
 A=\E[\Pi_C(g)g^{\mathsf T}],
 \qquad
 0\preceq A\preceq I_d,
 \qquad
 \tr A=\delta.
\end{equation}
In particular, the first and second Gaussian-chaos
components of $H$ have variances
\begin{equation}\label{eq:H-low-chaos-variances}
 \Var(H_1(g))=4\norm{m}^2,
 \qquad
 \Var(H_2(g))=2\tr(A^2).
\end{equation}
\end{proposition}

\begin{proof}
By \eqref{eq:projection-derivative}, the matrix $D\Pi_C(x)$
is symmetric and satisfies
$0\preceq D\Pi_C(x)\preceq I_d$ almost everywhere.
Taking expectations gives the corresponding properties of $A$.

Each coordinate of $\Pi_C$, together with its first weak
derivatives, belongs to $L^2(\R^d,\gamma_d)$.
Applying \eqref{eq:ibp-first} coordinate by coordinate
therefore gives $A=\E[\Pi_C(g)g^{\mathsf T}]$.
Taking the trace and using
\eqref{eq:projection-energy} and \eqref{eq:H-moments}, we obtain
\[
 \tr A
 =
 \E\ip{\Pi_C(g)}{g}
 =
 \E\norm{\Pi_C(g)}^2
 =
 \delta.
\]
This proves \eqref{eq:A-trace}.

By \eqref{eq:H-derivatives}, we have
$\nabla H=2\Pi_C$ and
$D^2H=2D\Pi_C$ almost everywhere.
These identities and the bounds on the projection and its
derivative show that $H$, its first derivatives, and its
second weak derivatives belong to $L^2(\R^d,\gamma_d)$.
Lemma~\ref{lem:second-chaos-projection} therefore gives
$H_0=\delta$, $H_1(x)=2\ip{m}{x}$, and
$H_2(x)=x^{\mathsf T}Ax-\tr A$.
Since $\tr A=\delta$, this proves
\eqref{eq:H-chaos-expansion}.

To relate the two expansions, we use
\eqref{eq:hermite-derivative-coefficients}.
For every $q\geq0$, that identity shows that the $q$-th
Gaussian-chaos component of $\partial_iH$ is
$\partial_iH_{q+1}$.
Indeed, a Hermite coefficient of total degree $q$ in
$\partial_iH$ is determined by the corresponding coefficient
of total degree $q+1$ in $H$.
Since $\nabla H=2\Pi_C$, comparison of the components of
order $q$ gives \eqref{eq:projection-chaos-derivative}.
Taking $q=0$ and $q=1$ yields
$(\Pi_C)_0=m$ and $(\Pi_C)_1(x)=Ax$,
which proves \eqref{eq:projection-chaos-expansion}.
The identity \eqref{eq:projection-chaos-derivative} holds
for the polynomial representatives of the individual
chaos components.

The first identity in \eqref{eq:H-low-chaos-variances}
follows from $H_1(g)=2\ip{m}{g}$.
For the second, let $\lambda_1,\ldots,\lambda_d$ be the
eigenvalues of $A$.
By orthogonal invariance of the Gaussian distribution,
$H_2(g)$ has the same distribution as
$\sum_{i=1}^d\lambda_i(g_i^2-1)$.
The summands are independent and have variances
$2\lambda_i^2$, so
$\Var(H_2(g))=2\sum_{i=1}^d\lambda_i^2=2\tr(A^2)$.
\end{proof}

The shift in the chaos order in
\eqref{eq:projection-chaos-derivative} is useful to keep
in mind. The vector $m$ determines the constant component
of the projection and the linear component of its squared
norm. The matrix $A$ determines the linear component of
the projection and the quadratic component of its squared
norm.

We next record how these coefficients behave under
polarity, orthogonal transformations, and Cartesian products.

\begin{proposition}[Transformation rules for the chaos coefficients]
\label{prop:chaos-transformations}
Let $Q$ be an orthogonal transformation of $\R^d$, and let
$D\subset\R^k$ be another nonempty closed convex cone.
Then
\begin{align}
 m_{C^\circ}&=-m_C,
 & A_{C^\circ}&=I_d-A_C,
 \label{eq:chaos-polarity}\\
 m_{QC}&=Qm_C,
 & A_{QC}&=QA_CQ^{\mathsf T},
 \label{eq:chaos-rotations}\\
 m_{C\times D}&=m_C\oplus m_D,
 & A_{C\times D}&=A_C\oplus A_D.
 \label{eq:chaos-products}
\end{align}
Here $\oplus$ denotes concatenation of vectors and the
block diagonal sum of matrices.

In particular, if $C$ is self-dual in the sense that
$C^\circ=-C$, then
\begin{equation}\label{eq:self-dual-chaos}
 A_C=\frac12I_d,
 \qquad
 \E V_C=\frac d2,
 \qquad
 H_2(x)=\frac12\bigl(\norm{x}^2-d\bigr).
\end{equation}
\end{proposition}

\begin{proof}
Moreau's decomposition \eqref{eq:moreau} gives
$\Pi_{C^\circ}(x)=x-\Pi_C(x)$.
Taking expectations at $g$ yields
$m_{C^\circ}=-m_C$.
Taking derivatives almost everywhere and then expectations
gives $A_{C^\circ}=I_d-A_C$.
This proves \eqref{eq:chaos-polarity}.

For an orthogonal transformation $Q$, the defining
minimization problem for the metric projection gives
$\Pi_{QC}(x)=Q\Pi_C(Q^{\mathsf T}x)$.
Since $Q^{\mathsf T}g$ has the same distribution as $g$,
taking expectations in this identity and in its derivative
gives \eqref{eq:chaos-rotations}.

Similarly, projection onto a Cartesian product separates
into its two coordinates:
$\Pi_{C\times D}(x,y)=\Pi_C(x)\oplus\Pi_D(y)$.
Taking means and average derivatives gives
\eqref{eq:chaos-products}.

Finally, setting $Q=-I_d$ in
\eqref{eq:chaos-rotations} gives $A_{-C}=A_C$.
If $C^\circ=-C$, then
\eqref{eq:chaos-polarity} consequently implies
$A_C=I_d-A_C$, and hence $A_C=\frac12I_d$.
The remaining identities in \eqref{eq:self-dual-chaos}
follow from \eqref{eq:A-trace} and
\eqref{eq:H-chaos-expansion}.
\end{proof}

\subsection{Geometric interpretations of the chaos coefficients}
\label{subsec:chaos-geometric-interpretations}

We retain the notation
$H(x)=\norm{\Pi_C(x)}^2$,
$m=\E\Pi_C(g)$,
$A=\E[D\Pi_C(g)]$, and
$\delta=\E V_C=\E H(g)$,
where $C\subset\R^d$ is a nonempty closed convex cone and
$d\geq1$.
Proposition~\ref{prop:H-chaos} identifies $m$ and $A$ as
coefficients in the Gaussian-chaos expansions of the projection
and its squared norm. We now examine their geometric meaning.

We begin by changing the mean and covariance of the Gaussian
vector. These perturbations give interpretations of $m$ and $A$
as derivatives of the expected squared projection norm.

\begin{proposition}[Translations and covariance perturbations]
\label{prop:projection-gaussian-perturbations}
As $x\to0$, we have
\begin{equation}\label{eq:projection-translation-expansion}
 \E H(g+x)
 =
 \delta+2\ip{m}{x}
 +x^{\mathsf T}Ax
 +o(\norm{x}^2).
\end{equation}
Moreover, for every symmetric matrix $B\in\R^{d\times d}$,
\begin{equation}\label{eq:projection-covariance-derivative}
 \left.
 \frac{\dd}{\dd t}
 \E H\bigl((I_d+tB)^{1/2}g\bigr)
 \right|_{t=0}
 =
 \tr(AB).
\end{equation}
The expression in
\eqref{eq:projection-covariance-derivative} is defined for all
sufficiently small real $t$, for which $I_d+tB$ is positive
definite.
\end{proposition}

\begin{proof}
The identity $\nabla H=2\Pi_C$ and the linear growth of
$\Pi_C$ allow us to differentiate the expectation with respect
to the translation. Since $\Pi_C$ is Lipschitz, its difference
quotients are bounded, and a second differentiation gives
\[
 \nabla_x\E H(g+x)=2\E\Pi_C(g+x),
 \qquad
 D_x^2\E H(g+x)=2\E[D\Pi_C(g+x)].
\]
For the second identity, differentiability of $\Pi_C$ almost
everywhere and dominated convergence apply at every fixed $x$.

The matrix-valued function
$x\mapsto\E[D\Pi_C(g+x)]$ is continuous.
Indeed, $D\Pi_C$ is bounded, and its expectation can be written
as an integral against a translated Gaussian density.
Translations of that density are continuous in $L^1$.
Thus $x\mapsto\E H(g+x)$ is twice continuously differentiable.
Its value, gradient, and Hessian at zero are
$\delta$, $2m$, and $2A$, respectively.
Taylor's theorem proves
\eqref{eq:projection-translation-expansion}.

For the covariance perturbation, diagonalization of $B$ gives
$(I_d+tB)^{1/2}=I_d+\frac12tB+O(t^2)$ as $t\to0$.
The difference quotients of
$H((I_d+tB)^{1/2}g)$ are dominated, for sufficiently small
$|t|$, by a constant multiple of $\norm{g}^2$.
This follows from the linear growth of $\nabla H$.
We may therefore differentiate under the expectation and obtain
\[
 \left.
 \frac{\dd}{\dd t}
 \E H\bigl((I_d+tB)^{1/2}g\bigr)
 \right|_{t=0}
 =
 \E\ip{\Pi_C(g)}{Bg}
 =
 \tr(AB),
\]
where the last equality follows from
\eqref{eq:A-trace}.
\end{proof}

For a unit vector $u$, the covariance matrix
$I_d+tuu^{\mathsf T}$ adds Gaussian variance in direction $u$
when $t>0$.
Equation~\eqref{eq:projection-covariance-derivative} shows that
the first-order change in the expected squared projection norm
is $\ip{Au}{u}$.
By \eqref{eq:A-trace}, this quantity lies between zero and one.

The endpoints of this interval have precise geometric meanings.
Recall that $C\cap(-C)$ is the lineality space of $C$, namely
the largest linear subspace contained in $C$.

\begin{proposition}[The mean projection and the extreme eigenspaces]
\label{prop:mean-projection-and-kernels}
The mean Gaussian projection satisfies
$m\in\relint(C)\cap\relint(-C^\circ)$, where the relative
interiors are taken in the respective linear spans.
The matrix $A$ satisfies
\begin{equation}\label{eq:projection-matrix-kernels}
 \ker A=(\spanop C)^\perp,
 \qquad
 \ker(I_d-A)=C\cap(-C).
\end{equation}
\end{proposition}

\begin{proof}
The support of the distribution of $\Pi_C(g)$ is exactly $C$.
Indeed, the projection takes values in $C$.
Conversely, if $y\in C$, then $\Pi_C(y)=y$, and continuity of
the projection implies that every neighborhood of $y$ receives
positive probability under the distribution of $\Pi_C(g)$.
Here we use that the Gaussian distribution assigns positive
probability to every nonempty open ball.

Closedness and convexity give $m\in C$.
Suppose that $m\notin\relint(C)$.
A supporting hyperplane in $\spanop C$ then gives a nonzero
vector $u\in\spanop C$ such that
$\ip{u}{z}\leq0$ for every $z\in C$ and
$\ip{u}{m}=0$.
The supporting level is zero because $C$ is a cone.
The nonpositive random variable
$\ip{u}{\Pi_C(g)}$ has expectation zero and therefore
vanishes almost surely.
Since the support of $\Pi_C(g)$ is $C$, it follows that
$\ip{u}{z}=0$ for every $z\in C$.
This contradicts $u\in\spanop C$ and $u\ne0$.
Thus $m\in\relint(C)$.
Applying the same argument to $C^\circ$ and using
$m_{C^\circ}=-m_C$ from \eqref{eq:chaos-polarity}
gives $m\in\relint(-C^\circ)$.

We next identify $\ker A$.
For any $u\in\R^d$, the Gaussian Poincar\'e inequality
\eqref{eq:poincare}, symmetry of $D\Pi_C$, and
\eqref{eq:projection-derivative} give
\[
 \Var\bigl(\ip{u}{\Pi_C(g)}\bigr)
 \leq
 \E\norm{D\Pi_C(g)u}^2
 \leq
 \E\ip{D\Pi_C(g)u}{u}
 =
 \ip{Au}{u}.
\]
If $Au=0$, the continuous function
$x\mapsto\ip{u}{\Pi_C(x)}$ is constant almost surely under
Gaussian measure.
Since Gaussian measure has full support, it is constant
everywhere. Its value at zero is zero.
Evaluating at $x\in C$, where $\Pi_C(x)=x$, shows that
$u\in(\spanop C)^\perp$.

Conversely, if $u\in(\spanop C)^\perp$, then
$\ip{u}{\Pi_C(g)}=0$ almost surely.
The identity $A=\E[\Pi_C(g)g^{\mathsf T}]$ gives
$u^{\mathsf T}A=0$, and symmetry of $A$ yields $Au=0$.
This proves the first identity in
\eqref{eq:projection-matrix-kernels}.

Applying that identity to $C^\circ$ and using
$A_{C^\circ}=I_d-A_C$, we obtain
\[
 \ker(I_d-A)
 =
 (\spanop C^\circ)^\perp
 =
 C\cap(-C).
\]
For the last equality, every vector in $C\cap(-C)$ is
orthogonal to $C^\circ$ by the definition of the polar cone.
Conversely, if $u$ is orthogonal to $C^\circ$, then
Moreau's decomposition gives
$\Pi_{C^\circ}(u)=\Pi_{C^\circ}(-u)=0$,
and hence both $u$ and $-u$ belong to $C$.
\end{proof}

Thus an eigenvalue zero corresponds to a direction orthogonal
to the cone, while an eigenvalue one corresponds to a direction
in its lineality space.
The following corollary describes the remaining part of the
cone and records two characterizations of linear subspaces.

\begin{corollary}[The orthogonal decomposition of the cone]
\label{cor:cone-orthogonal-decomposition}
The following conditions are equivalent:
\begin{equation}\label{eq:subspace-chaos-characterizations}
 C\text{ is a linear subspace}
 \quad\Longleftrightarrow\quad
 m=0
 \quad\Longleftrightarrow\quad
 A^2=A.
\end{equation}

More generally, let $L=\ker(I_d-A)$,
$N=\ker A$, and $E=(L\oplus N)^\perp$.
Then $C=L\oplus K$, where $K=C\cap L^\perp$ is a pointed
cone spanning $E$.
Here pointed means that $K\cap(-K)=\{0\}$.
Relative to the orthogonal decomposition
$\R^d=L\oplus E\oplus N$, we have
\begin{equation}\label{eq:cone-orthogonal-decomposition}
 A=
 \begin{pmatrix}
  I_L&0&0\\
  0&A_K&0\\
  0&0&0_N
 \end{pmatrix},
 \qquad
 m\in E.
\end{equation}
The matrix $A_K$ is computed by regarding $K$ as a cone in
the Euclidean space $E$ and using a standard Gaussian vector
in that space.
All eigenvalues of $A_K$ lie strictly between zero and one.
In particular,
\begin{equation}\label{eq:projection-matrix-ranks}
 \begin{aligned}
  \dim(\spanop C)&=\operatorname{rank}A,\\
  \dim(C\cap(-C))&=\dim\ker(I_d-A),\\
  \dim E&=\operatorname{rank}(A-A^2).
 \end{aligned}
\end{equation}
\end{corollary}

The characterization of linear subspaces by $m=0$ was also
noted by Goldstein, Nourdin, and Peccati
\cite[Remark~4.1(b)]{GNP17}.

\begin{proof}
By \eqref{eq:projection-matrix-kernels}, $L$ is the
lineality space of $C$ and $N=(\spanop C)^\perp$.
For every $x\in C$, both $\Pi_L(x)$ and $-\Pi_L(x)$
belong to $L\subset C$.
Since a convex cone is closed under addition,
$x-\Pi_L(x)\in C\cap L^\perp=K$.
This proves $C=L\oplus K$.

If $x\in K\cap(-K)$, then
$x\in C\cap(-C)=L$ and $x\in L^\perp$, so $x=0$.
Thus $K$ is pointed.
Moreover,
$\spanop C=N^\perp=L\oplus E$.
Taking linear spans in $C=L\oplus K$ gives
$\spanop K=E$.

For $x_L\in L$, $x_E\in E$, and $x_N\in N$, the projection
separates as
$\Pi_C(x_L+x_E+x_N)=x_L+\Pi_K(x_E)$.
Taking means and average derivatives proves
\eqref{eq:cone-orthogonal-decomposition}.
Since $E$ is orthogonal to the zero and unit eigenspaces
of $A$, every eigenvalue of $A_K$ lies in $(0,1)$.
The dimension formulas
\eqref{eq:projection-matrix-ranks} follow.
They also show that $A^2=A$ if and only if $E=\{0\}$,
which is equivalent to $C=L$ being a linear subspace.

Finally, if $m=0$, then
Proposition~\ref{prop:mean-projection-and-kernels} gives
$0\in\relint(C)$.
A relative neighborhood of zero in $\spanop C$ is therefore
contained in $C$.
Positive homogeneity then implies $C=\spanop C$,
so $C$ is a linear subspace.
Conversely, projection onto a linear subspace is a linear map
and has Gaussian mean zero.
This completes the proof of
\eqref{eq:subspace-chaos-characterizations}.
\end{proof}

In particular, the first Gaussian-chaos component of $H$
vanishes precisely when $C$ is a linear subspace.
For every other cone, $m$ supplies a distinguished direction
in the span of its pointed factor.
The second Gaussian-chaos matrix also detects linear
subspaces, through the condition $A^2=A$.

The trace identity $\tr A=\delta$ gives another interpretation
of the second chaos. It separates a component depending only
on the Gaussian radius from a component associated with the
traceless part of $A$.

\begin{remark}[A gap for the norm of the mean projection]
\label{rem:mean-projection-gap}
If $C$ is not a linear subspace, then
$\norm{m}\geq1/\sqrt{2\pi}$.
Indeed, by Proposition~\ref{prop:mean-projection-and-kernels}
and Corollary~\ref{cor:cone-orthogonal-decomposition},
the vector $e=m/\norm{m}$ belongs to
$C\cap(-C^\circ)$.
Consequently, Moreau's decomposition gives
$\ip{e}{\Pi_C(x)}\geq\max\{0,\ip{e}{x}\}$
for every $x\in\R^d$.
Taking Gaussian expectations proves the asserted bound.
Equality is attained, for example, by a ray.
\end{remark}

\begin{proposition}[The radial part of the second Gaussian chaos]
\label{prop:radial-quadratic-chaos}
The second Gaussian-chaos component of $H$ has the orthogonal
decomposition
\begin{equation}\label{eq:radial-quadratic-chaos}
 H_2(g)
 =
 \frac{\delta}{d}\bigl(\norm{g}^2-d\bigr)
 +
 g^{\mathsf T}
 \left(A-\frac{\delta}{d}I_d\right)g.
\end{equation}
The variances of the two terms are, respectively,
$2\delta^2/d$ and
$2\norm{A-\frac{\delta}{d}I_d}_{\HS}^2$.
Here the Hilbert--Schmidt norm is given by
$\norm{B}_{\HS}^2=\sum_{i,j}B_{ij}^2$.

Moreover,
\begin{equation}\label{eq:radial-conditional-expectation}
 \E\bigl[H(g)-\delta\mid\norm{g}\bigr]
 =
 \frac{\delta}{d}\bigl(\norm{g}^2-d\bigr).
\end{equation}
Thus the orthogonal projection of $H(g)-\delta$ onto the
space of square-integrable functions of $\norm{g}$ belongs
entirely to the second Gaussian chaos.
\end{proposition}

\begin{proof}
Write $g=RU$, where $R=\norm{g}$ and $U$ is uniform on the
unit sphere and independent of $R$.
Positive homogeneity of the projection gives
$H(RU)=R^2H(U)$.
Since $\E R^2=d$ and $\E H(g)=\delta$, we have
$\E H(U)=\delta/d$.
Consequently,
$\E[H(g)\mid R]=(\delta/d)R^2$,
which proves \eqref{eq:radial-conditional-expectation}.

The decomposition
\eqref{eq:radial-quadratic-chaos} follows from
$H_2(g)=g^{\mathsf T}Ag-\delta$.
To prove orthogonality, use
$\E[UU^{\mathsf T}]=I_d/d$ and $\tr A=\delta$ to obtain
\[
 \E\left[
  g^{\mathsf T}
  \left(A-\frac{\delta}{d}I_d\right)g
  \,\middle|\,R
 \right]
 =
 \frac{R^2}{d}
 \tr\left(A-\frac{\delta}{d}I_d\right)
 =
 0.
\]
The second term in
\eqref{eq:radial-quadratic-chaos} is therefore orthogonal
to every square-integrable function of $R$.

Since $\Var(R^2)=2d$, the first term has variance
$2\delta^2/d$.
Using orthogonality and
\eqref{eq:H-low-chaos-variances}, the variance of the
second term is
\[
 2\tr(A^2)-\frac{2\delta^2}{d}
 =
 2\norm{A-\frac{\delta}{d}I_d}_{\HS}^2.
\]
\end{proof}

For a self-dual cone, \eqref{eq:self-dual-chaos} gives
$A=I_d/2$ and $\delta=d/2$, so the traceless term in
\eqref{eq:radial-quadratic-chaos} vanishes.
For general cones, that term can be substantial.
For example, if $C=L$ is a linear subspace of dimension
$0<k<d$, then $A=\Pi_L$ and $\delta=k$.
After subtracting the radial term, the remaining quadratic
form has variance $2k(1-k/d)>0$, although $V_L$ is
deterministic.
This explains why the subtraction used below removes the
entire second Gaussian chaos.

Finally, the vector-valued expansion
\eqref{eq:projection-chaos-expansion} suggests approximating
the projection map by its constant and linear components.
The following identity identifies the optimal approximations
and their errors.

\begin{proposition}[Optimal linear and affine approximations]
\label{prop:projection-least-squares}
For every matrix $B\in\R^{d\times d}$ and vector $v\in\R^d$,
\begin{equation}\label{eq:affine-least-squares}
  \E\norm{\Pi_C(g)-Bg-v}^2
  =
  \tr(A-A^2)-\norm{m}^2
  +\norm{A-B}_{\HS}^2+\norm{m-v}^2.
\end{equation}
The unique optimal linear approximation to $\Pi_C$ in
Gaussian mean square is $x\mapsto Ax$.
The unique optimal affine approximation is $x\mapsto Ax+m$.

We denote the squared error of the optimal linear
approximation by
\begin{equation}\label{eq:a-definition}
 a
 =
 \E\norm{\Pi_C(g)-Ag}^2
 =
 \tr(A-A^2).
\end{equation}
Then, for every $B\in\R^{d\times d}$,
\begin{equation}\label{eq:least-squares}
 \E\norm{\Pi_C(g)-Bg}^2
 =
 a+\norm{A-B}_{\HS}^2,
\end{equation}
while the squared error of the optimal affine approximation
is $a-\norm{m}^2$.
In particular, $0\leq\norm{m}^2\leq a\leq d/4$.
\end{proposition}

\begin{proof}
Write
\[
 \Pi_C(g)-Bg-v
 =
 \bigl(\Pi_C(g)-Ag-m\bigr)
 +(A-B)g+(m-v).
\]
By \eqref{eq:projection-chaos-expansion}, the three terms
on the right are pairwise orthogonal in
$L^2(\R^d,\gamma_d;\R^d)$.
The same orthogonality gives
\[
 \begin{aligned}
  \E\norm{\Pi_C(g)-Ag-m}^2
  &=
  \E\norm{\Pi_C(g)}^2-\E\norm{Ag}^2-\norm{m}^2\\
  &=
  \delta-\tr(A^2)-\norm{m}^2\\
  &=
  \tr(A-A^2)-\norm{m}^2.
 \end{aligned}
\]
Taking the expected squared norm in the preceding
decomposition proves
\eqref{eq:affine-least-squares}.
Setting $v=0$ gives
\eqref{eq:a-definition} and \eqref{eq:least-squares}.
The nonnegative squared terms identify the unique minimizers.

Nonnegativity of the optimal affine approximation error
gives $\norm{m}^2\leq a$.
Finally, every eigenvalue of $A$ lies in $[0,1]$ and
contributes at most $1/4$ to $\tr(A-A^2)$.
Hence $a\leq d/4$.
\end{proof}

The vector $m$, the matrix $A$, and the quantity $a$
describe geometric and approximation properties of the
projection map.
In the next subsection, we establish their quantitative
relationships with the intrinsic-volume variance.

\subsection{The intrinsic-volume variance and the
Gaussian-chaos expansion}
\label{subsec:remainder-variance}

We retain the notation of the preceding subsections and write
$\tau^2=\Var(V_C)$.
In particular, $H(x)=\norm{\Pi_C(x)}^2$,
$m=\E\Pi_C(g)$, $A=\E[D\Pi_C(g)]$, and
$a=\tr(A-A^2)$.
We begin with an exact identity relating $\tau^2$ to the
Gaussian-chaos components of $H$.

\begin{proposition}[A signed representation of the
intrinsic-volume variance]
\label{prop:intrinsic-variance-chaoses}
The intrinsic-volume variance satisfies
\begin{align}
 \tau^2
 &=
 \Var(H(g))
 -\frac12\E\norm{\nabla H(g)}^2
 \nonumber\\
 &=
 \sum_{q=1}^{\infty}
 \left(1-\frac q2\right)\E[H_q(g)^2]
 \label{eq:tau-chaos-representation}\\
 &=
 2\norm{m}^2
 -\frac12\sum_{q=3}^{\infty}
 (q-2)\E[H_q(g)^2].
 \label{eq:tau-first-chaos}
\end{align}
The series in these identities converge absolutely.
\end{proposition}

\begin{proof}
By \eqref{eq:H-moments},
$\Var(H(g))=2\delta+\tau^2$.
The identity $\nabla H=2\Pi_C$ gives
$\E\norm{\nabla H(g)}^2=4\delta$.
Subtracting half of the latter quantity from the former
proves the first equality.

Equations~\eqref{eq:hermite-variance} and
\eqref{eq:hermite-gradient} then give
\eqref{eq:tau-chaos-representation}.
Absolute convergence follows from
$\sum_{q\geq1}q\E[H_q(g)^2]=4\delta<\infty$.
Finally, the coefficient of the second chaos is zero,
and \eqref{eq:H-low-chaos-variances} gives
$\E[H_1(g)^2]=4\norm{m}^2$.
This proves \eqref{eq:tau-first-chaos}.
\end{proof}

In \eqref{eq:tau-chaos-representation}, the first Gaussian
chaos has a positive coefficient, the second has coefficient
zero, and all higher chaoses have negative coefficients.
These signs describe an identity for $\tau^2$ in terms of
the chaos expansion of $H$; they do not constitute an
orthogonal variance decomposition of $V_C$.
In particular, \eqref{eq:tau-first-chaos} immediately gives
$\tau^2\leq2\norm{m}^2$.
The spectral argument underlying
\eqref{eq:tau-chaos-representation} is closely related
to the proof of the improved Gaussian Poincar\'e
inequality in \cite[Theorem~A.2]{GNP17}.
Goldstein, Nourdin, and Peccati applied that inequality
to the squared Gaussian projection norm to obtain
the upper bound $\tau^2\leq2\norm{m}^2$
\cite[Theorem~4.1]{GNP17}.
Here we retain the exact contributions of the
individual Gaussian-chaos components.

The vanishing coefficient of the second chaos explains why
we remove that component. We also remove the constant
component, so that the remaining function is centered.

\begin{definition}[The remainder after removing the
quadratic component]
\label{def:quadratic-chaos-remainder}
We define
\begin{equation}\label{eq:F-definition}
 F(x)
 =
 H(x)-\delta-H_2(x)
 =
 H(x)-x^{\mathsf T}Ax.
\end{equation}
\end{definition}

By \eqref{eq:H-chaos-expansion}, the Gaussian-chaos expansion
of $F$ is
\begin{equation}\label{eq:F-chaos-expansion}
 F(x)
 =
 H_1(x)+\sum_{q=3}^{\infty}H_q(x),
\end{equation}
with convergence in $L^2(\R^d,\gamma_d)$.
Thus $F_0=F_2=0$, while $F_1=H_1$ and $F_q=H_q$ for
$q\geq3$.
In particular, $\E F(g)=0$.
Since the zeroth and second components do not contribute
to \eqref{eq:tau-chaos-representation}, we also have
\begin{equation}\label{eq:tau-remainder-identity}
 \tau^2
 =
 \Var(F(g))
 -\frac12\E\norm{\nabla F(g)}^2.
\end{equation}

The following estimates show that the variance and gradient
energy of $F$ are controlled by $\tau^2$.
We recall that the operator norm of a matrix $B$ is
$\norm{B}_{\op}=\sup_{\norm{u}=1}\norm{Bu}$.

\begin{proposition}[Variance and derivative estimates
for the remainder]
\label{prop:removal}
The function $F$ is continuously differentiable, and its
gradient is $2$-Lipschitz.
Its weak Hessian satisfies
\begin{equation}\label{eq:F-hessian-op}
 \E[D^2F(g)]=0,
 \qquad
 \norm{D^2F(x)}_{\op}\leq2
 \quad\text{for almost every }x.
\end{equation}
Moreover,
\begin{align}
 \Var(F(g))
 &=\tau^2+2a,
 \label{eq:F-variance}\\
 \E\norm{\nabla F(g)}^2
 &=4a,
 \label{eq:F-energy}\\
 \E\norm{D^2F(g)}_{\HS}^2
 &\leq4a.
 \label{eq:F-hessian-HS}
\end{align}
Consequently,
\begin{equation}\label{eq:a-comparison}
 \frac12\tau^2\leq a\leq\frac32\tau^2,
 \qquad
 2\tau^2\leq\Var(F(g))\leq4\tau^2.
\end{equation}
In particular, $\tau^2\leq2a\leq d/2$.
\end{proposition}

\begin{proof}
The differentiability properties of $H$ give
\begin{equation}\label{eq:F-derivatives}
 \nabla F(x)=2\bigl(\Pi_C(x)-Ax\bigr),
 \qquad
 D^2F(x)=2\bigl(D\Pi_C(x)-A\bigr)
 \quad\text{almost everywhere}.
\end{equation}
Since $\E[D\Pi_C(g)]=A$, the expected Hessian of $F$ is zero.

Both $D\Pi_C(x)$ and $A$ are symmetric matrices between
zero and $I_d$.
It follows that
$-I_d\preceq D\Pi_C(x)-A\preceq I_d$
almost everywhere, which proves the operator-norm bound in
\eqref{eq:F-hessian-op}.
This bound also implies that $\nabla F$ is $2$-Lipschitz.
For completeness, one may convolve $\nabla F$ with a smooth
nonnegative approximate identity.
The resulting smooth maps have derivatives of operator norm
at most two and are therefore $2$-Lipschitz.
Passing to their locally uniform limit gives the same
bound for $\nabla F$.

The components of distinct Gaussian chaoses are orthogonal.
Hence \eqref{eq:F-chaos-expansion},
\eqref{eq:H-moments}, and
\eqref{eq:H-low-chaos-variances} give
\[
 \Var(F(g))
 =
 \Var(H(g))-\Var(H_2(g))
 =
 2\delta+\tau^2-2\tr(A^2)
 =
 \tau^2+2a.
\]
The gradient identity in \eqref{eq:F-derivatives} and
\eqref{eq:a-definition} give
$\E\norm{\nabla F(g)}^2=4a$.

For the Hessian, symmetry and
$\E[D\Pi_C(g)]=A$ yield
\[
 \begin{aligned}
  \E\norm{D^2F(g)}_{\HS}^2
  &=
  4\E\tr\bigl((D\Pi_C(g)-A)^2\bigr)\\
  &=
  4\E\tr\bigl(D\Pi_C(g)^2\bigr)-4\tr(A^2)\\
  &\leq
  4\tr(A-A^2)
  =
  4a.
 \end{aligned}
\]
The inequality follows from
$D\Pi_C(g)^2\preceq D\Pi_C(g)$.

The Gaussian Poincar\'e inequality \eqref{eq:poincare}
now gives $\tau^2+2a\leq4a$, and hence
$a\geq\tau^2/2$.
Since the second Gaussian-chaos component of $F$ vanishes,
the refined inequality \eqref{eq:hermite-gap} applies and
gives
\[
 4a
 \leq
 \Var(F(g))
 +\frac13\E\norm{D^2F(g)}_{\HS}^2
 \leq
 \tau^2+2a+\frac43a.
\]
Thus $a\leq3\tau^2/2$.
Substitution into \eqref{eq:F-variance} proves the remaining
bounds in \eqref{eq:a-comparison}.
Finally, $a\leq d/4$ by
Proposition~\ref{prop:projection-least-squares}, so
$\tau^2\leq2a\leq d/2$.
\end{proof}

These estimates hold without any restriction on the ratio
$\delta/\tau^2$.
They therefore give the variance control that was missing
from the direct comparison in
Subsection~\ref{subsec:intrinsic-volume-clt-difficulty}.

We next quantify the contribution of the first Gaussian
chaos to the variance that remains after the quadratic
component has been removed.

\begin{proposition}[Universal bounds for the first Gaussian chaos]
\label{prop:first-chaos-bounds}
For every nonempty closed convex cone,
\begin{equation}\label{eq:mean-projection-variance-comparison}
 \frac43\norm{m}^2
 \leq
 \tau^2
 \leq
 2\norm{m}^2.
\end{equation}
Equivalently,
$2\tau^2\leq\Var(H_1(g))\leq3\tau^2$.
Moreover,
\begin{equation}\label{eq:first-chaos-variance-fraction}
 \Var(H_1(g))
 \geq
 \frac34\Var(F(g)).
\end{equation}
\end{proposition}

\begin{proof}
By \eqref{eq:H-derivatives} and
\eqref{eq:projection-derivative},
\[
 \E\norm{D^2H(g)}_{\HS}^2
 =
 4\E\tr\bigl(D\Pi_C(g)^2\bigr)
 \leq
 4\tr A
 =
 4\delta.
\]
The Hermite identities therefore give
\begin{equation}\label{eq:higher-chaos-weighted-bound}
 \begin{aligned}
  &\sum_{q=3}^{\infty}
   (q-1)(q-2)\E[H_q(g)^2]=
   \E\norm{D^2H(g)}_{\HS}^2
   -2\E\norm{\nabla H(g)}^2
   +2\Var(H(g))\\
  &\qquad=
   \E\norm{D^2H(g)}_{\HS}^2
   -4\delta+2\tau^2
   \leq2\tau^2.
 \end{aligned}
\end{equation}
Since $q-1\geq2$ for $q\geq3$, it follows that
\begin{equation}\label{eq:higher-chaos-variance-bound}
 \sum_{q=3}^{\infty}(q-2)\E[H_q(g)^2]
 \leq\tau^2.
\end{equation}
Combining this with \eqref{eq:tau-first-chaos} gives
$2\tau^2\leq4\norm{m}^2\leq3\tau^2$,
which is equivalent to
\eqref{eq:mean-projection-variance-comparison}.

Also, because $q-2\geq1$ for $q\geq3$,
\eqref{eq:higher-chaos-variance-bound} implies
$\sum_{q\geq3}\E[H_q(g)^2]\leq\tau^2$.
Using \eqref{eq:tau-first-chaos} once more, we obtain
\[
  \E[H_1(g)^2]
  =
  2\tau^2+
  \sum_{q=3}^{\infty}(q-2)\E[H_q(g)^2] \geq
  3\sum_{q=3}^{\infty}\E[H_q(g)^2].
\]
Finally, \eqref{eq:F-chaos-expansion} gives
$\Var(F(g))=\E[H_1(g)^2]
+\sum_{q\geq3}\E[H_q(g)^2]$.
This proves \eqref{eq:first-chaos-variance-fraction}.
\end{proof}

Thus at least three quarters of the variance left after
the quadratic subtraction belongs to the first Gaussian
chaos. This gives a quantitative reason for retaining the
first component.

The preceding estimates also give two interpretations of
$\tau^2$ in terms of the matrix $A$.

\begin{corollary}[Linear approximation and the second Gaussian chaos]
\label{cor:geometric-variance}
The error of the optimal linear approximation satisfies
\begin{equation}\label{eq:geometric-comparison}
 \frac12\tau^2
 \leq
 \min_{B\in\R^{d\times d}}
 \E\norm{\Pi_C(g)-Bg}^2
 \leq
 \frac32\tau^2.
\end{equation}
Furthermore,
\begin{equation}\label{eq:second-chaos-variance-deficit}
 \tau^2
 \leq
 2\delta-\Var(H_2(g))
 =
 2a
 \leq
 3\tau^2.
\end{equation}
\end{corollary}

\begin{proof}
The minimum in \eqref{eq:geometric-comparison} equals $a$
by \eqref{eq:least-squares}, so the first assertion follows
from \eqref{eq:a-comparison}.
For the second, use
$\Var(H_2(g))=2\tr(A^2)$ and $\tr A=\delta$ to obtain
$2\delta-\Var(H_2(g))=2\tr(A-A^2)=2a$,
and apply \eqref{eq:a-comparison} again.
\end{proof}

Although the second chaos has coefficient zero in
\eqref{eq:tau-chaos-representation}, its matrix still
contains information about the intrinsic-volume variance.
Equation~\eqref{eq:second-chaos-variance-deficit} shows that
$\tau^2$ is comparable to the deficit of
$\Var(H_2(g))$ from $2\delta$.
Together with
\eqref{eq:mean-projection-variance-comparison} and
\eqref{eq:geometric-comparison}, this relates $\tau^2$
to the mean Gaussian projection, the optimal linear
approximation error, and the variance of the quadratic
Gaussian component.

We now characterize when one of the first two chaoses
accounts asymptotically for all the relevant variance.

\begin{proposition}[Asymptotic dominance of the first
and second Gaussian chaoses]
\label{prop:chaos-dominance}
Let $C_n\subset\R^{d_n}$ be a sequence of nonempty closed
convex cones, and let $g_n$ be a standard Gaussian vector
in $\R^{d_n}$.
Use subscripts $n$ for the quantities associated with $C_n$,
and denote the $q$-th Gaussian-chaos component of $H_n$
by $H_{n,q}$.

If $\delta_n>0$ for every $n$, then
\begin{equation}\label{eq:second-chaos-dominance}
 \frac{\Var(H_{n,2}(g_n))}
      {\Var(H_n(g_n))}
 \longrightarrow1
 \quad\Longleftrightarrow\quad
 \frac{\tau_n^2}{\delta_n}\longrightarrow0.
\end{equation}

If none of the cones $C_n$ is a linear subspace, then
\begin{equation}\label{eq:first-chaos-dominance}
  \frac{\Var(H_{n,1}(g_n))}
       {\Var(F_n(g_n))}
  \longrightarrow1\quad
  \Longleftrightarrow\quad
  \frac{\norm{m_n}^2}{a_n}\longrightarrow1\quad
  \Longleftrightarrow\quad
  \frac{a_n}{\tau_n^2}\longrightarrow\frac12.
\end{equation}
These equivalences do not require $\tau_n^2\to\infty$.
\end{proposition}

\begin{proof}
For the second-chaos statement, orthogonality gives
\[
 1-\frac{\Var(H_2(g))}{\Var(H(g))}
 =
 \frac{\Var(F(g))}{2\delta+\tau^2}
 =
 \frac{\tau^2+2a}{2\delta+\tau^2}.
\]
By \eqref{eq:a-comparison}, this quantity lies between
$2\tau^2/(2\delta+\tau^2)$ and
$4\tau^2/(2\delta+\tau^2)$.
For $\delta>0$, these bounds prove
\eqref{eq:second-chaos-dominance}.

We next prove \eqref{eq:first-chaos-dominance}.
For a cone that is not a linear subspace, we have $a>0$
by \eqref{eq:a-definition} and
\eqref{eq:subspace-chaos-characterizations}.
Equation~\eqref{eq:a-comparison} then gives $\tau^2>0$,
so all the ratios under consideration are defined.

The Hermite identities, together with
\eqref{eq:F-variance}, \eqref{eq:F-energy}, and
\eqref{eq:F-chaos-expansion}, give
\begin{align}
 2a-\tau^2
 &=
 \sum_{q=3}^{\infty}(q-1)\E[H_q(g)^2],
 \label{eq:remainder-energy-variance-gap}\\
 4\bigl(a-\norm{m}^2\bigr)
 &=
 \sum_{q=3}^{\infty}q\E[H_q(g)^2].
 \label{eq:affine-error-chaoses}
\end{align}
Since $q-1\leq q\leq\frac32(q-1)$ for $q\geq3$,
\begin{equation}\label{eq:affine-error-variance-gap}
 2a-\tau^2
 \leq
 4\bigl(a-\norm{m}^2\bigr)
 \leq
 \frac32(2a-\tau^2).
\end{equation}
After division by $a$, this proves the equivalence of the
second and third conditions in
\eqref{eq:first-chaos-dominance}.

To relate them to the first condition, observe that
\begin{equation}\label{eq:first-chaos-fraction-remainder}
 1-\frac{\Var(H_1(g))}{\Var(F(g))}
 =
 \frac{\displaystyle\sum_{q=3}^{\infty}\E[H_q(g)^2]}
      {\tau^2+2a}.
\end{equation}
The denominator lies between $2a$ and $4a$.
Thus the first condition in
\eqref{eq:first-chaos-dominance} is equivalent to
\[
 \frac1{a_n}\sum_{q\geq3}\E[H_{n,q}(g_n)^2]
 \longrightarrow0.
\]

By \eqref{eq:remainder-energy-variance-gap},
$2\sum_{q\geq3}\E[H_q(g)^2]\leq2a-\tau^2$.
Conversely, Cauchy--Schwarz and the Hessian estimate give
\[
 \begin{aligned}
  (2a-\tau^2)^2
  &=
  \left(
   \sum_{q=3}^{\infty}(q-1)\E[H_q(g)^2]
  \right)^2\\
  &\leq
  \left(\sum_{q=3}^{\infty}\E[H_q(g)^2]\right)
  \left(\sum_{q=3}^{\infty}(q-1)^2\E[H_q(g)^2]\right)\\
  &\leq
  \left(\sum_{q=3}^{\infty}\E[H_q(g)^2]\right)
  \E\norm{D^2F(g)}_{\HS}^2\\
  &\leq
  4a\sum_{q=3}^{\infty}\E[H_q(g)^2].
 \end{aligned}
\]
Consequently,
$\sum_{q\geq3}\E[H_{n,q}(g_n)^2]/a_n\to0$
if and only if $(2a_n-\tau_n^2)/a_n\to0$.
Together with \eqref{eq:affine-error-variance-gap},
this completes the proof.
\end{proof}

The middle condition in \eqref{eq:first-chaos-dominance}
has the geometric meaning established in
Proposition~\ref{prop:projection-least-squares}:
the error of the optimal affine approximation becomes
negligible relative to the error of the optimal linear
approximation.

The condition in \eqref{eq:second-chaos-dominance} also
connects with the normalization issue in
Subsection~\ref{subsec:intrinsic-volume-clt-difficulty}.
If, in addition, $\tau_n\to\infty$, then
$\tau_n^2/\delta_n\to0$ forces the absolute value of the
parameter in \eqref{eq:projection-laplace-parameter}
to tend to infinity for every fixed $t\ne0$.

We conclude with three examples.

\begin{example}[Linear subspaces]
\label{ex:subspace-remainder}
Let $C=L$ be a linear subspace, and let $P_L$ denote its
orthogonal projection matrix.
Then $m=0$ and $A=P_L$.
Since $H(x)=x^{\mathsf T}P_Lx$, the remainder $F$ is
identically zero and $a=0$.
Also, $V_L=\dim L$ almost surely, so $\tau^2=0$.
Thus the subtraction removes all Gaussian fluctuations
in this case.
\end{example}

\begin{example}[Adding a linear-subspace factor]
\label{ex:linear-factor-remainder}
Let $C\subset\R^d$ be a closed convex cone and let
$L\subset\R^k$ be a linear subspace.
By \eqref{eq:chaos-products},
\[
 m_{C\times L}=m_C\oplus0,
 \qquad
 A_{C\times L}=A_C\oplus P_L.
\]
The product formula \eqref{eq:polarity-product} gives
$V_{C\times L}\eqD V_C+\dim L$.
Consequently, taking the product with $L$ preserves
$\norm{m}^2$, $a$, and $\tau^2$, while increasing
$\delta$ by $\dim L$.

The cancellation is also exact at the level of the
remainder: for $x\in\R^d$ and $y\in\R^k$,
$F_{C\times L}(x,y)=F_C(x)$.
The subtraction therefore removes the Gaussian fluctuations
produced by the linear factor.
In particular, the estimates for $F$ are unaffected by
the addition of arbitrarily many linear directions.
\end{example}

\begin{example}[The positive orthant]
\label{ex:orthant-remainder}
For $C=\R_+^d$, the projection is obtained by replacing each
negative coordinate by zero.
The mean projection and its average derivative are
\[
 m=\frac1{\sqrt{2\pi}}(1,\ldots,1),
 \qquad
 A=\frac12I_d.
\]
The intrinsic-volume random variable has the binomial
distribution with parameters $d$ and $1/2$, by the product
formula \eqref{eq:polarity-product} applied to the half-lines.
Hence $\delta=d/2$ and $\tau^2=d/4$.
Also, $a=\tr(A-A^2)=d/4$.

The remainder takes the explicit form
\begin{equation}\label{eq:orthant-F}
 F(x)=\frac12\sum_{i=1}^d x_i|x_i|.
\end{equation}
The preceding formulas give
\[
 \Var(H_1(g))=\frac{2d}{\pi},
 \qquad
 \Var(H_2(g))=\frac d2,
 \qquad
 \Var(F(g))=\frac{3d}{4}.
\]
Thus the fraction of the variance of $F(g)$ belonging
to the first Gaussian chaos is
$8/(3\pi)\approx0.849$, independently of $d$.
The first chaos accounts for most of the remaining
variance, but this fraction does not tend to one as
the dimension increases.
\end{example}

The preceding results relate the intrinsic-volume variance
to the mean Gaussian projection, the error of the optimal
linear approximation to the projection map, and the variance
of the second Gaussian-chaos component of $H$.
In the next two subsections, we give these comparisons a
metric interpretation by introducing a distance between
convex cones and studying approximation by linear subspaces.

\subsection{The \texorpdfstring{$L^2$-$\Pi$-metric}{L2-Pi-metric} and the angle between cones}
\label{subsec:gaussian-projection-distance}

Let $\CCC_d$ denote the family of all nonempty closed
convex cones in $\R^d$, where $d\geq1$.
We compare two cones by comparing their projections of the
same Gaussian vector.
This gives a distance compatible with the Gaussian
mean-square approximation considered in the preceding
subsections.
When the cone varies, we write $\delta_C=\E V_C$ for its
statistical dimension.

\begin{definition}[$L^2$-$\Pi$-metric]
\label{def:gaussian-projection-distance}
For $C,D\in\CCC_d$, define
\begin{equation}\label{eq:gaussian-projection-distance}
 d_{\Pi,2}(C,D)
 =
 \left(
  \E\norm{\Pi_C(g)-\Pi_D(g)}^2
 \right)^{1/2},
\end{equation}
where the same standard Gaussian vector $g\in\R^d$ is
used in both projections.
\end{definition}

\begin{proposition}[Basic properties of the $L^2$-$\Pi$-metric]
\label{prop:gaussian-projection-distance}
The function $d_{\Pi,2}$ is a metric on $\CCC_d$.
For every orthogonal transformation $Q$ of $\R^d$,
\begin{equation}\label{eq:projection-distance-isometries}
 d_{\Pi,2}(QC,QD)=d_{\Pi,2}(C,D),
 \qquad
 d_{\Pi,2}(C^\circ,D^\circ)=d_{\Pi,2}(C,D).
\end{equation}
For $C_i,D_i\in\CCC_{d_i}$, $i=1,2$, we have
\begin{equation}\label{eq:projection-distance-products}
 d_{\Pi,2}(C_1\times C_2,D_1\times D_2)^2
 =
 d_{\Pi,2}(C_1,D_1)^2+d_{\Pi,2}(C_2,D_2)^2.
\end{equation}
Moreover,
\begin{equation}\label{eq:projection-distance-spherical}
 d_{\Pi,2}(C,D)^2
 =
 d\,\E\norm{\Pi_C(U)-\Pi_D(U)}^2,
\end{equation}
where $U$ is uniformly distributed on the unit sphere
in $\R^d$.
\end{proposition}

\begin{proof}
Each projection map belongs to
$L^2(\R^d,\gamma_d;\R^d)$.
Nonnegativity, symmetry, and the triangle inequality
therefore follow from the corresponding properties
of the norm in this Hilbert space.
If $d_{\Pi,2}(C,D)=0$, the two projection maps agree
almost surely under Gaussian measure.
Continuity and full support of Gaussian measure imply
that they agree everywhere.
Since a closed convex cone is the set of fixed points
of its projection map, this gives $C=D$.

The identity
$\Pi_{QC}(x)=Q\Pi_C(Q^{\mathsf T}x)$ and orthogonal
invariance of Gaussian measure prove the first equality
in \eqref{eq:projection-distance-isometries}.
The second follows from Moreau's decomposition,
which gives
$\Pi_{C^\circ}(x)-\Pi_{D^\circ}(x)
=-(\Pi_C(x)-\Pi_D(x))$.

Projection onto a Cartesian product acts separately
on the two coordinates.
The squared Euclidean norm is additive across these
coordinates, which proves
\eqref{eq:projection-distance-products}.

Finally, write $g=RU$, where $R=\norm{g}$ and $U$ is
uniform on the unit sphere and independent of $R$.
Positive homogeneity of the projection gives
$\Pi_C(RU)-\Pi_D(RU)
=R(\Pi_C(U)-\Pi_D(U))$.
Taking expected squared norms and using $\E R^2=d$
proves \eqref{eq:projection-distance-spherical}.
\end{proof}

The projection maps have an additional property:
they all lie on a sphere in the Gaussian Hilbert space.
The center of this sphere is the map $x\mapsto x/2$.

\begin{proposition}[A spherical representation of convex cones]
\label{prop:cone-spherical-embedding}
For $C\in\CCC_d$, put
$J_C(x)=\Pi_C(x)-x/2$.
Then
\begin{equation}\label{eq:projection-sphere-pointwise}
 \norm{J_C(x)}^2=\frac14\norm{x}^2
 \qquad\text{for every }x\in\R^d.
\end{equation}
Consequently,
\begin{equation}\label{eq:cone-spherical-embedding}
 \norm{J_C}_{L^2}^2=\frac d4,
 \qquad
 \norm{J_C-J_D}_{L^2}=d_{\Pi,2}(C,D),
 \qquad
 J_{C^\circ}=-J_C.
\end{equation}
Thus $C\mapsto J_C$ is an isometric embedding of
$(\CCC_d,d_{\Pi,2})$ into the sphere of radius
$\sqrt d/2$ in $L^2(\R^d,\gamma_d;\R^d)$,
and polarity corresponds to the antipodal map.
\end{proposition}

\begin{proof}
By \eqref{eq:projection-energy},
\[
  \norm{\Pi_C(x)-x/2}^2
  =
  \norm{\Pi_C(x)}^2-\ip{\Pi_C(x)}{x}
  +\frac14\norm{x}^2\\
  =
  \frac14\norm{x}^2.
\]
Taking expectations proves the first identity in
\eqref{eq:cone-spherical-embedding}.
The distance identity follows because the same map
$x\mapsto x/2$ is subtracted from both projections.
Finally, Moreau's decomposition gives
$\Pi_{C^\circ}(x)=x-\Pi_C(x)$,
and hence $J_{C^\circ}(x)=-J_C(x)$.
\end{proof}

\begin{corollary}[Metric characterization of polarity]
\label{cor:metric-characterization-polarity}
For every $C,D\in\CCC_d$,
\begin{equation}\label{eq:polar-distance-identity}
 d_{\Pi,2}(C,D)^2+d_{\Pi,2}(C^\circ,D)^2=d.
\end{equation}
In particular,
$\operatorname{diam}(\CCC_d,d_{\Pi,2})=\sqrt d$,
and
\begin{equation}\label{eq:diametral-cones}
 d_{\Pi,2}(C,D)=\sqrt d
 \quad\Longleftrightarrow\quad
 D=C^\circ.
\end{equation}
Thus every cone has a unique farthest cone, namely its polar.
\end{corollary}

\begin{proof}
By \eqref{eq:cone-spherical-embedding} and the
parallelogram identity,
\[
  d_{\Pi,2}(C,D)^2+d_{\Pi,2}(C^\circ,D)^2
  =
  \norm{J_C-J_D}_{L^2}^2
  +\norm{J_C+J_D}_{L^2}^2
  =
  2\norm{J_C}_{L^2}^2+2\norm{J_D}_{L^2}^2
  =
  d.
\]
The remaining assertions follow because
$d_{\Pi,2}(C^\circ,D)=0$ if and only if $D=C^\circ$.
\end{proof}

The polarity operation is therefore determined by the
metric structure alone.
Equation~\eqref{eq:polar-distance-identity} also gives
an exact stability statement:
\begin{equation}\label{eq:polar-distance-deficit}
 d-d_{\Pi,2}(C,D)^2
 =
 d_{\Pi,2}(C^\circ,D)^2.
\end{equation}
The deficit from the maximal squared distance is exactly
the squared distance to the polar cone.

\begin{corollary}[Statistical dimension as a squared distance]
\label{cor:statistical-dimension-distance}
For every $C\in\CCC_d$,
\begin{equation}\label{eq:statistical-dimension-distance}
 d_{\Pi,2}(C,\{0\})^2=\delta_C,
 \qquad
 d_{\Pi,2}(C,\R^d)^2=d-\delta_C.
\end{equation}
\end{corollary}

\begin{proof}
Since $\Pi_{\{0\}}=0$, the first identity is
$\E\norm{\Pi_C(g)}^2=\delta_C$.
For the second, use $\Pi_{\R^d}(g)=g$ and Moreau's
decomposition to obtain
$\E\norm{\Pi_C(g)-g}^2
=\E\norm{\Pi_{C^\circ}(g)}^2=d-\delta_C$.
\end{proof}

\begin{remark}[Normalization]
\label{rem:normalized-projection-distance}
The normalized metric
$\overline d_{\Pi,2}=d_{\Pi,2}/\sqrt d$
has diameter one.
Equation~\eqref{eq:projection-distance-spherical} identifies
its square with
$\E\norm{\Pi_C(U)-\Pi_D(U)}^2$,
and \eqref{eq:polar-distance-identity} becomes
\[
 \overline d_{\Pi,2}(C,D)^2
 +
 \overline d_{\Pi,2}(C^\circ,D)^2
 =
 1.
\]
In particular, for sequences
$C_n,D_n\subset\R^{d_n}$, even when the dimensions vary,
$\overline d_{\Pi,2}(C_n,D_n)\to1$ if and only if
$\overline d_{\Pi,2}(C_n^\circ,D_n)\to0$.
We shall continue to use the unnormalized metric
$d_{\Pi,2}$ in the variance comparisons.
\end{remark}

The spherical representation also gives an angle between
cones: it is the usual Hilbert-space angle between their
images under $C\mapsto J_C$.

\begin{definition}[Associated angular metric]
\label{def:gaussian-projection-angle}
For $C,D\in\CCC_d$, define
\begin{equation}\label{eq:gaussian-projection-angle}
  \theta_{\Pi,2}(C,D)
  =
  \arccos\left(
   \frac4d\,\E\ip{J_C(g)}{J_D(g)}
  \right)
  =
  \arccos\left(
   1-\frac2d\,d_{\Pi,2}(C,D)^2
  \right).
\end{equation}
\end{definition}

The two expressions agree by expanding
$\norm{J_C-J_D}_{L^2}^2$ and using
$\norm{J_C}_{L^2}^2=\norm{J_D}_{L^2}^2=d/4$.

\begin{proposition}[Properties of the angular metric]
\label{prop:gaussian-projection-angle}
The function $\theta_{\Pi,2}$ is a metric on
$\CCC_d$ with diameter $\pi$.
It satisfies
\begin{equation}\label{eq:projection-angle-identities}
 d_{\Pi,2}(C,D)
 =
 \sqrt d\sin\frac{\theta_{\Pi,2}(C,D)}2,
 \qquad
 \theta_{\Pi,2}(C^\circ,D)
 =
 \pi-\theta_{\Pi,2}(C,D).
\end{equation}
In particular,
$\theta_{\Pi,2}(C,D)=\pi$ if and only if $D=C^\circ$.
Moreover,
\begin{equation}\label{eq:statistical-dimension-angle}
 \frac{\delta_C}{d}
 =
 \sin^2\frac{\theta_{\Pi,2}(C,\{0\})}{2}.
\end{equation}
\end{proposition}

\begin{proof}
The function $\theta_{\Pi,2}$ is the restriction of the
angular metric on a Hilbert-space sphere.
To recall its triangle inequality, let $u,v,w$ be unit
vectors, and let $\alpha$ and $\beta$ be the angles between
$u,v$ and between $v,w$, respectively.
If $\alpha+\beta\geq\pi$, the assertion is immediate.
Otherwise, decomposing $u$ and $w$ into their components
parallel and orthogonal to $v$, and applying
Cauchy--Schwarz, gives
\[
 \ip{u}{w}
 \geq
 \cos\alpha\cos\beta-\sin\alpha\sin\beta
 =
 \cos(\alpha+\beta).
\]
Thus the angle between $u$ and $w$ is at most
$\alpha+\beta$.
The other metric properties follow directly from
\eqref{eq:gaussian-projection-angle}.

The first identity in
\eqref{eq:projection-angle-identities} follows from the
half-angle formula, and the second follows from
$J_{C^\circ}=-J_C$.
The diameter and its equality case follow from
\eqref{eq:diametral-cones}.
Finally, combine the first identity in
\eqref{eq:projection-angle-identities} with
\eqref{eq:statistical-dimension-distance}.
\end{proof}

\begin{example}[Linear subspaces and principal angles]
\label{ex:subspaces-principal-angles}
Let $L,M\subset\R^d$ be linear subspaces of dimensions
$k$ and $\ell$, and let $P_L,P_M$ be their orthogonal
projection matrices.

We recall a short definition of their principal angles.
Choose matrices $Q_L\in\R^{d\times k}$ and
$Q_M\in\R^{d\times\ell}$ whose columns are orthonormal
bases of $L$ and $M$.
The principal angles
$0\leq\beta_1\leq\cdots\leq\beta_r\leq\pi/2$,
where $r=\min\{k,\ell\}$, are defined by requiring
$\cos\beta_1,\ldots,\cos\beta_r$ to be the singular
values of $Q_L^{\mathsf T}Q_M$, listed in decreasing order.
Changes of orthonormal bases multiply this matrix on the
left and right by orthogonal matrices, so the angles do
not depend on the chosen bases.

Since the projection maps are linear,
\begin{equation}\label{eq:subspace-projection-distance}
 d_{\Pi,2}(L,M)^2
 =
 \norm{P_L-P_M}_{\HS}^2
 =
 |k-\ell|+2\sum_{j=1}^r\sin^2\beta_j.
\end{equation}
Indeed, the squared Hilbert--Schmidt norm is
$k+\ell-2\tr(P_LP_M)$, while
$\tr(P_LP_M)
=\norm{Q_L^{\mathsf T}Q_M}_{\HS}^2
=\sum_{j=1}^r\cos^2\beta_j$.
When one of the subspaces is $\{0\}$, the sum is empty.

For $k=\ell$, the distance in
\eqref{eq:subspace-projection-distance} is $\sqrt2$ times
the chordal distance
$(\sum_{j=1}^k\sin^2\beta_j)^{1/2}$
considered by Conway, Hardin, and Sloane
\cite[Sections~2 and~5]{CHS96}.
Under the identification of linear maps with their
matrices, the spherical representation
\eqref{eq:cone-spherical-embedding} assigns to $L$
the matrix $P_L-I_d/2$.
On $\operatorname{Gr}(k,d)$, this representation differs
only by a translation from the classical representation
by traceless matrices $P_L-(k/d)I_d$ introduced in
\cite{CHS96}.
The use of the same center $I_d/2$ for subspaces of all
dimensions makes polarity correspond to the antipodal map.

For two lines in $\R^d$, $d\geq2$, the single principal
angle $\beta$ is the usual acute angle between them:
if $u$ and $v$ are unit spanning vectors, then
$\cos\beta=|\ip{u}{v}|$.
In this case,
\begin{equation}\label{eq:line-projection-angle}
 d_{\Pi,2}(L,M)=\sqrt2\sin\beta,
 \qquad
 \theta_{\Pi,2}(L,M)
 =
 2\arcsin\left(\sqrt{\frac2d}\sin\beta\right).
\end{equation}
In particular, for lines in $\R^2$,
$\theta_{\Pi,2}(L,M)=2\beta$.
Thus the angular metric is related to the usual angle,
with a normalization that depends on the ambient dimension.

Let $\operatorname{Gr}(k,d)$ denote the family of all
$k$-dimensional linear subspaces of $\R^d$.
The principal-angle formula also gives
\begin{equation}\label{eq:grassmannian-diameters}
 \begin{aligned}
  \operatorname{diam}
  \bigl(\operatorname{Gr}(k,d),d_{\Pi,2}\bigr)
  &=
  \sqrt{2\min\{k,d-k\}},\\
  \operatorname{diam}
  \bigl(\operatorname{Gr}(k,d),\theta_{\Pi,2}\bigr)
  &=
  2\arcsin
  \sqrt{\frac{2\min\{k,d-k\}}d}.
 \end{aligned}
\end{equation}
To see the first equality when $2k\leq d$, use
$\sum_{j=1}^k\sin^2\beta_j\leq k$ and choose two
orthogonal $k$-dimensional subspaces to attain equality.
When $2k>d$, orthogonal complementation preserves
$d_{\Pi,2}$ and reduces the problem to dimension $d-k$.
The angular formula follows from
\eqref{eq:projection-angle-identities}.
These formulas include $k=0$ and $k=d$, when the
Grassmannian consists of a single point.
The full diameter of $\CCC_d$ is attained within
$\operatorname{Gr}(k,d)$ precisely when $d$ is even
and $k=d/2$.
\end{example}

\begin{remark}[The Hausdorff topology]
\label{rem:projection-distance-topology}
For fixed $d$, the metrics $d_{\Pi,2}$ and $\theta_{\Pi,2}$
induce the usual topology obtained by intersecting cones
with the closed unit ball
$B^d=\{x\in\R^d:\norm{x}\leq1\}$.
More precisely,
\begin{equation}\label{eq:cone-topology-equivalence}
 \begin{aligned}
  d_{\Pi,2}(C_n,C)\longrightarrow0
  &\quad\Longleftrightarrow\quad
  d_{\mathrm{Haus}}(C_n\cap B^d,C\cap B^d)
  \longrightarrow0\\
  &\quad\Longleftrightarrow\quad
  \sup_{x\in B^d}
  \norm{\Pi_{C_n}(x)-\Pi_C(x)}
  \longrightarrow0.
 \end{aligned}
\end{equation}
Here $d_{\mathrm{Haus}}$ is the Hausdorff distance between
nonempty compact sets.

For the equivalence between the last two conditions,
observe first that $\Pi_C(x)\in B^d$ whenever $x\in B^d$.
Hausdorff convergence of the truncated cones therefore
implies pointwise convergence of their projections on
$B^d$: every limit point of the minimizing points belongs
to $C\cap B^d$ and minimizes the distance to $x$, so
uniqueness identifies it as $\Pi_C(x)$.
The common Lipschitz constant one upgrades this convergence
to uniform convergence on $B^d$.
Conversely, uniform convergence of the projections gives
Hausdorff convergence by projecting each point of one
truncated cone onto the other.

Uniform convergence on $B^d$ implies convergence in
Gaussian $L^2$ by positive homogeneity.
For the converse, the differences
$\Pi_{C_n}-\Pi_C$ are uniformly $2$-Lipschitz.
If their convergence were not uniform on $B^d$,
compactness of the ball would provide a subsequence
and a fixed open ball on which their norms remain
bounded away from zero.
That ball has positive Gaussian measure, contradicting
convergence in Gaussian $L^2$.
The equivalence with convergence in $\theta_{\Pi,2}$
follows from \eqref{eq:projection-angle-identities}.

This topology is compact.
Indeed, the projection maps are $1$-Lipschitz and vanish
at zero, so Arzel\`a--Ascoli and positive homogeneity
give, from any sequence, a subsequence converging
uniformly on bounded sets to a map $P$.
Put $C=\{x:P(x)=x\}$.
The convergence implies that if $x_n\in C_n$ and
$x_n\to x$, then $x\in C$.
Conversely, every $x\in C$ is approximated by
$\Pi_{C_n}(x)\in C_n$.
Approximating two points of $C$ in this way and taking
nonnegative linear combinations shows that $C$ is a
closed convex cone.
Also, $P(x)\in C$ for every $x$.
For $y\in C$, the minimizing property gives
$\norm{x-\Pi_{C_n}(x)}
\leq\norm{x-\Pi_{C_n}(y)}$.
Passing to the limit yields
$\norm{x-P(x)}\leq\norm{x-y}$,
so $P=\Pi_C$.
Equation~\eqref{eq:cone-topology-equivalence} now proves
compactness.
\end{remark}

Finally, the chaos expansion of the projection gives an
orthogonal decomposition of the squared distance.
This connects the metric directly to the coefficients
studied in the preceding subsections.

\begin{proposition}[Chaos representation of the $L^2$-$\Pi$-metric]
\label{prop:projection-distance-chaoses}
For $C,D\in\CCC_d$,
\begin{equation}\label{eq:projection-distance-chaoses}
  d_{\Pi,2}(C,D)^2
  =
  \norm{m_C-m_D}^2
  +\norm{A_C-A_D}_{\HS}^2
  +
  \sum_{q=2}^{\infty}
  \E\norm{(\Pi_C)_q(g)-(\Pi_D)_q(g)}^2.
\end{equation}
In particular,
\begin{equation}\label{eq:projection-distance-low-chaoses}
 \norm{m_C-m_D}^2+\norm{A_C-A_D}_{\HS}^2
 \leq d_{\Pi,2}(C,D)^2.
\end{equation}
Thus $C\mapsto m_C$ and $C\mapsto A_C$ are
$1$-Lipschitz with respect to the Euclidean and
Hilbert--Schmidt norms, respectively.
\end{proposition}

\begin{proof}
Subtract the two expansions in
\eqref{eq:projection-chaos-expansion}.
The constant term is $m_C-m_D$, the linear term is
$(A_C-A_D)g$, and the remaining terms belong to the
Gaussian chaoses of orders at least two.
Orthogonality in
$L^2(\R^d,\gamma_d;\R^d)$ and the identity
$\E\norm{(A_C-A_D)g}^2=\norm{A_C-A_D}_{\HS}^2$
give \eqref{eq:projection-distance-chaoses}.
Dropping the nonnegative sum proves
\eqref{eq:projection-distance-low-chaoses}.
\end{proof}

The statistical dimension thus has the metric
interpretation
$\delta_C=d_{\Pi,2}(C,\{0\})^2$.

\begin{corollary}[Distance to the reflected cone]
\label{cor:distance-reflected-cone}
Let $C\subseteq\R^d$ be a nonempty closed convex
cone.
Recall that $m_C=\E\Pi_C(g)$,
$A_C=\E[\Pi_C(g)g^{\mathsf T}]$, and
\[
 a
 =
 \E\norm{\Pi_C(g)-A_Cg}^2
 =
 \tr(A_C-A_C^2).
\]
Then
\begin{equation}\label{eq:distance-reflected-cone}
 2\norm{m_C}
 \leq d_{\Pi,2}(C,-C)
 \leq2\sqrt a
 \leq2\sqrt3\,\norm{m_C}.
\end{equation}
Moreover, writing $\tau^2=\Var(V_C)$, we have
\begin{equation}\label{eq:variance-distance-reflected-cone}
 2\tau^2
 \leq d_{\Pi,2}(C,-C)^2
 \leq6\tau^2.
\end{equation}
\end{corollary}

\begin{proof}
The identity
$\Pi_{-C}(x)=-\Pi_C(-x)$ gives
\[
 d_{\Pi,2}(C,-C)^2
 =
 \E\norm{\Pi_C(g)+\Pi_C(-g)}^2.
\]
Since $g$ and $-g$ have the same distribution,
the vector inside the norm has expectation
$2m_C$.
Jensen's inequality therefore gives
$d_{\Pi,2}(C,-C)^2\geq4\norm{m_C}^2$.

For the upper bound, subtract the optimal
linear approximation from each projection.
The two linear terms cancel, and hence
\begin{align*}
 d_{\Pi,2}(C,-C)^2
 &=
 \E\norm{
   \bigl(\Pi_C(g)-A_Cg\bigr)
   +
   \bigl(\Pi_C(-g)+A_Cg\bigr)
 }^2\\
 &\leq
 2\E\norm{\Pi_C(g)-A_Cg}^2
 +
 2\E\norm{\Pi_C(-g)+A_Cg}^2\\
 &=4a.
\end{align*}
The last equality again follows from the
symmetry of the Gaussian distribution.

By Propositions~\ref{prop:removal}
and~\ref{prop:first-chaos-bounds},
$a\leq\frac32\tau^2\leq3\norm{m_C}^2$.
This proves \eqref{eq:distance-reflected-cone}.
Finally, the inequalities
$\tau^2\leq2\norm{m_C}^2$ and
$a\leq\frac32\tau^2$ give
\[
 2\tau^2
 \leq4\norm{m_C}^2
 \leq d_{\Pi,2}(C,-C)^2
 \leq4a
 \leq6\tau^2,
\]
which proves
\eqref{eq:variance-distance-reflected-cone}.
\end{proof}

For a linear subspace, the projection map has only a
first Gaussian-chaos component, so
\eqref{eq:projection-distance-chaoses} simplifies further.
In the next subsection, we use this observation to
relate the intrinsic-volume variance to the squared
distance from a cone to the family of all linear subspaces.

\subsection{Approximation by linear subspaces}
\label{subsec:approximation-by-subspaces}

Let $\LS_d$ denote the family of all linear subspaces
of $\R^d$, including $\{0\}$ and $\R^d$.
We write
$d_{\Pi,2}(C,\LS_d)
=\inf_{L\in\LS_d}d_{\Pi,2}(C,L)$.
The spectral description below will show that this infimum
is attained.

For a fixed cone $C$, we retain the notation
$A=\E[D\Pi_C(g)]$, $m=\E\Pi_C(g)$,
$\delta=\E V_C$, $\tau^2=\Var(V_C)$, and
$a=\tr(A-A^2)$.
If $P_L$ is the orthogonal projection matrix of a linear
subspace $L$, then \eqref{eq:least-squares} gives
\begin{equation}\label{eq:cone-subspace-distance}
  d_{\Pi,2}(C,L)^2
  =
  a+\norm{A-P_L}_{\HS}^2
  =
  \delta+\dim L-2\tr(AP_L).
\end{equation}
The second equality uses $P_L^2=P_L$,
$\tr P_L=\dim L$, and $\tr A=\delta$.
Thus the matrix $A$ determines the distances from $C$
to all linear subspaces.
Finding a nearest subspace amounts to finding a nearest
orthogonal projection matrix to $A$.

\begin{proposition}[Nearest linear subspaces]
\label{prop:nearest-subspaces}
Let $\lambda_1\geq\cdots\geq\lambda_d$ be the eigenvalues
of $A$.
Then
\begin{equation}\label{eq:nearest-subspace-distance}
 d_{\Pi,2}(C,\LS_d)^2
 =
 \sum_{i=1}^d\min\{\lambda_i,1-\lambda_i\}.
\end{equation}
The complete family of minimizing subspaces is
\begin{equation}\label{eq:nearest-subspace-minimizers}
 L=E_{>1/2}\oplus M,
 \qquad
 M\subseteq E_{=1/2},
\end{equation}
where $E_{>1/2}$ is the sum of the eigenspaces of $A$
with eigenvalues greater than $1/2$,
$E_{=1/2}=\ker(A-I_d/2)$, and $M$ is an arbitrary
linear subspace of $E_{=1/2}$.
In particular, the nearest subspace, with its dimension
unrestricted, is unique if and only if $1/2$ is not an
eigenvalue of $A$.

For a prescribed dimension $k\in\{0,\ldots,d\}$,
\begin{equation}\label{eq:fixed-dimension-subspace-distance}
 \min_{L\in\operatorname{Gr}(k,d)}
 d_{\Pi,2}(C,L)^2
 =
 \delta+k-2\sum_{i=1}^k\lambda_i.
\end{equation}
An optimal subspace is spanned by orthonormal eigenvectors
corresponding to the $k$ largest eigenvalues.
The sum is interpreted as zero when $k=0$.
\end{proposition}

\begin{proof}
Choose an orthonormal eigenbasis $e_1,\ldots,e_d$ of $A$
corresponding to the indicated eigenvalues.
For every orthogonal projection matrix $P_L$,
the quantities $\norm{P_Le_i}^2$ lie in $[0,1]$ and
sum to $\dim L$.
Equation~\eqref{eq:cone-subspace-distance} therefore gives
\[
 d_{\Pi,2}(C,L)^2
 =
 \sum_{i=1}^d
 \left(
  \lambda_i+(1-2\lambda_i)\norm{P_Le_i}^2
 \right).
\]
Each summand is at least
$\min\{\lambda_i,1-\lambda_i\}$.
Equality requires $P_Le_i=e_i$ when $\lambda_i>1/2$
and $P_Le_i=0$ when $\lambda_i<1/2$.
The eigenspace corresponding to $1/2$ imposes no further
restriction.
These conditions are equivalent to
\eqref{eq:nearest-subspace-minimizers}.
They prove both attainment and
\eqref{eq:nearest-subspace-distance}, as well as the
uniqueness assertion.

If $\dim L=k$, then
\[
 \tr(AP_L)
 =
 \sum_{i=1}^d\lambda_i\norm{P_Le_i}^2
 \leq
 \sum_{i=1}^k\lambda_i.
\]
Indeed, the coefficients in this weighted sum belong
to $[0,1]$ and sum to $k$, so the sum is maximized by
placing unit weight on the $k$ largest eigenvalues.
Equality is attained by the span of the corresponding
eigenvectors.
Substitution into \eqref{eq:cone-subspace-distance}
proves \eqref{eq:fixed-dimension-subspace-distance}.
\end{proof}

The distinction between approximation by a linear map
and approximation by a subspace projection is now explicit.
The former has error $a$, while the latter requires the
approximating matrix to be an orthogonal projection.
The additional restriction changes the error by at most
a universal factor.

\begin{theorem}[A metric interpretation of the
intrinsic-volume variance]
\label{thm:variance-distance-subspaces}
For every nonempty closed convex cone $C\subset\R^d$,
\begin{align}
 a
 &\leq d_{\Pi,2}(C,\LS_d)^2\leq2a,
 \label{eq:subspace-distance-a}\\
 \frac12\tau^2
 &\leq d_{\Pi,2}(C,\LS_d)^2\leq3\tau^2.
 \label{eq:variance-distance-subspaces}
\end{align}
\end{theorem}

\begin{proof}
For every $\lambda\in[0,1]$,
\[
 \lambda(1-\lambda)
 \leq
 \min\{\lambda,1-\lambda\}
 \leq
 2\lambda(1-\lambda).
\]
Summing over the eigenvalues of $A$ and using
$a=\tr(A-A^2)$ proves
\eqref{eq:subspace-distance-a}.
Combining this with \eqref{eq:a-comparison} gives
\eqref{eq:variance-distance-subspaces}.
\end{proof}

The statistical dimension is the squared distance to the
zero cone, by \eqref{eq:statistical-dimension-distance}.
The intrinsic-volume variance has a corresponding
interpretation: up to universal constants, it is the
squared distance to the family of all linear subspaces.
In particular, for a sequence of cones, the condition
$\tau_n^2\to\infty$ is equivalent to
$d_{\Pi,2}(C_n,\LS_{d_n})\to\infty$.

We collect the geometric comparisons established so far.

\begin{corollary}[Equivalent geometric measures of variance]
\label{cor:geometric-variance-comparisons}
For every nonempty closed convex cone,
\begin{equation}\label{eq:geometric-variance-comparisons}
  \tau^2
  \asymp
  \norm{m}^2
  \asymp
  \tr(A-A^2)
  \asymp
  2\delta-\Var(H_2(g))
  \asymp
  d_{\Pi,2}(C,\LS_d)^2.
\end{equation}
Here $\asymp$ denotes comparison above and below by
positive constants independent of the cone and the
ambient dimension.
\end{corollary}

\begin{proof}
Combine \eqref{eq:a-definition},
\eqref{eq:mean-projection-variance-comparison},
\eqref{eq:second-chaos-variance-deficit}, and
\eqref{eq:variance-distance-subspaces}.
\end{proof}

The metric comparison also gives a geometric interpretation
of second-chaos dominance.

\begin{corollary}[Relative approximation by linear subspaces]
\label{cor:relative-subspace-approximation}
Let $C_n\subset\R^{d_n}$ be a sequence of nonempty closed
convex cones with $\delta_n=\E V_{C_n}>0$, and write
$\tau_n^2=\Var(V_{C_n})$.
Let $H_n(x)=\norm{\Pi_{C_n}(x)}^2$, and denote its second
Gaussian-chaos component by $H_{n,2}$.
Then
\begin{equation}\label{eq:relative-subspace-approximation}
 \begin{aligned}
  \frac{\Var(H_{n,2}(g_n))}
       {\Var(H_n(g_n))}
  \longrightarrow1
  &\quad\Longleftrightarrow\quad
  \frac{\tau_n^2}{\delta_n}\longrightarrow0\\
  &\quad\Longleftrightarrow\quad
  \frac{d_{\Pi,2}(C_n,\LS_{d_n})}
       {d_{\Pi,2}(C_n,\{0\})}
  \longrightarrow0,
 \end{aligned}
\end{equation}
where $g_n$ is a standard Gaussian vector in $\R^{d_n}$.
\end{corollary}

\begin{proof}
The first equivalence is
\eqref{eq:second-chaos-dominance}.
For the second, divide
\eqref{eq:variance-distance-subspaces} by $\delta_n$
and use
$d_{\Pi,2}(C_n,\{0\})^2=\delta_n$.
\end{proof}

Thus the second chaos accounts asymptotically for all the
variance of the squared projection precisely when the
distance to linear subspaces is negligible relative to
$\sqrt{\delta_n}$.
The absolute distance to linear subspaces may still
diverge, as it does whenever $\tau_n^2\to\infty$.
This relative approximation gives a geometric interpretation
of the normalization issue discussed in
Subsection~\ref{subsec:intrinsic-volume-clt-difficulty}.

Polarity also determines the farthest linear subspaces.

\begin{proposition}[Farthest linear subspaces]
\label{prop:farthest-subspaces}
For every $L\in\LS_d$,
\begin{equation}\label{eq:complementary-subspace-distances}
 d_{\Pi,2}(C,L)^2+d_{\Pi,2}(C,L^\perp)^2=d.
\end{equation}
Consequently,
\begin{equation}\label{eq:farthest-subspace-distance}
  \max_{L\in\LS_d}d_{\Pi,2}(C,L)^2
  =
  d-d_{\Pi,2}(C,\LS_d)^2
  =
  \sum_{i=1}^d\max\{\lambda_i,1-\lambda_i\},
\end{equation}
where $\lambda_1,\ldots,\lambda_d$ are the eigenvalues
of $A$.
A subspace is farthest from $C$ if and only if its
orthogonal complement is nearest to $C$.
\end{proposition}

\begin{proof}
The polar of a linear subspace $L$ is $L^\perp$.
Equation~\eqref{eq:complementary-subspace-distances}
therefore follows from
\eqref{eq:polar-distance-identity} and symmetry of the metric.
Since orthogonal complementation is a bijection of
$\LS_d$, maximizing one term is equivalent to
minimizing the other.
Equation~\eqref{eq:nearest-subspace-distance} then gives
\eqref{eq:farthest-subspace-distance}.
\end{proof}

In particular, the intrinsic-volume variance is also
comparable to
$d-\max_{L\in\LS_d}d_{\Pi,2}(C,L)^2$.
It therefore measures both the squared distance to the
nearest subspace and the deficit of the largest squared
distance to a subspace from $d$.

We next consider a different extremal problem: how far
a cone can be from the entire family $\LS_d$.

\begin{proposition}[Cones farthest from linear subspaces]
\label{prop:cones-farthest-from-subspaces}
For every $d\geq1$,
\begin{equation}\label{eq:max-distance-to-subspaces}
 \max_{C\in\CCC_d}
 d_{\Pi,2}(C,\LS_d)^2
 =
 \frac d2.
\end{equation}
For a cone $C\in\CCC_d$, the following conditions
are equivalent:
\begin{enumerate}
 \item $d_{\Pi,2}(C,\LS_d)^2=d/2$.
 \item $A_C=I_d/2$.
 \item $d_{\Pi,2}(C,L)=\sqrt{d/2}$ for every
       $L\in\LS_d$.
 \item $\theta_{\Pi,2}(C,L)=\pi/2$ for every
       $L\in\LS_d$.
 \item $\langle J_C,J_L\rangle_{L^2}=0$ for every
       $L\in\LS_d$.
\end{enumerate}
Moreover,
\begin{equation}\label{eq:subspace-distance-deficit}
 \frac d2-d_{\Pi,2}(C,\LS_d)^2
 =
 \left\|A_C-\frac12I_d\right\|_*,
\end{equation}
where $\norm{\cdot}_*$ is the nuclear norm, defined as
the sum of the singular values.
Every self-dual cone satisfies the equivalent conditions
above.
\end{proposition}

\begin{proof}
For every $\lambda\in[0,1]$,
$\min\{\lambda,1-\lambda\}
=\frac12-|\lambda-\frac12|$.
Applying \eqref{eq:nearest-subspace-distance} gives
\[
 \frac d2-d_{\Pi,2}(C,\LS_d)^2
 =
 \sum_{i=1}^d\left|\lambda_i-\frac12\right|
 =
 \left\|A_C-\frac12I_d\right\|_*.
\]
The last equality uses symmetry of $A_C$.
This proves \eqref{eq:subspace-distance-deficit} and
shows that the squared distance is at most $d/2$,
with equality precisely when $A_C=I_d/2$.
The positive orthant has this matrix by
Example~\ref{ex:orthant-remainder}, so the maximum
in \eqref{eq:max-distance-to-subspaces} is attained.

If $A_C=I_d/2$, then $\delta_C=d/2$, and
\eqref{eq:cone-subspace-distance} gives
$d_{\Pi,2}(C,L)^2=d/2$ for every linear subspace $L$.
Conversely, this equidistance property implies that
the minimum over $L$ is $d/2$.
This proves the equivalence of the first three conditions.

By \eqref{eq:gaussian-projection-angle},
$d_{\Pi,2}(C,L)^2=d/2$ is equivalent both to
$\theta_{\Pi,2}(C,L)=\pi/2$ and to
$\langle J_C,J_L\rangle_{L^2}=0$.
Finally, if $C^\circ=-C$, then
\eqref{eq:self-dual-chaos} gives $A_C=I_d/2$.
\end{proof}

This applies to positive orthants and Lorentz cones,
as well as to the cone of positive semidefinite matrices
in the Euclidean space of real symmetric matrices
with inner product $\langle X,Y\rangle=\tr(XY)$.
The spectral condition $A_C=I_d/2$ is, however,
strictly more general than self-duality.

\begin{example}[A maximizer that is not self-dual]
\label{ex:non-self-dual-maximizer}
Consider the cone
\[
 C=
 \left\{
  (t,x,y)\in\R^3:
  t\geq0,\quad
  |x|\leq2^{-1/4}t,\quad
  |y|\leq2^{-1/4}t
 \right\}.
\]
Minimizing a linear functional over its square section
at $t=1$ shows that its positive dual is
\[
 -C^\circ
 =
 \left\{
  (t,x,y)\in\R^3:
  t\geq2^{-1/4}(|x|+|y|)
 \right\}.
\]

Let $Q$ be rotation through $\pi/4$ around the $t$-axis.
The inequalities defining $QC$ are
$|x+y|\leq2^{1/4}t$ and
$|x-y|\leq2^{1/4}t$.
Since
$\max\{|x+y|,|x-y|\}=|x|+|y|$,
we obtain $QC=-C^\circ$.

The cone is invariant under the two coordinate sign changes
$x\mapsto-x$ and $y\mapsto-y$, and under interchange of
$x$ and $y$.
By \eqref{eq:chaos-rotations}, these symmetries force
$A_C$ to be diagonal with equal $x$- and $y$-entries.
It is therefore invariant under $Q$.
On the other hand, reflection preserves the matrix,
while polarity replaces it by its complement to the
identity. Hence
\[
 A_C
 =
 A_{QC}
 =
 A_{-C^\circ}
 =
 I_3-A_C.
\]
Thus $A_C=I_3/2$, and
$d_{\Pi,2}(C,\LS_3)^2=3/2$.

Nevertheless, $C$ is not self-dual.
The point $(1,2^{-1/4},2^{-1/4})$ belongs to $C$,
but membership in $-C^\circ$ would require
$1\geq\sqrt2$.
Therefore $C^\circ\ne-C$.
\end{example}

For polyhedral cones, the matrix $A$ has a further
interpretation as the average of projection matrices
onto random face spans.
The underlying descriptions of intrinsic volumes and
projection regions are classical; see
McCoy and Tropp
\cite[Definition~2.1 and Fact~7.1]{MT14}.

Suppose that $C$ is polyhedral.
For $x\in\R^d$, let $L_C(x)$ denote the linear span
of the unique face of $C$ whose relative interior
contains $\Pi_C(x)$.
Thus $L_C(g)$ is a random linear subspace.
In the expectations below, the metric is evaluated
between the realized subspaces: by
\eqref{eq:subspace-projection-distance}, its square is
the squared Hilbert--Schmidt distance between their
orthogonal projection matrices.

We also recall that a Fr\'echet mean of a random element
of a metric space is a point minimizing its expected
squared distance.
Here the metric space in which the mean is sought is
$(\LS_d,d_{\Pi,2})$.

\begin{proposition}[Random face spans and Fr\'echet means]
\label{prop:random-face-spans}
Suppose that $C$ is polyhedral.
Then
\begin{equation}\label{eq:polyhedral-projection-derivative}
 D\Pi_C(x)=P_{L_C(x)}
 \qquad\text{for almost every }x\in\R^d,
\end{equation}
and
\begin{equation}\label{eq:random-face-span-moments}
 \dim L_C(g)\eqD V_C,
 \qquad
 A=\E P_{L_C(g)}.
\end{equation}
For every deterministic linear subspace $L$,
\begin{equation}\label{eq:random-face-distance}
 d_{\Pi,2}(C,L)^2
 =
 \E\,d_{\Pi,2}(L_C(g),L)^2.
\end{equation}
If $g'$ is an independent copy of $g$, then
\begin{equation}\label{eq:random-face-spread}
  a
  =
  \E\norm{P_{L_C(g)}-A}_{\HS}^2
  =
  \frac12\E\,d_{\Pi,2}(L_C(g),L_C(g'))^2.
\end{equation}
Consequently,
\begin{equation}\label{eq:random-face-variance-comparison}
 \tau^2
 \leq
 \E\,d_{\Pi,2}(L_C(g),L_C(g'))^2
 \leq
 3\tau^2.
\end{equation}
The nearest linear subspaces to $C$ are precisely the
Fr\'echet means of the distribution of $L_C(g)$ in
$(\LS_d,d_{\Pi,2})$.
\end{proposition}

\begin{proof}
For a face $G$ of $C$, the classical decomposition
of its projection region is
\[
 \Pi_C^{-1}(\relint G)=\relint G+N_G,
 \qquad
 N_G=C^\circ\cap(\spanop G)^\perp.
\]
The spans of $G$ and its normal face $N_G$ are
orthogonal complements.
Therefore
$\relint G+\relint N_G$ is open in $\R^d$.
On this set, Moreau's decomposition gives
$\Pi_C(x)=P_{\spanop G}x$.
The finitely many such open regions cover $\R^d$
except for polyhedral boundaries of Lebesgue measure zero.
This proves \eqref{eq:polyhedral-projection-derivative}.

Averaging the derivative identity gives
$A=\E P_{L_C(g)}$.
The identity in distribution in
\eqref{eq:random-face-span-moments} is the face
interpretation of the conic intrinsic volumes.
In particular,
$\E\dim L_C(g)=\delta$.

For a deterministic subspace $L$, use
\eqref{eq:subspace-projection-distance} to obtain
\[
  \E\,d_{\Pi,2}(L_C(g),L)^2
  =
  \E\norm{P_{L_C(g)}-P_L}_{\HS}^2
  =
  \delta+\dim L-2\tr(AP_L)
  =
  d_{\Pi,2}(C,L)^2.
\]
This proves \eqref{eq:random-face-distance}.

Since each $P_{L_C(g)}$ is an orthogonal projection
and has mean $A$,
\[
  \E\norm{P_{L_C(g)}-A}_{\HS}^2
  =
  \E\tr(P_{L_C(g)})-\tr(A^2)
  =
  \delta-\tr(A^2)
  =
  a.
\]
Independence of $g$ and $g'$ then gives
\[
  \E\,d_{\Pi,2}(L_C(g),L_C(g'))^2
  =
  \E\norm{P_{L_C(g)}-P_{L_C(g')}}_{\HS}^2
  =
  2\E\norm{P_{L_C(g)}-A}_{\HS}^2
  =
  2a.
\]
Equation~\eqref{eq:a-comparison} proves
\eqref{eq:random-face-variance-comparison}.

Finally, \eqref{eq:random-face-distance} shows that
minimizing $d_{\Pi,2}(C,L)^2$ over $L\in\LS_d$
is exactly the Fr\'echet-mean minimization problem
for $L_C(g)$.
\end{proof}

Since $\Var(\dim L_C(g))=\tau^2$,
\eqref{eq:random-face-variance-comparison} compares
the variability of the entire random face span with
the variability of its dimension.
For this distribution of random subspaces, fluctuations
in dimension therefore control their geometric
dispersion, including changes in orientation.

\begin{remark}[Hessian identities for polyhedral cones]
\label{rem:polyhedral-hessian-identities}
For a polyhedral cone,
\eqref{eq:polyhedral-projection-derivative} gives
$D\Pi_C(x)^2=D\Pi_C(x)$ almost everywhere.
Consequently,
$\E\norm{D^2H(g)}_{\HS}^2=4\delta$
and
$\E\norm{D^2F(g)}_{\HS}^2=4a$.
In particular, the final inequality in
\eqref{eq:higher-chaos-weighted-bound} is an equality
for polyhedral cones.
\end{remark}

The second Gaussian-chaos matrix determines the distances
from a cone to all linear subspaces.
In particular, all self-dual cones in a given ambient
dimension have the same distances to linear subspaces.
Their first and higher Gaussian-chaos components, and
hence their intrinsic-volume variances, can nevertheless
differ.
In the next subsection, we study circular cones and
compare, in particular, the Lorentz cone with the
positive orthant.

\subsection{Circular cones and first-chaos dominance}
\label{subsec:circular-cones}

The results of this subsection are not needed for the proof of the
main central limit theorem. They explain why the first Gaussian-chaos
component is retained when the second component is removed.
For circular cones, the first component accounts asymptotically for
all the variance of the remainder whenever the intrinsic-volume
variance diverges. Removing the first component as well leaves a
random variable with bounded variance; for a fixed opening angle,
its limiting distribution is a product of two independent Gaussian
variables.

The asymptotic formulas for the statistical dimension and
the intrinsic-volume variance of circular cones with fixed
opening angle were obtained by McCoy and Tropp
\cite[Subsection~6.3.1]{MT14}.
We give a derivation that also describes the Gaussian-chaos
components and their geometric meaning, and investigate
the case where the opening angle varies with the dimension.

For each $n\geq1$, let $\alpha_n\in(0,\pi/2)$ and consider the circular
cone
\begin{equation}\label{eq:circular-cones}
 C_n
 =
 \bigl\{(t,u)\in\R\times\R^n:
          t\geq\norm{u}\cot\alpha_n\bigr\}.
\end{equation}
The angle $\alpha_n$ is allowed to depend on $n$.
Write $g_n=(T,U_n)$, where $T$ is a standard Gaussian random variable,
$U_n$ is a standard Gaussian vector in $\R^n$, and $T$ and $U_n$ are
independent. We put $R_n=\norm{U_n}$.

For these cones, we use the notation
$H_n(x)=\norm{\Pi_{C_n}(x)}^2$ and denote its $q$th Gaussian-chaos
component by $H_{n,q}$. We also write
$m_n=\E\Pi_{C_n}(g_n)$ and $A_n=\E D\Pi_{C_n}(g_n)$.
Recall that
$\delta_n=\tr A_n$,
$\tau_n^2=\Var(V_{C_n})$,
$a_n=\tr(A_n-A_n^2)$, and
\begin{equation}\label{eq:circular-remainder}
 F_n(x)
 =
 H_n(x)-x^\top A_nx
 =
 H_{n,1}(x)+\sum_{q\geq3}H_{n,q}(x),
 \qquad
 H_{n,1}(x)=2\ip{m_n}{x}.
\end{equation}
The series in \eqref{eq:circular-remainder} converges in
$L^2(\R^{n+1},\gamma_{n+1})$.
Below, $\norm{\cdot}_2$ denotes the corresponding $L^2$-norm.

\begin{lemma}[Projection and Gaussian-chaos approximations]
\label{lem:circular-chaos-approximation}
For every $n\geq1$,
\begin{equation}\label{eq:circular-projection-expansion}
 H_n(g_n)
 =
 \cos^2\alpha_n\,T^2
 +\sin^2\alpha_n\,R_n^2
 +\sin(2\alpha_n)\,TR_n
 +D_n,
\end{equation}
where
\begin{equation}\label{eq:circular-projection-error}
 D_n
 =
 (T\sin\alpha_n-R_n\cos\alpha_n)_+^2
 -
 (-T\cos\alpha_n-R_n\sin\alpha_n)_+^2
\end{equation}
and $x_+=\max\{x,0\}$.
The correction satisfies $\abs{D_n}\leq T^2$, and hence
$\norm{D_n}_2\leq\sqrt3$.

Moreover,
\begin{equation}\label{eq:circular-chaos-approximation}
\begin{aligned}
 \norm{
   F_n(g_n)-\sin(2\alpha_n)\,TR_n
 }_2
 &\leq \norm{D_n}_2,\\
 \norm{
   H_{n,1}(g_n)-\sin(2\alpha_n)(\E R_n)T
 }_2
 &\leq \norm{D_n}_2,\\
 \norm{
   F_n(g_n)-H_{n,1}(g_n)
   -\sin(2\alpha_n)\,T(R_n-\E R_n)
 }_2
 &\leq \norm{D_n}_2.
\end{aligned}
\end{equation}
In particular, the variance contributed by the chaoses of orders
at least three is bounded uniformly in the dimension and the angle:
\begin{equation}\label{eq:circular-higher-chaos-bound}
 \Var\bigl(F_n(g_n)-H_{n,1}(g_n)\bigr)
 =
 \sum_{q\geq3}\E H_{n,q}(g_n)^2
 \leq (1+\sqrt3)^2.
\end{equation}
\end{lemma}

\begin{proof}
Rotational symmetry reduces the projection of $(t,u)$ to the
two-dimensional plane spanned by the axis of the cone and $u$.
Projection onto the resulting planar sector gives
\begin{equation}\label{eq:circular-projection-cases}
 H_n(T,U_n)
 =
 \begin{cases}
   T^2+R_n^2,
      & T\geq R_n\cot\alpha_n,\\[2mm]
   (T\cos\alpha_n+R_n\sin\alpha_n)^2,
      & -R_n\tan\alpha_n<T<R_n\cot\alpha_n,\\[2mm]
   0,
      & T\leq-R_n\tan\alpha_n.
 \end{cases}
\end{equation}
The expressions agree at the boundary points.
Subtracting the middle expression in
\eqref{eq:circular-projection-cases} yields
\eqref{eq:circular-projection-expansion} and
\eqref{eq:circular-projection-error}.

The two positive parts in \eqref{eq:circular-projection-error}
cannot both be nonzero. The first can be nonzero only when $T>0$,
and its square is then bounded by $T^2\sin^2\alpha_n$.
The second can be nonzero only when $T<0$, and its square is bounded
by $T^2\cos^2\alpha_n$.
Thus $\abs{D_n}\leq T^2$, which gives
$\norm{D_n}_2\leq(\E T^4)^{1/2}=\sqrt3$.

The quadratic polynomial
$\cos^2\alpha_n\,T^2+\sin^2\alpha_n\,R_n^2$
belongs to the sum of the zeroth and second Gaussian chaoses.
The function $TR_n$ is odd under $g_n\mapsto-g_n$, so its zeroth
and second components vanish.
Its first component is $(\E R_n)T$: the coefficient of $T$ is
$\E[T^2R_n]=\E R_n$, whereas the coefficients of the coordinates
of $U_n$ vanish by independence and $\E T=0$.
Applying the appropriate orthogonal projections to
\eqref{eq:circular-projection-expansion}, and using their contraction
property in $L^2$, proves
\eqref{eq:circular-chaos-approximation}.

Finally, the function $u\mapsto\norm{u}$ is $1$-Lipschitz.
The Gaussian Poincar\'e inequality therefore gives
$\Var(R_n)\leq1$.
Since $T$ and $R_n$ are independent,
$\norm{T(R_n-\E R_n)}_2=\sqrt{\Var(R_n)}\leq1$.
The last inequality in \eqref{eq:circular-chaos-approximation}
now gives \eqref{eq:circular-higher-chaos-bound}.
\end{proof}

\begin{proposition}[First-chaos dominance for circular cones]
\label{prop:circular-first-chaos}
For the cones in \eqref{eq:circular-cones},
\begin{equation}\label{eq:circular-variance-divergence}
 \tau_n^2\longrightarrow\infty
 \quad\Longleftrightarrow\quad
 n\sin^2(2\alpha_n)\longrightarrow\infty.
\end{equation}
Under these equivalent conditions,
\begin{equation}\label{eq:circular-first-chaos-dominance}
\begin{aligned}
 \frac{\Var(H_{n,1}(g_n))}{\Var(F_n(g_n))}
 &\longrightarrow1,\\
 \tau_n^2
 &\sim\frac n2\sin^2(2\alpha_n),\\
 a_n
 &\sim\frac n4\sin^2(2\alpha_n).
\end{aligned}
\end{equation}
\end{proposition}

\begin{proof}
Both $F_n(g_n)$ and $TR_n$ are centered, and
$\E T^2R_n^2=n$.
The first inequality in \eqref{eq:circular-chaos-approximation}
therefore implies
\begin{equation}\label{eq:circular-remainder-norm}
 \abs{
   \sqrt{\Var(F_n(g_n))}
   -\sqrt n\,\sin(2\alpha_n)
 }
 \leq\sqrt3.
\end{equation}
By \eqref{eq:F-variance} and \eqref{eq:a-comparison},
$2\tau_n^2\leq\Var(F_n(g_n))\leq4\tau_n^2$.
Together with \eqref{eq:circular-remainder-norm}, this proves
\eqref{eq:circular-variance-divergence}.

Assume now that these variances diverge.
Orthogonality of the chaos components and
\eqref{eq:circular-higher-chaos-bound} give
\[
 \frac{\Var(H_{n,1}(g_n))}{\Var(F_n(g_n))}
 =
 1-
 \frac{\Var(F_n(g_n)-H_{n,1}(g_n))}
      {\Var(F_n(g_n))}
 \longrightarrow1.
\]
By \eqref{eq:first-chaos-dominance}, it follows that
$a_n/\tau_n^2\to1/2$.
Consequently,
$\Var(F_n(g_n))=\tau_n^2+2a_n\sim2\tau_n^2$.
On the other hand, \eqref{eq:circular-remainder-norm} gives
$\Var(F_n(g_n))\sim n\sin^2(2\alpha_n)$.
Combining these relations proves the remaining assertions.
\end{proof}

The reason for the first-chaos dominance is visible in
\eqref{eq:circular-chaos-approximation}.
The radius $R_n$ has variance at most one, while
$(\E R_n)^2=n-\Var(R_n)$.
Thus the main part of $TR_n$ is $(\E R_n)T$.
The growing fluctuation comes from the Gaussian coordinate along the
axis of the cone; the additional fluctuation of the radius contributes
only a bounded amount to the variance.

For a fixed opening angle, the correction in
\eqref{eq:circular-projection-error} is exponentially small.
This gives more precise asymptotic formulas.

\begin{proposition}[Circular cones with a fixed opening angle]
\label{prop:circular-fixed-angle}
Suppose that $\alpha_n=\alpha\in(0,\pi/2)$ for every $n$.
Then
\begin{equation}\label{eq:circular-average-derivative}
 \norm{
   A_n-\operatorname{diag}
       (\cos^2\alpha,\sin^2\alpha\,I_n)
 }_{\HS}
 =
 o(n^{-1/2}).
\end{equation}
Moreover, as $n\to\infty$,
\begin{equation}\label{eq:circular-fixed-angle-asymptotics}
\begin{aligned}
 \delta_n
 &=n\sin^2\alpha+\cos^2\alpha+o(1),\\
 a_n
 &=\frac{n+1}{4}\sin^2(2\alpha)+o(1),\\
 \Var(F_n(g_n))
 &=n\sin^2(2\alpha)+o(1),\\
 \tau_n^2
 &=\frac{n-1}{2}\sin^2(2\alpha)+o(1),\\
 4\norm{m_n}^2
 &=\sin^2(2\alpha)(\E R_n)^2+o(1).
\end{aligned}
\end{equation}
\end{proposition}

\begin{proof}
By \eqref{eq:circular-projection-error}, the correction $D_n$
vanishes unless
$\abs{T}\geq R_n\min\{\tan\alpha,\cot\alpha\}$.
Together with $\abs{D_n}\leq T^2$, this gives
\begin{equation}\label{eq:circular-exponential-error}
\begin{aligned}
 \E D_n^2
 &\leq
 \E\bigl[T^4\eee^{T^2/4}\bigr]\,
 \E\eee^{-R_n^2\min\{\tan^2\alpha,\cot^2\alpha\}/4}\\
 &=12\sqrt2
   \left(
     1+\frac12\min\{\tan^2\alpha,\cot^2\alpha\}
   \right)^{-n/2}.
\end{aligned}
\end{equation}
Here we used independence, the Gaussian exponential moment of $T$,
and the fact that $R_n^2$ has the chi-square distribution with $n$
degrees of freedom.
Since $\alpha$ is fixed in $(0,\pi/2)$,
\eqref{eq:circular-exponential-error} implies
$\sqrt n\,\norm{D_n}_2\to0$.

Taking expectations in
\eqref{eq:circular-projection-expansion} proves the formula for
$\delta_n$.
To examine the second chaos, write
$B_n=\operatorname{diag}(\cos^2\alpha,\sin^2\alpha\,I_n)$
for the matrix of the quadratic term in that identity.
The second-chaos component of $D_n$ is
$g_n^\top(A_n-B_n)g_n-\tr(A_n-B_n)$.
Its squared $L^2$-norm is
$2\norm{A_n-B_n}_{\HS}^2$.
Consequently,
$\norm{A_n-B_n}_{\HS}\leq\norm{D_n}_2/\sqrt2$,
which proves \eqref{eq:circular-average-derivative}.

Both $A_n$ and $B_n$ lie between $0$ and $I_{n+1}$.
The eigenvalues of $I_{n+1}-A_n-B_n$ therefore lie in $[-1,1]$.
Using the Cauchy--Schwarz inequality for the Hilbert--Schmidt
inner product, we obtain
\begin{align*}
 \abs{a_n-\tr(B_n-B_n^2)}
 =
 \abs{
   \tr\bigl((I_{n+1}-A_n-B_n)(A_n-B_n)\bigr)
 }
 \leq
 \sqrt{n+1}\,\norm{A_n-B_n}_{\HS}
 =o(1).
\end{align*}
Since
$\tr(B_n-B_n^2)=(n+1)\sin^2(2\alpha)/4$,
this proves the formula for $a_n$.

The first two inequalities in
\eqref{eq:circular-chaos-approximation} give the formulas for
$\Var(F_n(g_n))$ and $4\norm{m_n}^2$.
Indeed, the comparison variables have $L^2$-norms at most $\sqrt n$,
whereas $\sqrt n\,\norm{D_n}_2\to0$.
Thus the resulting errors in their variances tend to zero.
We also use
$\Var(H_{n,1}(g_n))=4\norm{m_n}^2$.
Finally, the identity
$\tau_n^2=\Var(F_n(g_n))-2a_n$
gives the stated formula for $\tau_n^2$.
\end{proof}

We can now describe the remainder obtained by removing the first
component as well as the zeroth and second components.

\begin{corollary}[Removing the first chaos as well]
\label{cor:circular-higher-chaoses}
Suppose that $\alpha_n=\alpha\in(0,\pi/2)$ for every $n$.
Then
\begin{equation}\label{eq:circular-higher-chaos-limit}
 F_n(g_n)-H_{n,1}(g_n)
 \toD
 \frac{\sin(2\alpha)}{\sqrt2}\,Z_1Z_2,
\end{equation}
where $Z_1$ and $Z_2$ are independent standard Gaussian random
variables.
Furthermore,
\begin{equation}\label{eq:circular-higher-chaos-variance}
\begin{aligned}
 \Var\bigl(F_n(g_n)-H_{n,1}(g_n)\bigr)
 &\longrightarrow\frac12\sin^2(2\alpha),\\
 a_n-\norm{m_n}^2
 &\longrightarrow\frac38\sin^2(2\alpha).
\end{aligned}
\end{equation}
\end{corollary}

\begin{proof}
We first recall the limiting behavior of the Gaussian radius.
Since $R_n^2$ is a sum of $n$ independent squared standard Gaussian
variables, the ordinary central limit theorem and the law of large
numbers give
\[
 R_n-\sqrt n
 =
 \frac{(R_n^2-n)/\sqrt n}{R_n/\sqrt n+1}
 \toD \normal(0,1/2).
\]
Moreover,
$0\leq\sqrt n-\E R_n
=\Var(R_n)/(\sqrt n+\E R_n)\leq n^{-1/2}$,
because $\Var(R_n)\leq1$.
Thus centering by $\E R_n$ gives the same limit.

Corollary~\ref{lem:moment}, applied with $q=4$ to the
$1$-Lipschitz function $u\mapsto\norm{u}$, shows that
$\sup_n\E\abs{R_n-\E R_n}^4<\infty$.
The squared centered radii are therefore uniformly integrable.
It follows that
\begin{equation}\label{eq:gaussian-radius-asymptotics}
 R_n-\E R_n\toD\normal(0,1/2),
 \qquad
 \Var(R_n)\longrightarrow\frac12,
 \qquad
 (\E R_n)^2=n-\frac12+o(1).
\end{equation}

The last estimate in \eqref{eq:circular-chaos-approximation},
together with \eqref{eq:circular-exponential-error}, gives
\begin{equation}\label{eq:circular-higher-chaos-approximation}
 F_n(g_n)-H_{n,1}(g_n)
 =
 \sin(2\alpha)\,T(R_n-\E R_n)+o_{L^2}(1).
\end{equation}
For every $n$, the variables $T$ and $R_n-\E R_n$ are independent.
Their joint distribution therefore converges to that of
$(Z_1,Z_2/\sqrt2)$.
Equations~\eqref{eq:gaussian-radius-asymptotics} and
\eqref{eq:circular-higher-chaos-approximation} prove
\eqref{eq:circular-higher-chaos-limit}.

The comparison variable in
\eqref{eq:circular-higher-chaos-approximation} has variance
$\sin^2(2\alpha)\Var(R_n)$.
Its $L^2$-norm is bounded, so the $o_{L^2}(1)$ error also gives
convergence of the variances.
This proves the first limit in
\eqref{eq:circular-higher-chaos-variance}.

Finally, \eqref{eq:circular-fixed-angle-asymptotics} and
\eqref{eq:gaussian-radius-asymptotics} yield
\[
 a_n-\norm{m_n}^2
 =
 \frac{\sin^2(2\alpha)}4
 \bigl(n+1-(\E R_n)^2\bigr)+o(1)
 \longrightarrow\frac38\sin^2(2\alpha).
\]
\end{proof}

The limit in \eqref{eq:circular-higher-chaos-limit} is non-Gaussian.
After division by its own standard deviation, the remainder converges
to $Z_1Z_2$, whose fourth moment is $9$, whereas a standard Gaussian
variable has fourth moment $3$.
Thus, although $\tau_n^2\to\infty$, removing the first component
leaves a remainder with bounded variance and a non-Gaussian limit.
For circular cones, the first chaos contains all the fluctuations
of $F_n(g_n)$ whose variance diverges.

The second limit in
\eqref{eq:circular-higher-chaos-variance} has a geometric
interpretation as well.
By \eqref{eq:affine-least-squares},
$a_n-\norm{m_n}^2$ is the error of the best affine approximation to
the Gaussian projection.
For a fixed opening angle, this error remains bounded even though
the error $a_n$ of the best linear approximation diverges.

\begin{example}[The Lorentz cone and the orthant]
\label{ex:lorentz-orthant-comparison}
For $\alpha_n=\pi/4$, the cone $C_n$ in
\eqref{eq:circular-cones} is the Lorentz cone.
Write $d=n+1$ and suppress the subscript $n$ in the associated
quantities.
Since this cone is self-dual,
\eqref{eq:self-dual-chaos} and \eqref{eq:a-definition} give
$A=I_d/2$, $\delta=d/2$, and $a=d/4$.
By Proposition~\ref{prop:cones-farthest-from-subspaces},
its squared distance to the family of linear subspaces is
$d_{\Pi,2}(C_n,\LS_d)^2=d/2$, the largest possible value.

Proposition~\ref{prop:circular-fixed-angle} gives
$\tau^2=(d-2)/2+o(1)$.
In particular, the universal bound $\tau^2\leq d/2$, which follows
from $\tau^2\leq2a$ and $a\leq d/4$, is asymptotically sharp.
These cones also show that the constants in
$a\geq\tau^2/2$ and $\tau^2\leq2\norm{m}^2$ cannot be improved:
$a/\tau^2\to1/2$ and $\tau^2/(2\norm{m}^2)\to1$.

Table~\ref{tab:orthant-lorentz} compares the Lorentz cone with
the nonnegative orthant in the same ambient dimension.
The orthant values follow from
Example~\ref{ex:orthant-remainder}; the Lorentz-cone values follow
from the results above.
In each column, $H$, $H_1$, $H_2$, and $F$ refer to the cone in
that column.

The two cones have the same statistical dimension, the same
matrix $A$, and the same distance to every linear subspace.
Their first and higher Gaussian chaoses behave differently.
For the orthant, the chaoses of orders at least three contribute
a variance proportional to $d$.
For the Lorentz cone, their total variance converges to $1/2$,
and the first chaos accounts asymptotically for all the variance
of $F(g)$.
\end{example}

\begin{table}[ht]
\centering
\small
\renewcommand{\arraystretch}{1.3}
\begin{tabular}{lcc}
\hline
Quantity
 & Nonnegative orthant
 & Lorentz cone\\
\hline
$\delta$
 & $d/2$
 & $d/2$\\
$A$
 & $I_d/2$
 & $I_d/2$\\
$a$
 & $d/4$
 & $d/4$\\
$d_{\Pi,2}(C,\LS_d)^2$
 & $d/2$
 & $d/2$\\
$\tau^2$
 & $d/4$
 & $(d-2)/2+o(1)$\\
$\Var(H(g))$
 & $5d/4$
 & $3d/2-1+o(1)$\\
$\Var(H_1(g))$
 & $2d/\pi$
 & $d-3/2+o(1)$\\
$\Var(H_2(g))$
 & $d/2$
 & $d/2$\\
$\displaystyle\sum_{q\geq3}\E H_q(g)^2$
 & $\left(\frac34-\frac2\pi\right)d$
 & $\frac12+o(1)$\\
$a-\norm{m}^2$
 & $\left(\frac14-\frac1{2\pi}\right)d$
 & $\frac38+o(1)$\\
$\Var(H_2(g))/\Var(H(g))$
 & $2/5$
 & $\longrightarrow1/3$\\
$\Var(H_1(g))/\Var(F(g))$
 & $8/(3\pi)$
 & $\longrightarrow1$\\
\hline
\end{tabular}
\caption{Comparison of the nonnegative orthant and the Lorentz cone
in $\R^d$. All limits and error terms refer to $d\to\infty$.}
\label{tab:orthant-lorentz}
\end{table}

\subsection{A central limit theorem for the remainder}
\label{subsec:remainder-clt}

We conclude this section with a central limit theorem for the
centered remainder $F(g)$.
Recall that
$F(x)=\|\Pi_C(x)\|^2-x^\top A x$, where
$A=\E D\Pi_C(g)$, and that
$\Var(F(g))=\tau^2+2a$, where
$\tau^2=\Var(V_C)$ and $a=\operatorname{tr}(A-A^2)$.
Proposition~\ref{prop:removal} shows that the gradient of $F$
has Lipschitz constant at most two.
Consequently, when $\tau^2$ tends to infinity, division by
the standard deviation of $F(g)$ makes the Lipschitz constant
of the gradient tend to zero.
We shall also verify that the expected squared norm of the
normalized gradient remains bounded.
These are precisely the conditions needed to apply
Lemma~\ref{lem:gaussian-laplace}.

We include this consequence for completeness.
It describes the fluctuations that remain after the constant
and quadratic Gaussian-chaos components have been removed,
but it is not needed in the subsequent proof of
Theorem~\ref{thm:clt}.

\begin{proposition}[Asymptotic normality of the centered remainder]
\label{prop:F-clt}
Let $C_n\subseteq\R^{d_n}$ be a sequence of nonempty closed
convex cones, and let $g_n$ be a standard Gaussian random
vector in $\R^{d_n}$.
For each $n$, let $F_n$ and $a_n$ be defined by
\eqref{eq:F-definition} and \eqref{eq:a-definition}
for $C=C_n$, and write $\tau_n^2=\Var(V_{C_n})$.
If $\tau_n^2\to\infty$, then
\begin{equation}\label{eq:F-clt}
  \frac{F_n(g_n)}{\sqrt{\tau_n^2+2a_n}}
  \toD\normal(0,1).
\end{equation}
\end{proposition}

\begin{proof}
Since $\tau_n\to\infty$, we may restrict attention to indices
for which $\tau_n\geq2$.
Set $f_n(x)=F_n(x)/\sqrt{\tau_n^2+2a_n}$ for
$x\in\R^{d_n}$.
Proposition~\ref{prop:removal} and
\eqref{eq:F-variance} give
$\E f_n(g_n)=0$ and $\Var(f_n(g_n))=1$.

We next check the two bounds required by
Lemma~\ref{lem:gaussian-laplace}.
Equation~\eqref{eq:F-energy} and
Proposition~\ref{prop:removal} yield
\[
\begin{aligned}
  \E\|\nabla f_n(g_n)\|^2
  &=
  \frac{4a_n}{\tau_n^2+2a_n}
  \leq \frac32,\\
  \operatorname{Lip}(\nabla f_n)
  &\leq
  \frac{2}{\sqrt{\tau_n^2+2a_n}}
  \leq \frac{\sqrt2}{\tau_n}.
\end{aligned}
\]
Here the first inequality uses
$a_n\leq3\tau_n^2/2$, and the second uses
$a_n\geq\tau_n^2/2$; both bounds follow from
\eqref{eq:a-comparison}.
Thus the expected squared norm of the gradient is bounded
independently of $n$, while its Lipschitz constant tends
to zero.

The hypotheses of Lemma~\ref{lem:gaussian-laplace} therefore
hold with $M=\sqrt{3/2}$ and
$\varepsilon=2/\sqrt{\tau_n^2+2a_n}$.
In particular, $\varepsilon\leq1$ for the indices under
consideration.
Since $\Var(f_n(g_n))=1$,
equation~\eqref{eq:gaussian-laplace} gives constants
$r>0$ and $K<\infty$, independent of $n$, such that
\begin{equation}\label{eq:F-normalized-laplace}
  \left|
    \E\eee^{t f_n(g_n)}-\eee^{t^2/2}
  \right|
  \leq
  \frac{2K|t|^3}{\sqrt{\tau_n^2+2a_n}}
  \leq
  \frac{\sqrt2K|t|^3}{\tau_n},
  \qquad |t|\leq r.
\end{equation}
The moment-generating functions therefore converge to that
of a standard Gaussian random variable on a fixed
neighborhood of zero.
Lemma~\ref{lem:laplace-continuity}, stated and proved in
the appendix, yields \eqref{eq:F-clt}.
\end{proof}

The normalization in \eqref{eq:F-clt} is the standard
deviation of the remainder.
By \eqref{eq:F-variance} and \eqref{eq:a-comparison},
it satisfies
$2\tau_n^2\leq\Var(F_n(g_n))\leq4\tau_n^2$.
The central limit theorem for conic intrinsic volumes
requires convergence of
$(V_{C_n}-\delta_n)/\tau_n$, where
$\delta_n=\E V_{C_n}$.
To prove this convergence, the next section returns to
\eqref{eq:steiner-laplace} and incorporates the
quadratic-chaos subtraction into that identity.
The resulting comparison will give convergence of the
Laplace transforms with the required centering
$\delta_n$ and normalization $\tau_n$.

\section{Proof of the central limit theorem for conic intrinsic volumes}
\label{sec:proof-clt}

We now prove Theorem~\ref{thm:clt}.
Let $C_n\subseteq\R^{d_n}$ be a sequence of nonempty
closed convex cones, and write
$V_n=V_{C_n}$, $\delta_n=\E V_n$, and
$\tau_n^2=\Var(V_n)$.
Under the assumption $\tau_n^2\to\infty$, our goal,
already stated in \eqref{eq:laplace-clt-goal}, is to prove
\[
 \E\exp\left(
   t\frac{V_n-\delta_n}{\tau_n}
 \right)
 \longrightarrow\eee^{t^2/2}
 \qquad\text{for every }|t|\leq t_0,
\]
where $t_0>0$ is a universal constant.
By Lemma~\ref{lem:laplace-continuity}, convergence of these
moment-generating functions on a fixed neighborhood of
zero implies convergence in distribution to a standard
Gaussian random variable.

For the estimates leading to this limit, we work with
a single nonempty closed convex cone $C\subseteq\R^d$.
Let $g$ be a standard Gaussian vector in $\R^d$, and
write $\delta=\E V_C$ and $\tau^2=\Var(V_C)$.
We recall the notation from Section~\ref{sec:removal}:
$H(x)=\|\Pi_C(x)\|^2$, $A=\E[D\Pi_C(g)]$, and
$F(x)=H(x)-x^{\mathsf T}Ax$.
Thus $F$ is obtained from $H$ by removing its constant
and second Gaussian-chaos components.
We also write $a=\tr(A-A^2)$.
By \eqref{eq:a-definition} and \eqref{eq:least-squares},
this quantity is the error of the optimal linear
approximation to the projection:
$a=\E\|\Pi_C(g)-Ag\|^2$.

The following conclusions of
Proposition~\ref{prop:removal} will be central to the proof:
\[
  \E F(g)=0,
  \qquad
  \Var(F(g))=\tau^2+2a,
  \qquad
  \E\|\nabla F(g)\|^2=4a,
  \qquad
  \frac12\tau^2\leq a\leq\frac32\tau^2.
\]
In particular, the variance and gradient energy of $F$
are controlled by the intrinsic-volume variance,
independently of the statistical dimension and the
ambient dimension.

Proposition~\ref{prop:F-clt} establishes asymptotic
normality of the remainder after division by its
standard deviation $\sqrt{\tau^2+2a}$.
To obtain the required conclusion for $V_C$, we shall
incorporate the subtraction defining $F$ into the
Gaussian integral in \eqref{eq:steiner-laplace}.
The resulting exponent has first-order term $sF$ and
a second-order term that supplies the correction
needed for the normalization by $\tau$.

We first establish the exact transformed identity and
identify these two terms.
We then obtain remainder and derivative estimates
whose constants are independent of the cone and the
dimension.
After estimating the mean and variance of the
transformed exponent, we apply
Lemma~\ref{lem:gaussian-laplace} to obtain a quantitative
version of the desired Laplace-transform convergence.
The final subsections deduce the central and local
limit theorems and describe subsequential limits
when the intrinsic-volume variance remains bounded.

\subsection{The exponential change of variables}
\label{subsec:exponential-steiner}

Recall that $H(x)=\|\Pi_C(x)\|^2$ and
$F(x)=H(x)-x^{\mathsf T}Ax$, where
$A=\E[D\Pi_C(g)]$.
By \eqref{eq:A-trace}, the matrix $A$ is symmetric and
satisfies $0\preceq A\preceq I_d$ and $\tr A=\delta$.

The Laplace-transform identity
\eqref{eq:steiner-laplace} gives
\begin{equation}\label{eq:centered-steiner-laplace}
 \E\eee^{s(V_C-\delta)}
 =
 \eee^{-s\delta}
 \E\exp\left(
   \frac{1-\eee^{-2s}}2 H(g)
 \right),
 \qquad s\in\R.
\end{equation}
Since $(1-\eee^{-2s})/2=s+O(s^2)$ as $s\to0$,
the exponent inside the Gaussian expectation begins
with $sH$.
Our objective is to rewrite the integral so that
its exponent begins with $sF$.

A linear change of variables modifies the Gaussian
density by a quadratic term.
To obtain the subtraction in
$sF(x)=sH(x)-s x^{\mathsf T}Ax$, we seek a change
whose contribution to the exponent begins with
$-s x^{\mathsf T}Ax$.
We also want its Jacobian to cancel the centering
factor $\eee^{-s\delta}$ in
\eqref{eq:centered-steiner-laplace}.

Both properties hold for the substitution
$y=\eee^{sA}x$.
Indeed, $\eee^{sA}=I_d+sA+O(s^2)$ as $s\to0$,
so the change in the Gaussian density contributes
$-s x^{\mathsf T}Ax$ to first order.
Moreover, the trace identity gives
$\det(\eee^{sA})=\eee^{s\tr A}=\eee^{s\delta}$.
Because $A$ is symmetric, $\eee^{sA}$ is positive
definite and hence invertible for every real $s$.

\begin{proposition}[Exponential change of variables]
\label{prop:transformed-steiner}
For every $s\in\R$,
\begin{equation}\label{eq:transformed-steiner}
 \E\eee^{s(V_C-\delta)}
 =
 \E\eee^{\mathcal E_s(g)},
\end{equation}
where
\begin{equation}\label{eq:Es}
 \mathcal E_s(x)
 =
 \frac{1-\eee^{-2s}}2 H(\eee^{sA}x)
 -\frac12\left(
   \|\eee^{sA}x\|^2-\|x\|^2
 \right),
 \qquad x\in\R^d.
\end{equation}
Moreover, for every fixed $x\in\R^d$,
\begin{equation}\label{eq:Es-second-order}
 \mathcal E_s(x)
 =
 sF(x)-s^2\|\Pi_C(x)-Ax\|^2
 +O_x(|s|^3),
 \qquad s\to0.
\end{equation}
\end{proposition}

\begin{proof}
We first prove \eqref{eq:transformed-steiner}.
Writing the expectation in
\eqref{eq:centered-steiner-laplace} as a Gaussian
integral and substituting $y=\eee^{sA}x$, we obtain
\begin{align*}
 \E\eee^{s(V_C-\delta)}
 &=
 \frac{\eee^{-s\delta}}{(2\pi)^{d/2}}
 \int_{\R^d}
 \exp\left(
   \frac{1-\eee^{-2s}}2 H(y)
   -\frac{\|y\|^2}{2}
 \right)\dd y\\
 &=
 \frac1{(2\pi)^{d/2}}
 \int_{\R^d}
 \exp\left(
   \frac{1-\eee^{-2s}}2 H(\eee^{sA}x)
   -\frac{\|\eee^{sA}x\|^2}{2}
 \right)\dd x\\
 &=
 \E\eee^{\mathcal E_s(g)}.
\end{align*}
In the second equality, the Jacobian
$\det(\eee^{sA})=\eee^{s\delta}$ cancels the
centering factor.
For the last equality, we separate $-\|x\|^2/2$,
the exponent in the standard Gaussian density,
from the remaining terms.
All integrals are finite by
\eqref{eq:steiner-laplace}.

To prove \eqref{eq:Es-second-order}, we expand the
two terms in \eqref{eq:Es} separately.
Recall from \eqref{eq:H-derivatives} and
\eqref{eq:convex-Taylor} that $H$ is continuously
differentiable, $\nabla H(x)=2\Pi_C(x)$, and
\[
 0\leq
 H(x+h)-H(x)-2\langle\Pi_C(x),h\rangle
 \leq\|h\|^2,
 \qquad x,h\in\R^d.
\]
Apply this estimate with $h=(\eee^{sA}-I_d)x$.
Using $0\preceq A\preceq I_d$ and
$\|\Pi_C(x)\|\leq\|x\|$, we obtain
\[
 H(\eee^{sA}x)
 =
 H(x)+2s\langle\Pi_C(x),Ax\rangle
 +O(s^2\|x\|^2).
\]
Here and in the following expansions, the implicit
constants can be chosen universally for $|s|\leq1$.

For the second term, symmetry of $A$ gives
$\|\eee^{sA}x\|^2=x^{\mathsf T}\eee^{2sA}x$.
Expanding the matrix exponential therefore yields
\[
 \frac12\left(
   \|\eee^{sA}x\|^2-\|x\|^2
 \right)
 =
 s x^{\mathsf T}Ax+s^2\|Ax\|^2
 +O(|s|^3\|x\|^2).
\]
Combining these expansions with
$(1-\eee^{-2s})/2=s-s^2+O(|s|^3)$ gives
\[
 \begin{aligned}
  \mathcal E_s(x)
  ={}&
  s\bigl(H(x)-x^{\mathsf T}Ax\bigr)\\
  &+
  s^2\bigl(
    2\langle\Pi_C(x),Ax\rangle-H(x)-\|Ax\|^2
  \bigr)
  +O(|s|^3\|x\|^2)\\
  ={}&
  sF(x)-s^2\|\Pi_C(x)-Ax\|^2
  +O(|s|^3\|x\|^2).
 \end{aligned}
\]
The last equality uses $H(x)=\|\Pi_C(x)\|^2$.
This proves \eqref{eq:Es-second-order}.
The argument requires only the first derivative of $H$
and the remainder bound in \eqref{eq:convex-Taylor}.
\end{proof}

\begin{example}[A linear subspace]
\label{ex:subspace-transformed-steiner}
Let $C=L$ be a linear subspace, and let $P_L$ denote
its orthogonal projection matrix.
As recalled in Example~\ref{ex:subspace-remainder},
we have $A=P_L$ and $H(x)=\|P_Lx\|^2$.
Consequently,
$\eee^{sA}x=\eee^sP_Lx+P_{L^\perp}x$, and orthogonality
gives
\[
 H(\eee^{sA}x)=\eee^{2s}H(x),
 \qquad
 \|\eee^{sA}x\|^2-\|x\|^2
 =
 (\eee^{2s}-1)H(x).
\]
The two terms in \eqref{eq:Es} therefore cancel
exactly, so $\mathcal E_s(x)=0$ for every $s$ and $x$.
Thus the change of variables removes all the
chi-square fluctuations of the squared Gaussian
projection, consistently with
$V_L=\dim L$ almost surely.
\end{example}

For a general cone, \eqref{eq:Es-second-order}
identifies the correction that accompanies the
quadratic-chaos subtraction.
Recall from \eqref{eq:a-definition} that
$\E\|\Pi_C(g)-Ag\|^2=a$.
Also, \eqref{eq:F-variance} gives
$\E F(g)=0$ and $\Var(F(g))=\tau^2+2a$.
Consequently,
\[
  \E\left[-s^2\|\Pi_C(g)-Ag\|^2\right]
   +\frac12\Var(sF(g))
  =
   -as^2+\frac{s^2}{2}(\tau^2+2a)
   =
   \frac{s^2\tau^2}{2}.
\]
These are the leading contributions to the mean and
one half of the variance of the transformed exponent.
At $s=t/\tau$, their sum is exactly $t^2/2$.
Subsection~\ref{subsec:Es-moments} will justify this
approximation and estimate its error.

In particular, the second-order term in
\eqref{eq:Es-second-order} must be retained.
When $s=t/\tau$ and $\tau>0$, its expectation is
$-at^2/\tau^2$.
By \eqref{eq:a-comparison}, the ratio $a/\tau^2$
lies between $1/2$ and $3/2$, so this contribution
does not vanish as $\tau\to\infty$ for fixed $t\ne0$.

It remains to control the error in
\eqref{eq:Es-second-order} by quantities related
to $\tau^2$.
The proof above bounds that error by
$K|s|^3\|x\|^2$ with a universal constant $K$.
Its Gaussian expectation is bounded by
$Kd|s|^3$, which becomes $K|t|^3d/\tau^3$ at
$s=t/\tau$.
The assumption $\tau^2\to\infty$ alone does not
make this bound tend to zero.
The next subsection obtains the remainder and
derivative estimates needed to control the Gaussian
expectation independently of the ambient dimension.

\subsection{Uniform estimates for the transformed exponent}
\label{subsec:uniform-Es}

Recall from \eqref{eq:Es} that
\[
 \mathcal E_s(x)
 =
 \frac{1-\eee^{-2s}}2 H(\eee^{sA}x)
 -\frac12\left(
   \|\eee^{sA}x\|^2-\|x\|^2
 \right),
\]
where $H(x)=\|\Pi_C(x)\|^2$ and
$A=\E[D\Pi_C(g)]$.
We continue to write
$F(x)=H(x)-x^{\mathsf T}Ax$,
$a=\tr(A-A^2)$, and $\tau^2=\Var(V_C)$.

The proof of \eqref{eq:Es-second-order} bounds the
error in the second-order expansion by
$K|s|^3\|x\|^2$.
Although the constant is universal, the Gaussian
expectation of this bound involves the ambient
dimension.
We now replace $\|x\|^2$ by quantities whose
Gaussian expectations are controlled by $\tau^2$.

The appropriate quantities are
$|F(x)|$, $\|\Pi_C(x)-Ax\|^2$, and
$\langle(A-A^2)x,x\rangle$.
Indeed, \eqref{eq:F-variance},
\eqref{eq:a-definition}, and
\eqref{eq:a-comparison} give
\[
 \begin{aligned}
  \E|F(g)|
  &\leq\sqrt{\tau^2+2a}\leq2\tau,\\
  \E\|\Pi_C(g)-Ag\|^2
  &=a,\\
  \E\langle(A-A^2)g,g\rangle
  &=\tr(A-A^2)=a,
 \end{aligned}
\]
where $\tau^2/2\leq a\leq3\tau^2/2$.
Thus all three expectations are controlled
independently of the ambient dimension.
The same comparison will also provide the derivative
bounds needed to apply
Lemma~\ref{lem:gaussian-laplace}.

\begin{proposition}[Uniform estimates for the transformed exponent]
\label{prop:Es-estimates}
There exists a universal constant $K<\infty$ such that,
for every nonempty closed convex cone $C\subseteq\R^d$
and every $|s|\leq1$, the function $\mathcal E_s$
defined in \eqref{eq:Es} satisfies
\begin{equation}\label{eq:Es-expansion}
 \begin{aligned}
  &\left|
    \mathcal E_s(x)-sF(x)
    +s^2\|\Pi_C(x)-Ax\|^2
  \right|\\
  &\qquad\leq
  K|s|^3\left(
    |F(x)|+\|\Pi_C(x)-Ax\|^2
    +\langle(A-A^2)x,x\rangle
  \right),
  \qquad x\in\R^d.
 \end{aligned}
\end{equation}
Moreover, $\mathcal E_s$ is continuously differentiable,
and
\begin{equation}\label{eq:Es-derivative-bounds}
 \operatorname{Lip}(\nabla\mathcal E_s)\leq K|s|,
 \qquad
 \E\|\nabla\mathcal E_s(g)\|^2\leq Ks^2a.
\end{equation}
Here $\operatorname{Lip}(\nabla\mathcal E_s)$ denotes
the global Lipschitz constant of the gradient.
\end{proposition}

\begin{proof}
Throughout the proof, $K$ denotes a universal constant
that may increase from one occurrence to the next.
Fix $|s|\leq1$.
All pointwise estimates below hold for every
$x\in\R^d$.

\parhead{Comparison of the two transformations.}

The main idea is to compare $\eee^{sA}x$ with
$\eee^s\Pi_C(x)+\Pi_{C^\circ}(x)$.
The projection of the latter point is known exactly.

Recall from Proposition~\ref{prop:moreau} that
$x=\Pi_C(x)+\Pi_{C^\circ}(x)$, where
$\Pi_C(x)\in C$, $\Pi_{C^\circ}(x)\in C^\circ$,
and the two vectors are orthogonal.
Conversely, if a vector is written as an orthogonal
sum of an element of $C$ and an element of $C^\circ$,
then these two summands are its projections onto
the respective cones.
Since $\eee^s>0$, this applies to
$\eee^s\Pi_C(x)+\Pi_{C^\circ}(x)$ and gives
\begin{equation}\label{eq:known-projection}
 \Pi_C\bigl(
   \eee^s\Pi_C(x)+\Pi_{C^\circ}(x)
 \bigr)
 =
 \eee^s\Pi_C(x),
 \qquad s\in\R.
\end{equation}

Using Moreau's decomposition again, we can write
the difference between the two transformations as
\begin{equation}\label{eq:exponential-displacement}
 \begin{aligned}
  \eee^{sA}x
  &=
  \eee^s\Pi_C(x)+\Pi_{C^\circ}(x)+z_s(x),\\
  z_s(x)
  &=
  -(\eee^s-1)(\Pi_C(x)-Ax)-\Delta_sx,
 \end{aligned}
\end{equation}
where
$\Delta_s=I_d+(\eee^s-1)A-\eee^{sA}$.
The first term in this displacement contains the
error of the linear approximation $x\mapsto Ax$
to the projection.
The second is the error obtained by replacing
the matrix exponential by its linear interpolant
between the eigenvalues zero and one.
We first estimate this matrix error.

\parhead{The matrix interpolation error.}

For every twice continuously differentiable
function $f$ on $[0,1]$, the error of linear
interpolation satisfies
\begin{equation}\label{eq:linear-interpolation-bound}
 \left|
   (1-\lambda)f(0)+\lambda f(1)-f(\lambda)
 \right|
 \leq
 \frac{\lambda(1-\lambda)}2
 \sup_{0\leq u\leq1}|f''(u)|,
 \qquad 0\leq\lambda\leq1.
\end{equation}
To see this, subtract the affine interpolant from $f$.
The resulting function vanishes at both endpoints.
Subtracting
$\frac12\sup_{[0,1]}|f''|\,\lambda(1-\lambda)$
from this function makes it convex, so its values
are nonpositive.
Adding the same quadratic function makes it
concave, so its values are nonnegative.
These two comparisons prove
\eqref{eq:linear-interpolation-bound}.

To separate the second-order term from the
remaining error, apply
\eqref{eq:linear-interpolation-bound} to
$f(\lambda)=\eee^{s\lambda}-s^2\lambda^2/2$.
Its second derivative is
$f''(\lambda)=s^2(\eee^{s\lambda}-1)$,
whose absolute value is at most $K|s|^3$
for $|s|\leq1$ and $0\leq\lambda\leq1$.
Consequently,
\begin{equation}\label{eq:exponential-interpolation}
 \left|
   1+(\eee^s-1)\lambda-\eee^{s\lambda}
   -\frac{s^2}{2}\lambda(1-\lambda)
 \right|
 \leq
 K|s|^3\lambda(1-\lambda),
 \qquad 0\leq\lambda\leq1.
\end{equation}
In particular,
$|1+(\eee^s-1)\lambda-\eee^{s\lambda}|
\leq Ks^2\lambda(1-\lambda)$.

Recall from \eqref{eq:A-trace} that $A$ is symmetric
and $0\preceq A\preceq I_d$.
We may therefore diagonalize $A$ and apply
\eqref{eq:exponential-interpolation} to each
of its eigenvalues.
This gives
\begin{equation}\label{eq:exponential-matrix-bounds}
 \begin{aligned}
  \|\Delta_sx\|
  &\leq
  Ks^2\sqrt{\langle(A-A^2)x,x\rangle},\\
  \left|
    \langle x,\Delta_sx\rangle
    -\frac{s^2}{2}\langle(A-A^2)x,x\rangle
  \right|
  &\leq
  K|s|^3\langle(A-A^2)x,x\rangle.
 \end{aligned}
\end{equation}
For the first estimate, we also used
$\lambda^2(1-\lambda)^2
\leq\lambda(1-\lambda)/4$ on $[0,1]$.

The factor $\lambda(1-\lambda)$ is essential:
it makes the interpolation error vanish at
both endpoints.
Accordingly, $\Delta_s$ vanishes on the zero
and unit eigenspaces of $A$.
By \eqref{eq:projection-matrix-kernels}, these
eigenspaces are $(\spanop C)^\perp$ and
$C\cap(-C)$, respectively.
Thus these directions make no contribution
to the matrix error.

Combining \eqref{eq:exponential-displacement}
and \eqref{eq:exponential-matrix-bounds}, we obtain
\begin{equation}\label{eq:exponential-displacement-bounds}
 \begin{aligned}
  \|z_s(x)\|
  &\leq
  K\left(
    |s|\|\Pi_C(x)-Ax\|
    +s^2\sqrt{\langle(A-A^2)x,x\rangle}
  \right),\\
  \left|
    \|z_s(x)\|^2-s^2\|\Pi_C(x)-Ax\|^2
  \right|
  &\leq
  K|s|^3\left(
    \|\Pi_C(x)-Ax\|^2
    +\langle(A-A^2)x,x\rangle
  \right).
 \end{aligned}
\end{equation}
The first estimate follows immediately from
$|\eee^s-1|\leq K|s|$.

For the second, expand the squared norm in
\eqref{eq:exponential-displacement}.
The coefficient $(\eee^s-1)^2$ differs from
$s^2$ by at most $K|s|^3$.
The absolute value of the cross term is bounded by
\[
 2|\eee^s-1|\,
  \|\Pi_C(x)-Ax\|\,\|\Delta_sx\|
  \leq
  K|s|^3\left(
    \|\Pi_C(x)-Ax\|^2
    +\langle(A-A^2)x,x\rangle
  \right),
\]
using the first estimate in
\eqref{eq:exponential-matrix-bounds} and
$2uv\leq u^2+v^2$.
Finally,
$\|\Delta_sx\|^2
\leq Ks^4\langle(A-A^2)x,x\rangle$.
Since $|s|\leq1$, this term also has the required
third-order bound.
This proves
\eqref{eq:exponential-displacement-bounds}.

\parhead{The Taylor estimate.}

Recall from \eqref{eq:H-derivatives} and
\eqref{eq:convex-Taylor} that
$\nabla H=2\Pi_C$ and
\[
 0\leq
 H(y+z)-H(y)-2\langle\Pi_C(y),z\rangle
 \leq\|z\|^2,
 \qquad y,z\in\R^d.
\]
We apply this estimate at
$\eee^s\Pi_C(x)+\Pi_{C^\circ}(x)$,
with displacement $z_s(x)$.
By \eqref{eq:known-projection}, the projection
at the expansion point is $\eee^s\Pi_C(x)$.
Thus
\begin{equation}\label{eq:exponential-H-Taylor}
 \begin{aligned}
  H(\eee^{sA}x)
  &=
  \eee^{2s}H(x)
  +2\eee^s\langle\Pi_C(x),z_s(x)\rangle
  +\rho_s(x),\\
  0&\leq\rho_s(x)\leq\|z_s(x)\|^2,
 \end{aligned}
\end{equation}
where $\rho_s(x)$ is the Taylor remainder.

We substitute this formula into \eqref{eq:Es}.
Orthogonality in Moreau's decomposition gives
\[
 \bigl\|
   \eee^s\Pi_C(x)+\Pi_{C^\circ}(x)
 \bigr\|^2-\|x\|^2
 =
 (\eee^{2s}-1)H(x).
\]
The terms that contain neither $z_s(x)$ nor
$\rho_s(x)$ therefore cancel.
Expanding the remaining squared norm gives
\[
  \mathcal E_s(x)
  =
  \left\langle
    -\eee^{-s}\Pi_C(x)-\Pi_{C^\circ}(x),
    z_s(x)
  \right\rangle
  -\frac12\|z_s(x)\|^2
  +\frac{1-\eee^{-2s}}2\rho_s(x).
\]

To substitute for $z_s(x)$, first observe that
$-\eee^{-s}\Pi_C(x)-\Pi_{C^\circ}(x)
=-x+(1-\eee^{-s})\Pi_C(x)$.
We also recall from \eqref{eq:projection-energy}
that $\langle x,\Pi_C(x)\rangle=H(x)$.
Together with the definition of $F$, this gives
\[
 \begin{aligned}
  \langle x,\Pi_C(x)-Ax\rangle
  &=F(x),\\
  2\langle\Pi_C(x),\Pi_C(x)-Ax\rangle
  &=
  F(x)+\|\Pi_C(x)-Ax\|^2
  +\langle(A-A^2)x,x\rangle.
 \end{aligned}
\]
The second equality follows by expanding
$\|\Pi_C(x)-Ax\|^2$.

Using these identities and
$(\eee^s-1)(1-\eee^{-s})=2(\cosh(s)-1)$,
we obtain
\begin{align}
 \mathcal E_s(x)
 ={}&
 \sinh(s)F(x)+\langle x,\Delta_sx\rangle
 \notag\\
 &-(\cosh(s)-1)
 \left(
   \|\Pi_C(x)-Ax\|^2
   +\langle(A-A^2)x,x\rangle
 \right)
 \notag\\
 &-(1-\eee^{-s})
   \langle\Pi_C(x),\Delta_sx\rangle
 \notag\\
 &-\frac12\|z_s(x)\|^2
 +\frac{1-\eee^{-2s}}2\rho_s(x).
 \label{eq:Es-geometric-identity}
\end{align}
In particular, the coefficient of $F(x)$ is
$(\eee^s-1)-(\cosh(s)-1)=\sinh(s)$.

We now estimate the terms in
\eqref{eq:Es-geometric-identity}.
For the mixed term containing $\Pi_C(x)$,
write $\Pi_C(x)=Ax+(\Pi_C(x)-Ax)$.
The matrices $A$ and $\Delta_s$ are diagonal
in the same orthonormal basis.
The scalar estimate following
\eqref{eq:exponential-interpolation}, together
with $0\leq\lambda\leq1$, therefore gives
$|\langle Ax,\Delta_sx\rangle|
\leq Ks^2\langle(A-A^2)x,x\rangle$.
Consequently,
\[
  |\langle\Pi_C(x),\Delta_sx\rangle|
  \leq
  |\langle Ax,\Delta_sx\rangle|
  +\|\Pi_C(x)-Ax\|\,\|\Delta_sx\|\\
  \leq
  Ks^2\left(
    \|\Pi_C(x)-Ax\|^2
    +\langle(A-A^2)x,x\rangle
  \right).
\]
Here we used
\eqref{eq:exponential-matrix-bounds} and
$2uv\leq u^2+v^2$.
Since $|1-\eee^{-s}|\leq K|s|$, multiplication
by the coefficient in
\eqref{eq:Es-geometric-identity} gives the
required third-order bound.

For the Taylor remainder,
\eqref{eq:exponential-H-Taylor} and
\eqref{eq:exponential-displacement-bounds} give
\[
  \left|
    \frac{1-\eee^{-2s}}2\rho_s(x)
  \right|
  \leq K|s|\,\|z_s(x)\|^2
  \leq
  K|s|^3\left(
    \|\Pi_C(x)-Ax\|^2
    +\langle(A-A^2)x,x\rangle
  \right).
\]

It remains to collect the contributions of
the other terms in
\eqref{eq:Es-geometric-identity}.
Use \eqref{eq:exponential-matrix-bounds},
\eqref{eq:exponential-displacement-bounds},
and the scalar estimates
$|\sinh(s)-s|\leq K|s|^3$ and
$|\cosh(s)-1-s^2/2|\leq K|s|^4$.
The contributions of orders one and two are
\[
 \begin{aligned}
  &sF(x)
  +\frac{s^2}{2}\langle(A-A^2)x,x\rangle\\
  &\quad
  -\frac{s^2}{2}\left(
    \|\Pi_C(x)-Ax\|^2
    +\langle(A-A^2)x,x\rangle
  \right)
  -\frac{s^2}{2}\|\Pi_C(x)-Ax\|^2\\
  &\qquad=
  sF(x)-s^2\|\Pi_C(x)-Ax\|^2.
 \end{aligned}
\]
Every error in these replacements is bounded
by a universal constant times
\[
 |s|^3\left(
   |F(x)|+\|\Pi_C(x)-Ax\|^2
   +\langle(A-A^2)x,x\rangle
 \right).
\]
For the fourth-order error from the hyperbolic
cosine, we use $|s|^4\leq|s|^3$.
Together with the bounds for the mixed term
and the Taylor remainder, this proves
\eqref{eq:Es-expansion}.

\parhead{The derivative estimates.}

Since $H$ is continuously differentiable,
so is $\mathcal E_s$.
Differentiating \eqref{eq:Es} and recalling
$\nabla H=2\Pi_C$, we obtain
\begin{equation}\label{eq:Es-gradient}
 \nabla\mathcal E_s(x)
 =
 (I_d-\eee^{2sA})x
 +(1-\eee^{-2s})\eee^{sA}
   \Pi_C(\eee^{sA}x).
\end{equation}

Recall that the operator norm of a matrix $B$ is
$\|B\|_{\op}=\sup_{\|u\|=1}\|Bu\|$.
By Proposition~\ref{prop:projection-properties},
the projection is $1$-Lipschitz:
$\|\Pi_C(u)-\Pi_C(v)\|\leq\|u-v\|$.
Applying this inequality to
\eqref{eq:Es-gradient} yields
\[
 \operatorname{Lip}(\nabla\mathcal E_s)
 \leq
 \|I_d-\eee^{2sA}\|_{\op}
 +|1-\eee^{-2s}|\,
   \|\eee^{sA}\|_{\op}^2
 \leq K|s|.
\]
The last inequality follows from
$0\preceq A\preceq I_d$ and $|s|\leq1$.
In particular, $\|\eee^{sA}\|_{\op}\leq\eee$.
This proves the first bound in
\eqref{eq:Es-derivative-bounds}.

For the expected squared norm of the gradient,
we need a bound involving $a$.
We therefore combine the two terms in
\eqref{eq:Es-gradient} before taking norms.
The comparison
\eqref{eq:exponential-displacement},
the identity \eqref{eq:known-projection},
and the Lipschitz property of $\Pi_C$ give
\begin{equation}\label{eq:exponential-projection-error}
 \left\|
   \Pi_C(\eee^{sA}x)-\eee^s\Pi_C(x)
 \right\|
 \leq\|z_s(x)\|.
\end{equation}

Insert
$\Pi_C(\eee^{sA}x)
=\eee^s\Pi_C(x)
+(\Pi_C(\eee^{sA}x)-\eee^s\Pi_C(x))$
into \eqref{eq:Es-gradient}, and then write
$\Pi_C(x)=Ax+(\Pi_C(x)-Ax)$.
Collecting the terms gives
\begin{align}
 \nabla\mathcal E_s(x)
 ={}&
 2\eee^{sA}
   \bigl(A\sinh(s)-\sinh(sA)\bigr)x
 \notag\\
 &+2\sinh(s)\eee^{sA}(\Pi_C(x)-Ax)
 \notag\\
 &+(1-\eee^{-2s})\eee^{sA}
   \bigl(
     \Pi_C(\eee^{sA}x)-\eee^s\Pi_C(x)
   \bigr),
 \label{eq:Es-gradient-decomposition}
\end{align}
where $\sinh(sA)=(\eee^{sA}-\eee^{-sA})/2$.
Here we used
$I_d-\eee^{2sA}=-2\eee^{sA}\sinh(sA)$
and $(1-\eee^{-2s})\eee^s=2\sinh(s)$.

The first term in
\eqref{eq:Es-gradient-decomposition} is another
matrix interpolation error.
Apply \eqref{eq:linear-interpolation-bound}
to $f(\lambda)=\sinh(s\lambda)$.
Since $f''(\lambda)=s^2\sinh(s\lambda)$,
its absolute value is at most $K|s|^3$ on
$[0,1]$.
It follows that
\[
 |\lambda\sinh(s)-\sinh(s\lambda)|
 \leq
 K|s|^3\lambda(1-\lambda),
 \qquad 0\leq\lambda\leq1.
\]
Diagonalizing $A$, as in the proof of
\eqref{eq:exponential-matrix-bounds}, gives
\[
 \bigl\|
   (A\sinh(s)-\sinh(sA))x
 \bigr\|
 \leq
 K|s|^3
 \sqrt{\langle(A-A^2)x,x\rangle}.
\]

Apply this estimate, together with
\eqref{eq:exponential-projection-error}
and \eqref{eq:exponential-displacement-bounds},
to \eqref{eq:Es-gradient-decomposition}.
Using $\|\eee^{sA}\|_{\op}\leq\eee$,
$|\sinh(s)|\leq K|s|$, and
$|1-\eee^{-2s}|\leq K|s|$, we obtain
\begin{equation}\label{eq:Es-pointwise-gradient}
 \|\nabla\mathcal E_s(x)\|
 \leq
 K\left(
   |s|\|\Pi_C(x)-Ax\|
   +|s|^3\sqrt{\langle(A-A^2)x,x\rangle}
 \right).
\end{equation}
The last term in
\eqref{eq:Es-gradient-decomposition} also produces
a bound proportional to
$s^2\|\Pi_C(x)-Ax\|$.
It is absorbed into the first term on the
right-hand side of
\eqref{eq:Es-pointwise-gradient}, because
$s^2\leq|s|$ for $|s|\leq1$.

Finally, recall from \eqref{eq:a-definition} that
$\E\|\Pi_C(g)-Ag\|^2=a$.
Also,
$\E\langle(A-A^2)g,g\rangle
=\tr(A-A^2)=a$.
Squaring \eqref{eq:Es-pointwise-gradient},
using $(u+v)^2\leq2u^2+2v^2$, and taking
Gaussian expectations therefore gives
\[
 \begin{aligned}
  \E\|\nabla\mathcal E_s(g)\|^2
  &\leq
  K\left(
    s^2\E\|\Pi_C(g)-Ag\|^2
    +s^6\E\langle(A-A^2)g,g\rangle
  \right)\\
  &=
  K(s^2+s^6)a
  \leq2Ks^2a.
 \end{aligned}
\]
Increasing the universal constant proves the
second bound in
\eqref{eq:Es-derivative-bounds}.
\end{proof}

The comparison uses Moreau's decomposition, the
Lipschitz continuity of the projection, and the
Taylor estimate for the continuously differentiable
function $H$.
In the next subsection, we apply
\eqref{eq:Es-expansion} with $s=t/\tau$ to estimate
the mean and variance of $\mathcal E_{t/\tau}(g)$.

\subsection{The mean and variance of the transformed exponent}
\label{subsec:Es-moments}

We now estimate the mean and variance of
$\mathcal E_{t/\tau}(g)$.
Let $C\subseteq\R^d$ be a nonempty closed convex cone,
and let $g$ be a standard Gaussian vector in $\R^d$.
Recall that
$F(x)=\|\Pi_C(x)\|^2-x^{\mathsf T}Ax$, where
$A=\E[D\Pi_C(g)]$.
We write $\tau^2=\Var(V_C)$ and recall from
\eqref{eq:a-definition} that
$a=\tr(A-A^2)=\E\|\Pi_C(g)-Ag\|^2$.
Furthermore, \eqref{eq:F-variance} and
\eqref{eq:a-comparison} give
$\E F(g)=0$, $\Var(F(g))=\tau^2+2a$, and
$\tau^2/2\leq a\leq3\tau^2/2$.

We shall work with $\tau\geq2$ and $|t|\leq1$.
The second-order term in \eqref{eq:Es-expansion},
evaluated at $s=t/\tau$, is
$-t^2\|\Pi_C(g)-Ag\|^2/\tau^2$.
Its expectation is $-at^2/\tau^2$.
We shall show that its centered part has
$L^2$ norm at most a universal constant times
$t^2/\tau$.
Thus its random fluctuations become small as
$\tau\to\infty$, while its expectation must be
retained.

This leads us to compare $\mathcal E_{t/\tau}(g)$
with $tF(g)/\tau-at^2/\tau^2$.
The following proposition gives separate estimates
for the mean and variance.
Their combination identifies the exponent that
will enter the Laplace-transform approximation
in the next subsection.

\begin{proposition}[Mean and variance of the transformed exponent]
\label{prop:Es-moments}
There exists a universal constant $K<\infty$ such that,
for every nonempty closed convex cone $C\subseteq\R^d$
with $\tau\geq2$ and every $|t|\leq1$,
\begin{align}
 \left|
   \E\mathcal E_{t/\tau}(g)
   +\frac{at^2}{\tau^2}
 \right|
 &\leq
 \frac{K|t|^3}{\tau},
 \label{eq:Es-mean}\\
 \left|
   \Var(\mathcal E_{t/\tau}(g))
   -\frac{t^2}{\tau^2}(\tau^2+2a)
 \right|
 &\leq
 \frac{K|t|^3}{\tau}.
 \label{eq:Es-variance}
\end{align}
Consequently,
\begin{equation}\label{eq:Es-mean-variance}
 \left|
   \E\mathcal E_{t/\tau}(g)
   +\frac12\Var(\mathcal E_{t/\tau}(g))
   -\frac{t^2}{2}
 \right|
 \leq
 \frac{K|t|^3}{\tau}.
\end{equation}
\end{proposition}

\begin{proof}
Throughout the proof, $K$ denotes a universal constant
that may increase from one occurrence to the next.

\parhead{Moment bounds for the terms in
\eqref{eq:Es-expansion}.}

We first control the fluctuations of
$\|\Pi_C(g)-Ag\|^2$.
Recall the Gaussian Poincar\'e inequality
\eqref{eq:poincare}:
$\Var(f(g))\leq\E\|\nabla f(g)\|^2$
whenever $f$ and its first weak derivatives
belong to $L^2(\R^d,\gamma_d)$.

Since $\Pi_C$ is Lipschitz, the function
$x\mapsto\|\Pi_C(x)-Ax\|^2$ is locally Lipschitz.
Its weak gradient is
\[
 \nabla\bigl(\|\Pi_C(x)-Ax\|^2\bigr)
 =
 2(D\Pi_C(x)-A)(\Pi_C(x)-Ax)
 \quad\text{almost everywhere}.
\]
Here we used the symmetry of $D\Pi_C(x)$ and $A$.
By \eqref{eq:projection-derivative} and
\eqref{eq:A-trace}, both matrices lie between
zero and $I_d$ in the symmetric-matrix order.
Their difference therefore lies between $-I_d$
and $I_d$, so
$\|D\Pi_C(x)-A\|_{\op}\leq1$
almost everywhere.

The function has at most quadratic growth,
and its weak gradient has at most linear growth.
Both are therefore square integrable under
Gaussian measure.
Applying \eqref{eq:poincare} and recalling
$\E\|\Pi_C(g)-Ag\|^2=a$, we obtain
\begin{equation}\label{eq:squared-error-variance}
 \begin{aligned}
  \Var\bigl(\|\Pi_C(g)-Ag\|^2\bigr)
  &\leq
  4\E\bigl\|
    (D\Pi_C(g)-A)(\Pi_C(g)-Ag)
  \bigr\|^2\\
  &\leq
  4\E\|\Pi_C(g)-Ag\|^2
  =4a.
 \end{aligned}
\end{equation}

We also need to control
$\langle(A-A^2)g,g\rangle$.
Its expectation is $\tr(A-A^2)=a$.
Diagonalizing the symmetric matrix $A-A^2$,
and using the independence of the Gaussian
coordinates together with $\Var(g_i^2)=2$,
gives
\begin{equation}\label{eq:matrix-error-variance}
 \Var\bigl(\langle(A-A^2)g,g\rangle\bigr)
 =
 2\tr\bigl((A-A^2)^2\bigr)
 \leq\frac a2.
\end{equation}
For the inequality, recall that the eigenvalues
of $A$ belong to $[0,1]$.
The eigenvalues of $A-A^2$ consequently belong
to $[0,1/4]$, and hence
$\tr((A-A^2)^2)\leq\tr(A-A^2)/4=a/4$.

For any square-integrable random variable $X$,
the triangle inequality gives
$\|X\|_{L^2}\leq|\E X|+\sqrt{\Var(X)}$.
Using this observation,
\eqref{eq:squared-error-variance},
\eqref{eq:matrix-error-variance}, and the
identities for $F$, we obtain
\begin{equation}\label{eq:Es-coefficient-L2}
 \begin{aligned}
  \|F(g)\|_{L^2}
  &=
  \sqrt{\tau^2+2a}\leq2\tau,\\
  \left\|\,\|\Pi_C(g)-Ag\|^2\,\right\|_{L^2}
  &\leq a+2\sqrt a,\\
  \bigl\|\langle(A-A^2)g,g\rangle\bigr\|_{L^2}
  &\leq a+\sqrt{a/2}.
 \end{aligned}
\end{equation}
The first line uses $\E F(g)=0$ and
$a\leq3\tau^2/2$.

\parhead{The mean and an $L^2$ approximation.}

Since $\tau\geq2$ and $|t|\leq1$, the estimate
\eqref{eq:Es-expansion} applies with $s=t/\tau$.
It gives
\[
 \begin{aligned}
  &\left|
    \mathcal E_{t/\tau}(g)
    -\frac t\tau F(g)
    +\frac{t^2}{\tau^2}\|\Pi_C(g)-Ag\|^2
  \right|\\
  &\qquad\leq
  \frac{K|t|^3}{\tau^3}
  \left(
    |F(g)|+\|\Pi_C(g)-Ag\|^2
    +\langle(A-A^2)g,g\rangle
  \right).
 \end{aligned}
\]
Taking $L^2$ norms and using
\eqref{eq:Es-coefficient-L2}, we obtain
\begin{equation}\label{eq:Es-Taylor-L2}
  \left\|
    \mathcal E_{t/\tau}(g)
    -\frac t\tau F(g)
    +\frac{t^2}{\tau^2}\|\Pi_C(g)-Ag\|^2
  \right\|_{L^2}
  \leq
  \frac{K|t|^3}{\tau^3}
  \bigl(\tau+a+\sqrt a\bigr)
  \leq\frac{K|t|^3}{\tau}.
\end{equation}
The last inequality follows from
$a\leq3\tau^2/2$ and $\tau\geq2$.

The absolute value of the expectation of a
random variable is at most its $L^2$ norm.
Since $\E F(g)=0$ and
$\E\|\Pi_C(g)-Ag\|^2=a$,
equation~\eqref{eq:Es-Taylor-L2} therefore proves
\eqref{eq:Es-mean}.
Using the exact expectation of the second-order
term preserves the cubic factor $|t|^3$
in this estimate.

For the variance, we replace the random
second-order term by its expectation.
Equation~\eqref{eq:squared-error-variance}
gives
$\bigl\|\,\|\Pi_C(g)-Ag\|^2-a\,\bigr\|_{L^2}
\leq2\sqrt a$.
Combining this with
\eqref{eq:Es-Taylor-L2} yields
\begin{equation}\label{eq:Es-L2-approximation}
 \begin{aligned}
  &\left\|
    \mathcal E_{t/\tau}(g)
    -\frac t\tau F(g)
    +\frac{at^2}{\tau^2}
  \right\|_{L^2}\leq
  \frac{K|t|^3}{\tau}
  +\frac{t^2}{\tau^2}
   \left\|\,\|\Pi_C(g)-Ag\|^2-a\,\right\|_{L^2}\\
  &\qquad\leq
  \frac{K|t|^3}{\tau}
  +\frac{2t^2\sqrt a}{\tau^2}
  \leq\frac{Kt^2}{\tau}.
 \end{aligned}
\end{equation}
The last inequality uses
$a\leq3\tau^2/2$ and $|t|\leq1$.
This error is quadratic in $t$.
In the variance comparison below, it will be
combined with the standard deviation of
$tF(g)/\tau$, which is at most $2|t|$.
This supplies the additional factor $|t|$
needed in \eqref{eq:Es-variance}.

\parhead{The variance.}

For any square-integrable random variables
$X$ and $Y$, we have
\[
 |\Var(X)-\Var(Y)|
 \leq
 \|X-Y\|_{L^2}
 \left(
   2\sqrt{\Var(Y)}+\|X-Y\|_{L^2}
 \right).
\]
Indeed, expanding the variance of $Y+(X-Y)$
and applying Cauchy--Schwarz gives
\[
 \begin{aligned}
  |\Var(X)-\Var(Y)|
  &\leq
  2\sqrt{\Var(Y)}\sqrt{\Var(X-Y)}
  +\Var(X-Y)\\
  &\leq
  2\sqrt{\Var(Y)}\,\|X-Y\|_{L^2}
  +\|X-Y\|_{L^2}^2.
 \end{aligned}
\]

Apply this inequality with
$X=\mathcal E_{t/\tau}(g)+at^2/\tau^2$
and $Y=tF(g)/\tau$.
Adding a constant does not change the variance
of $X$.
By \eqref{eq:F-variance} and
\eqref{eq:a-comparison},
$\Var(Y)=t^2(\tau^2+2a)/\tau^2\leq4t^2$.
Consequently,
\eqref{eq:Es-L2-approximation} gives
\[
  \left|
    \Var(\mathcal E_{t/\tau}(g))
    -\frac{t^2}{\tau^2}(\tau^2+2a)
  \right|\leq
  \frac{Kt^2}{\tau}
  \left(
    4|t|+\frac{Kt^2}{\tau}
  \right)
  \leq\frac{K|t|^3}{\tau}.
\]
In the last step, we used $|t|\leq1$ and
$\tau\geq2$.
This proves \eqref{eq:Es-variance}.

Finally, the leading contributions to the mean
and one half of the variance satisfy
\[
 -\frac{at^2}{\tau^2}
 +\frac{t^2}{2\tau^2}(\tau^2+2a)
 =
 \frac{t^2}{2}.
\]
Combining \eqref{eq:Es-mean} and
\eqref{eq:Es-variance} therefore proves
\eqref{eq:Es-mean-variance}, after increasing
the universal constant.
\end{proof}

The contribution of $2a$ to the variance of $F(g)$
is thus compensated by the mean of the second-order
correction in the transformed exponent.
To obtain the desired Laplace-transform estimate,
it remains to compare
$\E\eee^{\mathcal E_{t/\tau}(g)}$ with
$\exp\bigl(
 \E\mathcal E_{t/\tau}(g)
 +\frac12\Var(\mathcal E_{t/\tau}(g))
\bigr)$.
The next subsection makes this comparison using
Lemma~\ref{lem:gaussian-laplace} and the derivative
bounds \eqref{eq:Es-derivative-bounds}.

\subsection{The quantitative Laplace-transform estimate}
\label{subsec:quantitative-laplace}

We now combine the results of the preceding three
subsections to estimate the Laplace transform of
$(V_C-\delta)/\tau$.
Let $C\subseteq\R^d$ be a nonempty closed convex cone,
let $g$ be a standard Gaussian vector in $\R^d$,
and write $\delta=\E V_C$ and $\tau^2=\Var(V_C)$.
Recall from \eqref{eq:transformed-steiner} that
\[
 \E\eee^{s(V_C-\delta)}
 =
 \E\eee^{\mathcal E_s(g)},
 \qquad s\in\R.
\]
Thus it suffices to estimate the exponential moment
of the transformed exponent.

Two properties of this exponent will be used.
First, \eqref{eq:Es-derivative-bounds} gives
\[
 \operatorname{Lip}(\nabla\mathcal E_s)\leq K|s|,
 \qquad
 \E\|\nabla\mathcal E_s(g)\|^2\leq Ks^2a,
 \qquad |s|\leq1,
\]
where $a=\tr(A-A^2)$ and $A=\E[D\Pi_C(g)]$.
Recall also from \eqref{eq:a-comparison} that
$a\leq3\tau^2/2$.
Second, \eqref{eq:Es-mean-variance} shows that
\[
 \left|
   \E\mathcal E_{t/\tau}(g)
   +\frac12\Var(\mathcal E_{t/\tau}(g))
   -\frac{t^2}{2}
 \right|
 \leq\frac{K|t|^3}{\tau},
 \qquad \tau\geq2,\quad |t|\leq1.
\]
The remaining step is to approximate the exponential
moment by that of a Gaussian random variable with
the same mean and variance.
We shall obtain this approximation from
Lemma~\ref{lem:gaussian-laplace}, using the derivative
bounds above.

\begin{theorem}[A local Laplace-transform estimate]
\label{thm:laplace-main}
There exist universal constants $t_0>0$ and
$K_0<\infty$ with the following property.
Let $C\subseteq\R^d$ be a nonempty closed convex cone,
and write $\delta=\E V_C$ and $\tau^2=\Var(V_C)$.
If $\tau\geq2$, then
\begin{equation}\label{eq:laplace-main}
 \left|
   \E\exp\left(
     t\frac{V_C-\delta}{\tau}
   \right)
   -\eee^{t^2/2}
 \right|
 \leq
 \frac{K_0|t|^3}{\tau},
 \qquad |t|\leq t_0.
\end{equation}
The constants $t_0$ and $K_0$ do not depend on
the cone or the ambient dimension.
\end{theorem}

\begin{proof}
We first recall the part of
Lemma~\ref{lem:gaussian-laplace} that we shall use.
Taking $M=1$ in that lemma gives universal constants
$r>0$ and $K_1<\infty$ with the following property.
If $f\in C^1(\R^d)$ satisfies
$\E f(g)=0$, $\E\|\nabla f(g)\|^2\leq1$,
and $\operatorname{Lip}(\nabla f)\leq\varepsilon$
for some $0\leq\varepsilon\leq1$, then
\[
 \left|
   \E\eee^{u f(g)}
   -\eee^{u^2\Var(f(g))/2}
 \right|
 \leq
 K_1\varepsilon|u|^3,
 \qquad |u|\leq r.
\]
We shall apply this estimate with
$\varepsilon=1/\tau$.

\parhead{Normalization and Gaussian approximation.}

Let $K_*$ be a universal constant for which
\eqref{eq:Es-derivative-bounds} holds.
Choose a universal constant
$L\geq\max\{1,K_*,\sqrt{3K_*/2}\}$,
and set $t_0=\min\{1,r/L\}$.
These choices are independent of the cone and
the ambient dimension.

The estimate \eqref{eq:laplace-main} is immediate
when $t=0$.
Fix $0<|t|\leq t_0$.
We divide the centered exponent by $Lt$ and
apply the preceding inequality at parameter $u=Lt$.
This normalization will make the derivative bounds
independent of $t$, while the parameter $Lt$
will supply the cubic factor $|t|^3$.
Accordingly, define
\[
 f(x)
 =
 \frac{
   \mathcal E_{t/\tau}(x)
   -\E\mathcal E_{t/\tau}(g)
 }{Lt},
 \qquad x\in\R^d.
\]
The function $f$ is continuously differentiable
and satisfies $\E f(g)=0$.

Since $|t|\leq1$ and $\tau\geq2$, the derivative
estimates \eqref{eq:Es-derivative-bounds} apply
with $s=t/\tau$.
Together with $a\leq3\tau^2/2$, they give
\[
 \begin{aligned}
  \E\|\nabla f(g)\|^2
  &\leq
  \frac{K_*a}{L^2\tau^2}
  \leq
  \frac{3K_*}{2L^2}
  \leq1,\\
  \operatorname{Lip}(\nabla f)
  &=
  \frac{\operatorname{Lip}
    (\nabla\mathcal E_{t/\tau})}{L|t|}
  \leq
  \frac{K_*}{L\tau}
  \leq\frac1\tau.
 \end{aligned}
\]
Thus the hypotheses of
Lemma~\ref{lem:gaussian-laplace} hold with
$M=1$ and $\varepsilon=1/\tau\leq1$.
Moreover, $|Lt|\leq r$ by the choice of $t_0$.

Applying \eqref{eq:gaussian-laplace} at parameter
$u=Lt$, we obtain
\begin{equation}\label{eq:centered-Es-laplace}
  \Bigg|
   \E\exp\left(
     \mathcal E_{t/\tau}(g)
     -\E\mathcal E_{t/\tau}(g)
   \right)
   -
   \exp\left(
     \frac12\Var(\mathcal E_{t/\tau}(g))
   \right)
  \Bigg|
  \leq
  \frac{K_1L^3|t|^3}{\tau}.
\end{equation}
Here we used
$(Lt)f(g)=\mathcal E_{t/\tau}(g)
-\E\mathcal E_{t/\tau}(g)$ and
$(Lt)^2\Var(f(g))
=\Var(\mathcal E_{t/\tau}(g))$.
Both identities hold for either sign of $t$.

\parhead{Restoring the mean.}

In the remainder of the proof, $K$ denotes a
universal constant that may increase from one
occurrence to the next.
Recall from \eqref{eq:Es-mean} that
\[
 \left|
   \E\mathcal E_{t/\tau}(g)
   +\frac{at^2}{\tau^2}
 \right|
 \leq\frac{K|t|^3}{\tau}.
\]
Since $a/\tau^2\leq3/2$, it follows that
\[
 \left|\E\mathcal E_{t/\tau}(g)\right|
 \leq
 \frac{at^2}{\tau^2}
 +\frac{K|t|^3}{\tau}
 \leq
 \frac32t^2+\frac{K|t|^3}{\tau}
 \leq K.
\]
The final bound is uniform because $|t|\leq1$
and $\tau\geq2$.

Multiplying \eqref{eq:centered-Es-laplace}
by $\exp(\E\mathcal E_{t/\tau}(g))$ therefore
changes the error bound by at most a universal
factor.
Absorbing the fixed constants $K_1$ and $L$
into $K$, we obtain
\begin{equation}\label{eq:Es-laplace-approximation}
  \Bigg|
   \E\eee^{\mathcal E_{t/\tau}(g)}
   -
   \exp\left(
     \E\mathcal E_{t/\tau}(g)
     +\frac12\Var(\mathcal E_{t/\tau}(g))
   \right)
  \Bigg|
  \leq
  \frac{K|t|^3}{\tau}.
\end{equation}

\parhead{Completion of the proof.}

By \eqref{eq:Es-mean-variance},
\[
 \left|
   \E\mathcal E_{t/\tau}(g)
   +\frac12\Var(\mathcal E_{t/\tau}(g))
   -\frac{t^2}{2}
 \right|
 \leq\frac{K|t|^3}{\tau}.
\]
Both
$\E\mathcal E_{t/\tau}(g)
+\frac12\Var(\mathcal E_{t/\tau}(g))$
and $t^2/2$ consequently lie in a fixed bounded
interval, independently of $C$, $d$, $\tau$,
and $t$ in the ranges under consideration.
The mean value theorem for the exponential
function therefore gives
\[
  \left|
    \exp\left(
      \E\mathcal E_{t/\tau}(g)
      +\frac12\Var(\mathcal E_{t/\tau}(g))
    \right)
    -\eee^{t^2/2}
  \right|
  \leq\frac{K|t|^3}{\tau}.
\]
Combining this estimate with
\eqref{eq:Es-laplace-approximation}, and then
using \eqref{eq:transformed-steiner} with
$s=t/\tau$, yields
\[
  \left|
    \E\eee^{t(V_C-\delta)/\tau}
    -\eee^{t^2/2}
  \right|
  =
  \left|
    \E\eee^{\mathcal E_{t/\tau}(g)}
    -\eee^{t^2/2}
  \right|\\
  \leq\frac{K|t|^3}{\tau}.
\]
Choosing a universal constant $K_0$ that dominates
the final constant proves
\eqref{eq:laplace-main}.
\end{proof}

The estimate \eqref{eq:laplace-main} completes the
comparison of Laplace transforms.
Its uniformity in the cone and the ambient
dimension allows it to be applied to any sequence
whose intrinsic-volume variance tends to infinity.
In the next subsection, we use it to prove
Theorem~\ref{thm:clt} and then deduce the local
central limit theorem.

\subsection{Proofs of the central and local limit theorems}
\label{subsec:proof-central-local-clt}

We now deduce Theorems~\ref{thm:clt} and
\ref{thm:local-clt} from the quantitative
Laplace-transform estimate established in the
preceding subsection.
Recall that $C_n\subseteq\R^{d_n}$ is a sequence
of nonempty closed convex cones and that
$V_n=V_{C_n}$, $\delta_n=\E V_n$, and
$\tau_n^2=\Var(V_n)$.

\begin{proof}[Proof of Theorem~\ref{thm:clt}]
Since $\tau_n^2\to\infty$, we have $\tau_n\geq2$
for all sufficiently large $n$.
Theorem~\ref{thm:laplace-main} gives universal constants
$t_0>0$ and $K_0<\infty$ such that
\[
 \left|
   \E\exp\left(t\frac{V_n-\delta_n}{\tau_n}\right)
   -\eee^{t^2/2}
 \right|
 \leq\frac{K_0|t|^3}{\tau_n},
 \qquad |t|\leq t_0.
\]
Consequently, the moment-generating functions of the
standardized variables converge to that of a standard
Gaussian random variable on a fixed neighborhood of
zero.
By Lemma~\ref{lem:laplace-continuity}, stated and
proved in the appendix, this implies
\[
 \frac{V_n-\delta_n}{\tau_n}\toD\normal(0,1),
\]
which proves \eqref{eq:main-clt}.
\end{proof}

The passage from the central limit theorem to the
local central limit theorem uses two additional
properties of conic intrinsic volumes:
their log-concavity and the absence of gaps in
their support.

\begin{proof}[Proof of Theorem~\ref{thm:local-clt}]
A classical result of Bender
\cite[Lemma~2]{Bender73} states that for integer-valued random variables $X_n$
with means $\mu_n$ and standard deviations
$\sigma_n\to\infty$, whose probability mass functions
are log-concave and have consecutive supports,
the central limit theorem
\[
 \frac{X_n-\mu_n}{\sigma_n}\toD\normal(0,1)
\]
implies the uniform local central limit theorem
\[
 \sup_{k\in\mathbb Z}
 \left|
   \sigma_n\PP\{X_n=k\}
   -
   \frac{1}{\sqrt{2\pi}}
   \exp\left(-\frac{(k-\mu_n)^2}{2\sigma_n^2}\right)
 \right|
 \longrightarrow0.
\]

We apply this result to $X_n=V_n$, with
$\mu_n=\delta_n$ and $\sigma_n=\tau_n$.
By assumption, $\tau_n\to\infty$, and the required
central limit theorem holds by Theorem~\ref{thm:clt}.

Recall from \eqref{eq:intrinsic-variable} that
$\PP\{V_n=k\}=\upsilon_k(C_n)$.
The recent result of Fu and Wang \cite{FuWang26},
recalled in \eqref{eq:log-concavity}, gives
\[
 \upsilon_k(C_n)^2
 \geq
 \upsilon_{k-1}(C_n)\upsilon_{k+1}(C_n),
 \qquad k\in\mathbb Z,
\]
where the intrinsic volumes are extended by zero
outside $\{0,\ldots,d_n\}$.
Thus the probability mass functions of $V_n$
are log-concave.

Furthermore, by \cite[Theorem~4.5.1]{Schneider22},
for every nonempty closed convex cone $C$,
\[
 \upsilon_k(C)>0
 \quad\Longleftrightarrow\quad
 \dim(C\cap(-C))\leq k\leq\dim C,
 \qquad k\in\mathbb Z.
\]
In particular, the support of $V_n$ consists of
consecutive integers.

All the hypotheses of Bender's local limit theorem
are therefore satisfied.
Its conclusion is precisely
\eqref{eq:local-clt}.
\end{proof}

Theorems~\ref{thm:clt} and \ref{thm:local-clt}
establish asymptotic normality and its local version
under the assumption that the intrinsic-volume
variance tends to infinity.
In the next subsection, we investigate what happens
when this variance remains bounded and establish
the necessity of the variance-divergence condition
for asymptotic normality.

\subsection{Subsequential limits at bounded variance}
\label{subsec:proof-bounded-variance}

We now prove Proposition~\ref{prop:bounded-variance},
which describes the possible subsequential limits of
the standardized intrinsic-volume random variables
when their variances remain bounded.
The proof also establishes a universal positive lower
bound on the intrinsic-volume variance of a cone that
is not a linear subspace.

\begin{proof}[Proof of Proposition~\ref{prop:bounded-variance}]
Recall that $V_n=V_{C_n}$, $\delta_n=\E V_n$, and
$\tau_n^2=\Var(V_n)$.
By assumption, $0<\tau_n^2\leq M$ for every $n$.
Since $\tau_n^2>0$, none of the cones $C_n$ is a
linear subspace.
The conic Gauss--Bonnet identity \eqref{eq:gauss-bonnet}
therefore gives
\[
 \PP\{V_n\text{ is even}\}
 =
 \PP\{V_n\text{ is odd}\}
 =
 \frac12.
\]

Let $V_n'$ be an independent copy of $V_n$.
Since both variables are integer-valued, their squared
difference is at least one whenever they have different
parity.
By independence and the parity balance,
\[
  2\tau_n^2
  =
  \E[(V_n-V_n')^2]
  \geq
  \PP\{V_n-V_n'\text{ is odd}\}
  =
  \frac12.
\]
Together with the assumed upper bound, this gives
\begin{equation}\label{eq:bounded-variance-range}
 \frac14\leq\tau_n^2\leq M.
\end{equation}

To obtain a tight sequence of integer-valued random
variables, we center $V_n$ by the integer part of its
expectation rather than by the expectation itself.
More precisely, put
\[
 k_n=\lfloor\delta_n\rfloor,
 \qquad
 W_n=V_n-k_n,
 \qquad
 \alpha_n=\delta_n-k_n\in[0,1).
\]
Then $W_n$ is integer-valued, $\E W_n=\alpha_n$,
and $\Var(W_n)=\tau_n^2$.
In particular,
$\E W_n^2=\tau_n^2+\alpha_n^2\leq M+1$.
Thus the second moments of $W_n$ are uniformly bounded,
which implies that the sequence $(W_n)$ is tight.

Starting from an arbitrary subsequence, tightness of
$(W_n)$, boundedness of $(\alpha_n)$, and
\eqref{eq:bounded-variance-range} allow us to extract
a further subsequence along which
\[
 W_n\toD W,
 \qquad
 \alpha_n\to\alpha,
 \qquad
 \tau_n\to\tau,
\]
where $\alpha\in[0,1]$ and $\tau\in[1/2,\sqrt M]$.
Since every $W_n$ is integer-valued and $\mathbb Z$
is closed in $\R$, the Portmanteau theorem implies
that $W$ is also integer-valued.
As $\tau>0$, Slutsky's theorem gives, along this
further subsequence,
\[
 \frac{V_n-\delta_n}{\tau_n}
 =
 \frac{W_n-\alpha_n}{\tau_n}
 \toD
 \frac{W-\alpha}{\tau}.
\]

It remains to verify the parity balance of $W$.
Since $W_n$ differs from $V_n$ by a deterministic
integer, the even and odd values of $W_n$ also
have equal probabilities.
Equivalently, $\E\cos(\pi W_n)=0$.
Since the function $x\mapsto\cos(\pi x)$ is bounded
and continuous, weak convergence gives
$\E\cos(\pi W)=0$.
As $W$ is integer-valued, this is equivalent to
\[
 \PP\{W\text{ is even}\}
 =
 \PP\{W\text{ is odd}\}
 =
 \frac12.
\]
In particular, $W$ cannot be almost surely constant.
Hence $(W-\alpha)/\tau$ has a discrete, nondegenerate
distribution.

Finally, if a subsequence of the standardized variables
converges weakly, the argument above applied to that
subsequence yields a further subsequence with a limit
of the stated form.
By uniqueness of the weak limit, the original
subsequential limit has this form as well.
This proves all the assertions of the proposition.
\end{proof}

The variance-gap argument at the beginning of the proof
applies to every nonempty closed convex cone that is not
a linear subspace, irrespective of its ambient dimension.
It shows that
\begin{equation}\label{eq:intrinsic-volume-variance-gap}
 \Var(V_C)\geq\frac14
 \qquad
 \text{if $C$ is not a linear subspace}.
\end{equation}
The bound is sharp, since it is attained by a ray
and by a halfspace.

Proposition~\ref{prop:bounded-variance} also establishes
the necessity of the variance-divergence condition
in Theorem~\ref{thm:clt}.
Indeed, suppose that $\tau_n>0$ eventually and that
the standardized intrinsic-volume random variables
converge in distribution to a standard Gaussian
random variable.
If $\tau_n^2$ did not tend to infinity, there would
exist a subsequence along which the variances
remain bounded.
By Proposition~\ref{prop:bounded-variance}, this
subsequence would contain a further subsequence
converging to a discrete, nondegenerate distribution,
contradicting the assumed Gaussian limit.
Together with Theorem~\ref{thm:clt}, this proves
the equivalence \eqref{eq:clt-equivalence}.

\appendix

\section{Continuity lemma for Laplace transforms}

We record the precise continuity principle used at this point.

\begin{lemma}[Local Laplace convergence]\label{lem:laplace-continuity}
Suppose that $Y_n$ are real random variables and that, for some $r>0$,
\[
 \E\eee^{tY_n}\longrightarrow\eee^{t^2/2}
 \qquad\text{for every }|t|<r.
\]
Then $Y_n\toD\normal(0,1)$.
\end{lemma}

\begin{proof}
Fix $0<b<r$. Discarding finitely many indices if necessary, convergence
of the transforms at $b$ and $-b$ gives
\[
 \sup_n\E\eee^{b|Y_n|}
 \leq\sup_n\bigl(\E\eee^{bY_n}+\E\eee^{-bY_n}\bigr)<\infty.
\]
Markov's inequality therefore gives a uniform exponential tail bound,
so the sequence of laws is tight.

Take a subsequence converging weakly to a law $Y$. For $|t|<b/2$,
\[
 \sup_n\E(\eee^{tY_n})^2
 \leq\sup_n\E\eee^{b|Y_n|}<\infty.
\]
Thus $\eee^{tY_n}$ is uniformly integrable. Truncation of the continuous
function $y\mapsto\eee^{ty}$ allows passage of its expectation through
weak convergence, giving
\[
 \E\eee^{tY}=\eee^{t^2/2}\qquad(|t|<b/2).
\]
A probability law whose moment-generating function is finite near zero
is uniquely determined by that function. Hence $Y$ is standard normal.
Every weakly convergent subsequence has the same limit; tightness then
implies convergence of the whole sequence.
\end{proof}

\section*{Statement on the use of artificial intelligence}

The central limit theorem for conic intrinsic volumes
(Theorem~\ref{thm:clt}) was originally proved by ChatGPT
(GPT-6).
The proof presented in this paper is the result of a
subsequent extended mathematical dialogue between the
human author and ChatGPT.
The human author devoted considerable effort to understanding
and critically examining the original proof, clarifying its
arguments, separating its 
ingredients, identifying implicit assumptions, eliminating
ambiguities, simplifying individual arguments, and reorganizing
the proof into a coherent and transparent mathematical exposition.

This dialogue led to a substantial development of the mathematical
content of the paper beyond the original central limit theorem.
Among the results obtained in the course of these discussions
are the systematic analysis of the first two nonconstant
Wiener-chaos components of the squared Gaussian projection norm,
the geometric interpretations of the vector $m_C$ and the
matrix $A_C$, the characterization of optimal linear and affine
approximations to the metric projection, and the introduction
and study of the $L^2$-$\Pi$-metric $d_{\Pi,2}$.
Further developments included universal comparisons between
the intrinsic-volume variance and geometric quantities
associated with the cone, as well as geometric characterizations
of asymptotic normality.

ChatGPT played a central role throughout the mathematical
development and writing of the paper.
Its contributions included proposing and refining mathematical
ideas, suggesting proof strategies and generalizations, carrying
out technical calculations, deriving and checking mathematical
identities and inequalities, constructing examples, identifying
connections with related areas of probability and convex geometry,
researching relevant literature, and assisting with the
organization and writing of the manuscript.

The work developed iteratively.
The human author formulated research questions, supplied
mathematical direction and additional ideas, questioned proposed
arguments, suggested alternative approaches, and made decisions
concerning the scope, structure, and exposition of the paper.
Suggestions generated by ChatGPT were critically examined,
corrected, refined, or omitted from the final manuscript.
The organization and exposition of the paper were developed
through repeated discussions and revisions.

Mathematical arguments and claims included in the final
manuscript were critically examined and independently verified
by the human author, who takes full responsibility for the
contents of the paper.

\bibliographystyle{amsplain}
\bibliography{references}

\end{document}